\RequirePackage{fix-cm}
\documentclass[11pt,a4paper]{article}
\usepackage[margin=29mm]{geometry}
\usepackage{amsmath,amssymb,amsthm,mathtools,xcolor}
\usepackage{tikz}
\usepackage[hidelinks]{hyperref}
\definecolor{navy}{HTML}{173E4C}
\allowdisplaybreaks[1]
\newtheorem{theorem}{Theorem}[section]

\newtheorem{lemma}[theorem]{Lemma}
\newtheorem{corollary}[theorem]{Corollary}
\theoremstyle{definition}\newtheorem{definition}{Definition}[subsection]
\theoremstyle{remark}\newtheorem{remark}[theorem]{Remark}
\newcommand{\Div}{\operatorname{div}}\newcommand{\R}{\mathbb R}
\newcommand{\dd}{\mathrm d}
\newcommand{\avg}[1]{\langle #1\rangle}
\newcommand{\manuscripttitle}{Global classical solutions for ternary Darcy--Maxwell--Stefan diffusion on impermeable domains}
\hypersetup{pdftitle={\manuscripttitle},pdfsubject={Quasilinear parabolic systems},pdfkeywords={Maxwell--Stefan diffusion, Darcy flow, global classical solutions, Neumann boundary conditions}}
\title{\manuscripttitle}
\author{Dieter Bothe\\[3pt]
\small Department of Mathematics, Technische Universit\"at Darmstadt, Darmstadt, Germany\\
\small \href{mailto:bothe@mma.tu-darmstadt.de}{bothe@mma.tu-darmstadt.de}}
\hypersetup{pdfauthor={Dieter Bothe}}
\date{28 November 2026}
\begin{document}
\maketitle
\begin{abstract}
We prove existence and uniqueness of global classical solutions for ternary Darcy--Maxwell--Stefan diffusion on bounded \(C^{2,\alpha}\) domains with impermeable walls in dimensions two and three. On a fixed positive initial total-density band, small mean oscillation of the partial densities or uniform smallness of two minority species suffices for positive compatible data. The thresholds are independent of initial gradients and species minima. A binary characteristic region permits finite-amplitude variation of two abundant species when the third is dilute. For pairwise distinct inverse masses and pairwise distinct effective frictions we classify the invariant curves of positive states on which the system reduces to the heat equation. In two open coefficient regimes, every fixed smooth positive compatible profile gives a global solution after a sufficiently small common scaling of its densities, preserving its composition field. Further low-density classes arise from small Gibbs mean oscillation, confinement to a pressure strip, and perturbations of aligned masses. A shifted entropy and Neumann reflection reduce possible gradient blow-up to rigidity of ancient solutions. When at least two inverse masses differ, an integrating factor identifies the composition-driven pressure flux that must concentrate at a terminal gradient singularity on a positive density band. Uniform maximal \(L^p\) estimates connect fast compositional relaxation to slow scalar diffusion of total density.
\end{abstract}



\section{Introduction and main results}
\label{sec:introduction}
In an isothermal ideal mixture, the pressure is proportional to the total molar density, whereas the Darcy velocity transports total mass. Species of different molar masses can therefore diffuse relative to one another at zero net mass flux and still change the local pressure. Even an initially uniform pressure can develop a gradient and set the mixture in motion. We ask when this feedback between diffusion and bulk transport preserves classical regularity for all time.

Two features of the initial state lead to different global results. At finite total density, the relevant quantities are the mean oscillation of the partial densities and the abundance of the minority species. At low total density, a separation of time scales allows us to prescribe any fixed smooth positive compatible initial profile in two open regimes of the mass and friction parameters.

Let \(\Omega\subset\R^d\), \(d\in\{2,3\}\), be bounded and connected, and assume \(\partial\Omega\in C^{2,\alpha}\) for \(0<\alpha<1\).

We consider the planar and three-dimensional porous-medium models. The signed-source tests in Sections~\ref{sec:full-pressure} and~\ref{sec:frozen-affine-cylinder} use the logarithmic capacity estimate specific to the plane.

\subsection{Equations and solution class}
\label{sec:model}
Let \(r_i(x,t)\) be the partial mass densities and \(m_i>0\) the inverse molar masses. At each \((x,t)\) set
\[
 \rho=\sum_{i=1}^3r_i,\qquad C=m\cdot r,\qquad
 y=r/\rho,\qquad z_i=m_i r_i/C,\qquad
 E=\mathbf1^\perp,\qquad M_d=\operatorname{diag}(m_i).
\]
We call variations in the mass fractions $y$ at fixed total density compositional variations. The dimensional pressure is $RT\,C$ at fixed temperature. The constant factor is included in the Darcy mobility $\lambda>0$. The symmetric frictions $f_{ij}=f_{ji}>0$ are constant. The effective frictions are $g_{ij}=f_{ij}/(m_i m_j)$, related to the conventional diffusivities by \eqref{eq:physicaldiff}. Coincidences are allowed unless a result imposes an ordering.

With $P_y=I-y\mathbf1^\top$, define the resistance for a flux vector $J$ by
\[
 (B_y^fJ)_i=\sum_{j\ne i} f_{ij}(y_jJ_i-y_iJ_j),\qquad
 L_y=(B_y^f|_E)^{-1}.
\]
Fix constants \(\lambda,\nu>0\). The mass balances and constrained flux law are
\begin{equation}\label{eq:model}
 \partial_t r=\Div(\lambda r\nabla C+\nu J),\qquad
 B_y^fJ=M_d\nabla r-y\nabla C,\qquad
 J\in E \quad\text{in each spatial direction}.
\end{equation}
In each spatial direction the physical relative mass flux is $-\nu J$, and the barycentric velocity is $-\lambda\nabla C$. In particular $\rho_t=\Div(\lambda\rho\nabla C)$. Let \(n(x)\) be the outward unit normal. The wall condition is
\begin{equation}\label{eq:noflux}
 (\lambda r_i\nabla C+\nu J_i)(x,t)\cdot n(x)=0
 \qquad x\in\partial\Omega,\quad t>0,\quad i=1,2,3.
\end{equation}
When $\rho>0$, summation gives $\partial_n C=0$, hence $J_i\cdot n=0$ for all $i$. The normal resistance equation then gives $\partial_n r_i=0$. Conversely, componentwise zero normal derivative implies \eqref{eq:noflux}. Integration of each species balance then gives conservation of each spatial mean \(\avg{r_i}=|\Omega|^{-1}\int_\Omega r_i\,\dd x\).

\begin{definition}\label{def:classical}
Fix $p>d+2$. Write \(X_\gamma=W^{2-2/p,p}_{\rm N}(\Omega;\R^3)\), where the subscript denotes the homogeneous Neumann trace. A global classical solution satisfies, for every finite \(T>0\),
\[
 r\in W^{1,p}(0,T;L^p(\Omega;\R^3))
 \cap L^p(0,T;W^{2,p}_{\rm N}(\Omega;\R^3))
 \cap C([0,T];X_\gamma).
\]
It also belongs to \(C^1((0,\infty);C(\bar\Omega;\R^3))
\cap C((0,\infty);C^2(\bar\Omega;\R^3))\), and its partial densities are
strictly positive on \(\bar\Omega\). The initial trace \(r^0\in X_\gamma\)
is strictly positive on \(\bar\Omega\).
The results with prescribed initial total density state any additional
regularity required of the initial data.
\end{definition}

\subsection{Global existence and continuation}
For \(x\in\bar\Omega\) and \(R>0\) set \(D_R(x)=\Omega\cap B_R(x)\), with \(B_R(x)\) the Euclidean ball. Write \((v)_D=|D|^{-1}\int_D v\,\dd x\). The intrinsic mean-oscillation seminorm is
\[
 [v]_{\mathrm{BMO}_2(\Omega)}
 =\sup_{\substack{x\in\bar\Omega\\R>0}}
 \left(\frac{1}{|D_R(x)|}\int_{D_R(x)}
        |v-(v)_{D_R(x)}|^2\,\dd x'\right)^{1/2}.
\]

\begin{theorem}[Principal global existence results]\label{thm:main}
Fix $p>d+2$ and positive inverse masses and frictions as in Section~\ref{sec:model}.
\begin{enumerate}
\item Fix \(0<L_0<U_0<\infty\). For positive compatible data \(r^0\in X_\gamma\) with \(L_0\le\sum_i r_i^0\le U_0\), there is a unique global classical solution if either \([r^0]_{\mathrm{BMO}_2(\Omega)}\) or \(\|\sum_{j\ne i}r_j^0\|_\infty\) for some fixed species \(i\) is sufficiently small. The thresholds are independent of the initial gradients and individual species minima.
\item Consider parameters in a sufficiently small open neighborhood of a reference point from either family below.
\begin{align*}
 &f_{ij}^*=a\,m_i^*m_j^*,\quad a>0,\quad
                 m_1^*,m_2^*,m_3^*>0;\\
 &m_2^*=m_3^*=m_b,\quad f_{12}^*=f_{13}^*=a,\quad f_{23}^*=b,
   \quad m_1^*,m_b,a,b>0,\quad m_1^*\ne m_b,\quad a\ne b .
\end{align*}
For every fixed smooth positive compatible profile \(r^0\), the system with initial data \(\delta r^0\) has a unique global classical solution for all sufficiently small \(\delta>0\). The threshold can depend on \(r^0\). The second family of open sets allows arbitrarily large mass and friction contrasts independently.
\end{enumerate}
\end{theorem}

\begin{proof}[Proof of Theorem~\ref{thm:main}]
The two alternatives in part~1 follow from Theorem~\ref{thm:bmo-global} and part~1 of Theorem~\ref{thm:finite-density}, respectively. Theorem~\ref{thm:all-shapes} proves part~2 by establishing rigidity at the reference coefficients, passing to open neighborhoods, and applying Lemma~\ref{lem:comparison}.
\end{proof}

The first assertion controls initial variation on all spatial scales. In particular, the mean-oscillation condition allows a fixed pointwise excursion and arbitrarily large initial gradients (Section~\ref{sec:finite-bmo}). There is also a finite-density regime with two abundant species of distinct inverse masses. Set \(\rho_{12}=r_1+r_2\) and \(C_{12}=m_1r_1+m_2r_2\). Two characteristic curves of the binary system and the lower edge of its physical state space bound a convex region in the \((\rho_{12},C_{12})\)-plane. For two such nested regions, initial values in the inner region and a sufficiently dilute third species ensure confinement to the outer region and global classical existence. Figure~\ref{fig:physical-lens} shows the range of compositions this permits (Theorem~\ref{thm:finite-density}).

In the second assertion, multiplying all partial densities by \(\delta\) preserves the mass-fraction field \(y^0=r^0/\rho^0\) pointwise. Thus low density imposes no smallness condition on the composition oscillation of the chosen profile. Corollary~\ref{cor:constant-d} gives the corresponding result for constant conventional binary diffusivities in the stated coefficient regimes.

Three additional classes are described in Section~\ref{sec:low-density}. Small mean oscillation of the Gibbs composition suffices at any fixed smooth positive total-density profile (Theorem~\ref{thm:gibbs-bmo}). For the effective-friction orders of Theorem~\ref{thm:frozen-pressure-strip}, small initial deviation of the specific pressure $m\cdot y^0$ from a selected level and a nearly constant prescribed total suffice. Near aligned inverse masses and frictions, Theorem~\ref{thm:aligned} allows perturbations of any fixed smooth positive constant-total reference shape, including shapes with large composition gradients.

A third result classifies exact reductions to a scalar equation. Assume the inverse masses are pairwise distinct, and order the species by their opposite effective frictions \(\ell_i=g_{jk}\), \(\{i,j,k\}=\{1,2,3\}\), so that \(\ell_1<\ell_2<\ell_3\). A smooth regular curve of strictly positive states is called invariant if the evolution preserves it locally for every smooth scalar profile on the curve. Such a curve exists precisely when the molar mass \(1/m_2\) is a strict extremum among the three molar masses. Every maximal curve is then an affine segment of constant \(\rho\) and \(C\), along which composition obeys a scalar Neumann heat equation. Profiles with range compactly inside the segment have open neighborhoods of global data in the strong trace space (Theorem~\ref{thm:pressure-heat-line}).

\subsection{Structure of the proofs}

\paragraph{Ancient solutions and impermeable walls.}
Suppose a spatial gradient becomes unbounded while the total density stays between positive bounds. Rescaling at a succession of gradient maxima gives a limiting solution on all of \(\R^d\) and for all past times, with a unit gradient at the chosen point. At a wall we first flatten the boundary and reflect the no-flux equation. A shifted entropy controls squared gradients even when a species vanishes in the limit, and local parabolic estimates give convergence of gradients up to the terminal time of the rescaling (Section~\ref{sec:boundary-estimates}). The remaining step depends on the mixture. For example, in the binary region the curvature of the characteristic curves yields scalar entropies whose dissipation forces every bounded-gradient ancient solution confined to that region to be constant. This rules out the unit gradient at the rescaled point (Lemma~\ref{lem:binary-ancient}).

\paragraph{A potential for the pressure flux.}
In molar fractions \(z\), the pressure balance admits the decomposition
\[
 C_t=\Div\big((\lambda C+\nu k(z))\nabla C+F\big),
 \qquad F=\nu C\,\alpha(z)\nabla\Phi(z).
\]
Here \(k\ge0\). When at least two inverse masses differ, an integrating factor gives a scalar potential \(\Phi\) with \(\alpha\) of fixed sign and bounded away from zero in absolute value (Section~\ref{sec:full-pressure}). This formula makes the composition-driven flux \(F\) a test for regularity. If \(\rho\) stays in a compact subinterval of \((0,\infty)\) near time \(T\), then at every terminal gradient singularity \((x_0,T)\) one must have
\[
 \int_{T-R^2}^T\!\int_{\Omega\cap B_R(x_0)} |F|^2\,\dd x\,\dd t
 \;\ge\;\epsilon_0 R^d
\]
for all sufficiently small \(R\) (Theorem~\ref{thm:full-pressure-continuation}). The frozen system has a related rigidity property. When the prescribed total density is constant and at least two inverse masses differ, a bounded ancient whole-space solution in the local energy class whose pressure is independent of time has constant composition. This yields local gradient tests based on pressure variation. A separate test uses the smallest eigenvalue of the composition covariance on its two-dimensional tangent plane (Theorem~\ref{thm:frozen-affine-cylinder}).

\paragraph{Fast composition and slow total density.}
After the common density scaling, the fast reference evolution holds the initial total density fixed while the species diffuse. In Gibbs coordinates, vanishing microscopic flux corresponds to spatially constant differences of chemical potentials. The conserved species means select the Gibbs profile for each prescribed total-density field. On the slow time scale this profile determines the conductivity in a scalar Neumann diffusion equation for the total density. Uniform maximal \(L^p\) estimates compare the full system with the two reference evolutions on successive finite intervals. Conservation also recovers the spatially constant mode of the chemical-potential differences, which is lost when estimating their mean-zero parts. Once the slow reference is near homogeneous equilibrium, a stability radius independent of the small time-scale ratio gives global continuation (Lemma~\ref{lem:comparison}).

For arbitrary positive coefficients, Theorem~\ref{thm:all-shapes} gives the same low-density conclusion for every fixed initial shape provided that all bounded-gradient ancient solutions of the constant-total frozen system on the closed composition simplex are constant. With $h$ defined by $f_{ij}=h_i+h_j$, Theorem~\ref{thm:interior-rigidity} proves rigidity on every compact interior state range when $m,h,\mathbf1$ are linearly independent. Under this independence assumption, any nonconstant bounded-gradient ancient solution must satisfy $\inf_{x,\tau}\min_i y_i(x,\tau)=0$.

\subsection{Relation to Darcy-scale and Maxwell--Stefan models}

In \eqref{eq:model} the pressure gradient determines the common velocity, and the projected species forces determine the relative flux.

\paragraph{Darcy-scale constitutive laws.}
Brust, Hopf, Fuhrmann, Cheilytko, Wullenkord and Sattler \cite{brust} give a Darcy-scale model for ideal gases in a porous medium. Their isothermal equations use an additive sum of convective and diffusive mass fluxes, the ideal gas law, and the full Maxwell--Stefan driving force including pressure diffusion. In the notation here their molar and mass fractions are $z_i$ and $y_i$. The identity
\[
 \nabla z_i+(z_i-y_i)\nabla\log C
       = C^{-1}(m_i\nabla r_i-y_i\nabla C)
\]
identifies their force with ours. Their equations (2.22), (2.24)--(2.26) in the WIAS version give this comparison directly. With constant porosity absorbed in time, scalar constant permeability and viscosity, and no reactions, their isotropic relation gives \eqref{eq:model} upon taking the usual molar binary diffusivity in the form
\begin{equation}\label{eq:physicaldiff}
 D_{ij}^{\rm eff}(r)=\frac{\nu\rho}{C}\,\frac{m_im_j}{f_{ij}} .
\end{equation}
Indeed insert $j_i=-\nu J_i$ for the relative mass flux into
$-\sum_{j\ne i}z_i z_j(D_{ij}^{\rm eff})^{-1}(j_i/r_i-j_j/r_j)$.
This yields the displayed projected force using the normalization in \eqref{eq:model}. In dimensional Darcy's law the scalar factor is $\lambda=KRT/\eta$ before any porosity rescaling. Thus constant $f_{ij}$ in the present resistance convention corresponds to the common state factor $\rho/C$ in \eqref{eq:physicaldiff}. Brust et al.\ establish thermodynamic consistency and steady-state properties and study the model numerically with a finite-volume method.

The low-density limit in this paper concerns the specified continuum constitutive law. At fixed composition, \eqref{eq:physicaldiff} is independent of the overall density scale. For comparison, the isothermal Chapman--Enskog dilute-gas law $D_{ij}\propto C^{-1}$ corresponds to $f_{ij}=\rho\widehat f_{ij}$ with constant $\widehat f_{ij}$. With physical densities $\delta r$ and fast time $\tau=\nu t_{\rm phys}/\delta$, the normalized microscopic operator then acquires the factor $1/\rho$ inside its flux and the Darcy term the factor $\delta^2/\nu$. Extending the comparison estimates to this density-dependent resistance is a separate problem. At low pressure, collisions with pore walls can also contribute through Knudsen diffusion \cite{mason}.

Earlier porous-gas models of Stockie, Promislow and Wetton \cite{stockie} and Weber and Newman \cite{weber} also combine pressure-driven Darcy convection with multicomponent gas diffusion. The Maxwell--Stefan closure in the appendix of \cite{stockie}, for example, uses mole-fraction gradients. Its computational boundary setting and that driving-force choice differ from the present projected force and impermeable walls. Binning et al.\ \cite{binning} give a related reactive porous-column model with gravity and gas--solid friction. The dusty-gas formulation \cite{mason} describes gas--solid interaction by treating the stationary porous matrix as an additional species. Here the convex free energy and the zero-sum microscopic flux are for the three mobile species.

\paragraph{Rigorous results for related evolution equations.}
Druet and J\"ungel \cite{druetjungel} study the same common Darcy velocity on a bounded impermeable domain. Their species flux is pure transport by this velocity. Under their equation-of-state assumptions, a scalar volume variable satisfies a uniformly parabolic equation along the positive branch and the remaining fractions solve transport equations. This yields global classical and weak results for that parabolic--hyperbolic system. The ideal linear law $C=m\cdot r$ belongs to their equation-of-state class. Adding a Maxwell--Stefan flux changes the scalar reduction, since the molar-density balance then contains a nonzero microscopic flux.

J\"ungel, Miky\v{s}ka and Zamponi \cite{jmz} study Darcy transport with a symmetric Onsager diffusion matrix and van der Waals mixture pressure. They prove global weak existence and study the approach to equilibrium, with equilibrium boundary conditions in the analytic problem. For comparison, write \eqref{eq:model} as
\[
 r_t=\Div\!\left[(\lambda rr^\top+
  \nu L_y(\operatorname{diag}(r)-rr^\top/\rho))\nabla\mu\right],
       \qquad \mu_i=m_i\log r_i .
\]
Their estimates use a uniform transverse lower bound for the chemical-potential diffusion matrix. The Maxwell--Stefan mobility in this display loses that bound at a missing-species face, although the physical density matrix remains normally elliptic on finite positive total-density bands.

Ostrowski and Rohde \cite{ostrowski} retain a momentum equation for each species and include gas--solid drag. Their Darcy-scale system consists of hyperbolic balance laws with relaxation. They prove whole-space global classical solutions for small perturbations of equilibrium and a relative-entropy comparison with a given smooth solution of a parabolic relaxation limit. Even with common matrix drag, the limiting force system contains that drag in the resistance for relative as well as total motion.

The projected chemical-potential force follows the continuum thermodynamic formulation of Bothe and Dreyer \cite{bothedreyer}. Druet \cite{druetideal} uses the ideal-gas Maxwell--Stefan force, including the pressure projection, in a global generalized weak theory with compressible Navier--Stokes momentum, additional high-pressure constitutive growth and low-pressure transport assumptions. Bothe and Druet \cite{bothedruet} prove local strong existence for compressible class-one models with a Navier--Stokes barycentric velocity and near-equilibrium existence on a prescribed finite horizon. For isobaric Maxwell--Stefan diffusion, entropy methods give global weak existence \cite{jungelstelzer}, and local strong Neumann theory is available in \cite{hmpw}. The ternary isobaric preprint \cite{ternary} uses two Stieltjes roots to prove global strong existence at constant total and equal molar masses. Bothe's additive-friction preprint \cite{botheadditive} develops the Stieltjes coordinates and mixed-production cancellation for isobaric diffusion. We use the separated-root estimate of \cite{ternary} in one frozen reference class.

\section{Quasilinear structure and local theory}
\label{sec:local-theory}
We denote the full diffusion matrix in physical density variables by
\begin{equation}
 A_{\rm full}(r)=\lambda r\,m^\top+\nu L_yP_yM_d .
\end{equation}
On a compact band $0<L\le\rho\le U$ this matrix and its derivatives extend smoothly through each individual species face. For each positive \(y\), the resistance is positive on \(E\) in the \(y^{-1}\) metric, since
\[
 \langle v,B_y^f v\rangle_{y^{-1}}
   =\sum_{i<j} f_{ij}y_i y_j
                \left|\frac{v_i}{y_i}-\frac{v_j}{y_j}\right|^2,\qquad
 \langle v,w\rangle_{y^{-1}}=\sum_i v_iw_i/y_i .
\]
The same expression with all $f_{ij}=1$ equals $\langle v,v\rangle_{y^{-1}}$ for $v\in E$. At an inactive row $i$ the resistance has only its positive diagonal $\sum_{j\ne i}f_{ij}y_j$, and reduces on the active tangent block to its resistance there. Thus it is invertible on $E$ through the faces. At the same state the $i$th row of $A_{\rm full}$ has only the positive diagonal $\nu m_i/\sum_{j\ne i}f_{ij}y_j$. Smoothness gives an $O(r_i)$ factor in its off-diagonal entries nearby.

The free energy in the chosen units is $H(r)=\sum_i m_i r_i(\log r_i-1)$. Its derivative $\mu_i=m_i\log r_i$ gives $\sum_i r_i\nabla\mu_i=\nabla C$. Set $d_r=M_d\nabla r-y\nabla C$. The entropy production at positive states is
\begin{equation}\label{eq:entropy}
 D(r,\nabla r)=\lambda|\nabla C|^2+
       \frac{\nu}{\rho}\langle d_r,L_y d_r\rangle_{y^{-1}} .
\end{equation}
The factor $\lambda|\nabla C|^2$ is the dissipation of the common Darcy velocity. The second term follows from $J\in E$ and
$\nabla\mu_i=(d_r)_i/r_i+\nabla C/\rho$. In the interior $D^2H\,A_{\rm full}$ is symmetric positive definite. For example the entropy mobility is the sum of $\lambda r r^\top$ and $\nu L_y(\operatorname{diag}(r)-rr^\top/\rho)$, the latter symmetric positive on $E$ by the resistance metric. At a face the active block has the same positive spectrum as the smaller mixture, and the missing rows supply the positive diagonals. This proves uniform normal ellipticity on any closed finite positive total-density band. The ellipticity constants and the coefficient bounds can be chosen uniformly on compact families of positive parameters, including families with coincident masses or frictions. On the band \eqref{eq:entropy} controls $\sum_i|\nabla r_i|^2/r_i$ by decomposing $M_d\nabla r=d_r+y\nabla C$ orthogonally in the $y^{-1}$ metric. The identity has zero entropy flux on the wall and its spatial integral is the entropy dissipation. The local entropy flux is bounded by $C D^{1/2}$, uniformly through a species face, from the bounds on $\sum_i r_i|\log r_i|^2$ and the resistance estimate.

To normalize a low-density shape, write $r^{\rm phys}=\delta r$, $\tau=\nu t_{\rm phys}$, and $\varepsilon=\delta/\nu$. In normalized variables
\[
 r_\tau=\Div\big(L_yP_y M_d\nabla r+\varepsilon\lambda r\nabla C\big).
\]
For $\varepsilon=0$ total density is a prescribed function $\rho^0(x)$, and
\begin{equation}\label{eq:frozen}
 \rho^0 y_\tau=\Div\big(\rho^0 A(y)\nabla y+K(y)\nabla\rho^0\big),\qquad
 A(y)=L_yP_yM_d|_E,\qquad K(y)=L_y(M_dy-y\,m\cdot y).
\end{equation}
Both factors of $A$ are positive self-adjoint on $E$ in the $y^{-1}$ metric at interior states, and the active and missing blocks again give normal ellipticity up to the faces. Since $\partial_n\rho^0=0$, zero microscopic normal flux implies $\partial_n r=0$. In fact, applying $M_d^{-1}$ to the normal force relation and summing gives $(\sum y_i/m_i)\partial_n C=0$. Thus \eqref{eq:frozen} has $\partial_n y=0$. Its entropy dissipation consists of the microscopic term in \eqref{eq:entropy} without $\nu$, and controls the square of the force uniformly on the prescribed density band.

\paragraph{Local existence, regularity and positivity.}
Normally elliptic second-order systems of the form here have the common homogeneous Neumann domain $W^{2,p}_{\rm N}(\Omega)$. To check the complementing condition at a frozen point, the decaying half-line solutions of the resolvent problem are
$e^{-x_d(|\xi'|^2 I+\zeta A_0^{-1})^{1/2}}v$, for a normally elliptic matrix $A_0$, $\mathrm{Re}\,\zeta\ge0$ and $|\zeta|+|\xi'|>0$. The stable square root is invertible, so the homogeneous normal condition annihilates their amplitude. An orthogonal coordinate change or a flat chart normalization preserves this argument. Apply \cite[Theorem 2.1]{dhp} with equal time and space integrability exponent \(p\). The choice \(p>d+2\) gives \(X_\gamma\hookrightarrow C^{1,\beta}(\bar\Omega)\) for \(0<\beta<1-(d+2)/p\). Smooth state dependence then gives local existence, uniqueness and continuous dependence in \(X_\gamma=W_{\rm N}^{2-2/p,p}\), as in \cite[Theorem 3.2]{hmpw}. The same applies to the reduced frozen system. Classical boundary regularity follows from Lemma~\ref{lem:classical-upgrade}. Compactness of these Sobolev traces is used at restart.

A species is inactive at a state where its density is zero. On a neighborhood of the closed density band use the same algebraic formula for the flux, with $L_y=(B_y^f|_E)^{-1}$. The inverse persists there by compactness. Even for nearby signed states with $r_i=0$, the $i$th resistance row is diagonal, so the off-diagonal entries of the $i$th diffusion row are smoothly divisible by $r_i$. The diagonal is positive near that face. Patching it to a positive function away from the face and putting the difference into the drift gives the scalar row
\[
 (r_i)_t=\Div(a_i\nabla r_i+r_i V_i),\qquad a_i\ge c>0 .
\]
The same construction applies to the frozen flux in $r$ variables at prescribed total. On a strong strip $V_i$ is bounded by the spatial gradient bound. Testing with the negative part, with its negative sign, and using zero normal flux gives
\[
 \tfrac12\tfrac{\dd}{\dd t}\|r_i^-\|_2^2+
 \tfrac c2\|\nabla r_i^-\|_2^2
 \le C\|V_i\|_\infty^2\|r_i^-\|_2^2 .
\]
The diagonal adjustment and compensating drift preserve the row flux pointwise and hence its normal trace. Thus a local restart from a nonnegative trace preserves nonnegativity. The strong trace gives positivity for a short time from positive data. Subsequently the scalar divergence-form Harnack inequality applies to nonnegative components on buffered strips, since $a_i$ is uniformly positive and bounded and $V_i$ is bounded. At a wall use the variational reflection below, reflecting the normal drift component oddly. A zero at positive time then forces the component to vanish on its connected past, in contradiction with its positive conserved integral. The same backward implication applies to a bounded-gradient entire limit. These arguments use no uniform lower bound on individual species in a limiting density band.

\subsection{Classical derivatives at an impermeable wall}

\begin{lemma}[Classical regularity for positive times]\label{lem:classical-upgrade}
Let \(\partial\Omega\in C^{2,\alpha}\) for \(0<\alpha<1\). Let \(r\) be a strong solution of the physical system on a finite strip with its state range in a compact nonnegative set of positive total density. For each \(0<\beta<\min\{\alpha,1-(d+2)/p\}\), on every compact positive-time subinterval one has
\[
 r\in C_t^1(C(\overline\Omega))\cap C_t(C^{2,\beta}(\overline\Omega)).
\]
In particular the physical equations and the normal boundary condition hold pointwise up to the wall at positive times.
\end{lemma}

\begin{proof}
In nondivergence form, on a neighborhood of this state range where the principal matrix remains normally elliptic,
\[
 r_t=A_{\rm full}(r)\Delta r
      +\sum_k \big(DA_{\rm full}(r)[\partial_k r]\big)\partial_k r,
 \qquad \partial_n r=0 .
\]

Write \(A=A_{\rm full}\) and \(N(r)=\sum_k(DA(r)[\partial_k r])\partial_k r\).
Fix a positive time \(t_*\) and a short interval \(I=[t_0,t_1]\) with \(t_0<t_*<t_1\).
Use the strong space \(W^{1,p}(I;X_0)\cap L^p(I;X_1)\) with fixed initial trace \(r(t_0)\), where \(X_0=L^p(\Omega;\R^3)\) and \(X_1=W^{2,p}_{\rm N}(\Omega;\R^3)\). The residual for the equation with right-hand side multiplied by \(1+\theta\) is continuously differentiable on this affine space near \(r\), by the embedding of \(X_\gamma\) into \(C^1\) and smoothness of \(A\). Its derivative in the solution at \(\theta=0\), on zero-initial-trace variations \(v\), is
\[
 v_t-A(r)\Delta v-DA(r)[v]\Delta r-DN(r)[v].
\]
The frozen operator \(\partial_t-A(r(t_*))\Delta\) has a bounded zero-initial-trace inverse. The maximal-regularity and \(C(I;X_\gamma)\) trace bounds for this inverse are uniform as \(I\) shrinks, by translating the lower face, extending the forcing by zero to a fixed time interval and restricting its solution. In particular,
\begin{align*}
 \|DA(r)[v]\Delta r\|_{L^p(I;X_0)}
   &\le C\|\Delta r\|_{L^p(I;X_0)}\|v\|_{C(I;X_\gamma)},\\
 \|DN(r)[v]\|_{L^p(I;X_0)}
   &\le C|I|^{1/p}\|v\|_{C(I;X_\gamma)} .
\end{align*}
The latter bound follows by differentiating \(N\), using the uniform \(C^1\) bounds for \(r\) and the \(C^1\) trace embedding for \(v\). The remaining principal difference has operator norm bounded by \(C\sup_I\|r(t)-r(t_*)\|_{C(\bar\Omega)}\). This supremum and both perturbation factors tend to zero as \(I\) shrinks about \(t_*\). Hence the linearization is invertible by perturbation, and the implicit function theorem gives differentiable dependence on \(\theta\) in the strong space and its continuous trace.

On a common shorter interval still containing \(t_*\), uniqueness identifies the perturbed solution with \(r(t_0+(1+\theta)(t-t_0))\). The difference quotient converges in \(L^p(I';X_0)\) on this shorter interval \(I'\) to \((t-t_0)r_t(t)\), by the \(W^{1,p}\) time chain rule and continuity of time dilations in \(L^p\). The implicit-function derivative gives this expression a continuous \(X_\gamma\) representative. Dividing by \(t-t_0>0\) near \(t_*\) shows that \(r_t\) is continuous in \(X_\gamma\) there. Since \(t_*\) was arbitrary, this holds on every compact positive-time subinterval.

The matrix \(A_{\rm full}(r)\) is invertible and smooth there. Hence
\[
 \Delta r=A_{\rm full}(r)^{-1}
     \left(r_t-\sum_k (DA_{\rm full}(r)[\partial_k r])\partial_k r\right)
\]
holds for almost every time and has a right-hand side continuous in time with values in \(C^\beta(\overline\Omega)\). At times when the equation holds, apply the scalar elliptic Neumann Schauder estimate to each component with its \(C^0\) norm controlling the additive constants. Applied to differences at such times, the estimate makes \(r\) Cauchy in \(C^{2,\beta}(\overline\Omega)\) along any sequence converging to a fixed time. By continuity of the strong trace this limit is \(r\) at that time. The same difference estimate gives continuity in \(C^{2,\beta}(\overline\Omega)\). This is the classical positive-time upgrade of the Sobolev trace, as in \cite[Theorem 3.3]{hmpw}. Only the original first-order trace is required of the initial data. The same argument applies to the frozen equation when $\rho^0$ is prescribed in $C^{2,\alpha}$, since its additional flux is a fixed smooth function of the state and the first two spatial derivatives of $\rho^0$.
\end{proof}

\begin{remark}
On compact \(C^2\) boundaries the global existence results hold in the Sobolev strong class of Definition~\ref{def:classical}. The positive-time \(C^2(\bar\Omega)\) conclusion uses the \(C^{2,\alpha}\) wall regularity of Lemma~\ref{lem:classical-upgrade}.

The compact $C^2$ wall gives a $C^{1,1}$ reflected chart and a $C^2$ distance collar, as used in Lemmas~\ref{lem:wall-record} and~\ref{lem:neumann-weight}. These bounds suffice for the localized strong estimates and the Neumann restart by \cite[Theorem 2.1]{dhp}. The additional wall exponent enters the elliptic Schauder step of Lemma~\ref{lem:classical-upgrade}.
\end{remark}

\section{Boundary estimates and compactness}
\label{sec:boundary-estimates}

\paragraph{Conventions.}
In estimates, \(C\) also denotes a positive constant that may change from line to line. Where both meanings occur in a single calculation, \(C_{\mathrm{est}}\) denotes the constant.

A pole is the spatial center and upper time of a backward cylinder. A physical cap is its spatial ball intersected with $\Omega$. A buffered subcylinder has positive clearance from the artificial lateral boundary and the lower time face. A gradient record is a pole selected at a running maximum, or a comparable localized maximum, of the scaled gradient so that earlier nearby gradients are controlled. A strong branch is a solution in the mixed Sobolev class of Definition~\ref{def:classical} on each compact subinterval of its existence interval. By a local-energy solution we mean a weak solution of the divergence identity in $L^\infty_{\rm loc}(L^2_{\rm loc})\cap L^2_{\rm loc}(H^1_{\rm loc})$, with the natural variational no-flux condition when a wall is present.

\begin{lemma}[Shifted entropy on the closed state range]\label{lem:shifted-entropy}
Fix a compact positive total-density band and a compact family of positive parameters. A fixed sufficiently small $\gamma>0$ gives a smooth strictly convex entropy
\[
 H_\gamma(r)=\sum_i m_i(r_i+\gamma)(\log(r_i+\gamma)-1)
\]
such that the quadratic form of $D^2 H_\gamma(r)A_{\rm full}(r)$ is uniformly coercive on the closed nonnegative band. The analogous estimate for the frozen operator holds on the tangent plane $E$. The relative entropy
\[
 H_\gamma(r\mid a)=H_\gamma(r)-H_\gamma(a)-DH_\gamma(a)\cdot(r-a)
\]
is uniformly equivalent to $|r-a|^2$ for $r,a$ in the band. For constant $a$, its unforced entropy flux is bounded by $C|r-a||\nabla r|$.
\end{lemma}

\begin{proof}
Put $D_\gamma=\operatorname{diag}\sqrt{m_i/(r_i+\gamma)}$ and consider the real quadratic form of $Q_\gamma=D_\gamma A_{\rm full}D_\gamma^{-1}$. If no positive uniform lower bound held for small $\gamma$, there would be a sequence $\gamma_j\downarrow0$, admissible states and parameters, and real unit vectors $v_j$ with $v_j^\top Q_{\gamma_j}v_j\le 1/j$. Pass to a subsequence on which the states, parameters and vectors converge. The off-diagonal entries of a row becoming inactive tend to zero by the $O(r_i)$ row bound and $r_i/\sqrt{r_i+\gamma_j}\le\sqrt{r_i}$. Entries of an inactive column in active rows tend to zero by the factor $\sqrt{r_i+\gamma_j}$. The limit is block diagonal, with a symmetric positive entropy-normalized active block and positive inactive diagonals. Its quadratic form is positive on the limiting unit vector, a contradiction.

For the frozen operator, use $L_yP_yM_d$ in the conjugation and choose the unit vectors in $D_{\gamma_j}E$. Their limits are orthogonal to $(\sqrt{r_i/m_i})_i$ at the limiting state. The active component is then in the entropy-normalized active tangent space, where the active block is positive. The inactive diagonals control the remaining directions, giving the same contradiction. Undo the conjugation using the uniform positive lower bound of $D_\gamma$ on the band. For fixed $\gamma$, the uniform Hessian bounds along segments in the band give the relative-entropy comparison. The flux estimate follows by combining the Lipschitz bound for $DH_\gamma$ with the bounded diffusion matrix. The estimates pass to invariant faces and fixed affine coordinates on their state planes.
\end{proof}

\begin{lemma}[Boundary gradient records and terminal traces]\label{lem:wall-record}
Consider the full system in a fixed closed nonnegative band $0<L\le\rho\le U$, or the frozen equation on the closed composition simplex. The parameters are fixed or range in a compact subset of their positive parameter space. In the frozen case the prescribed densities lie in a bounded $C^2$ family with a common positive lower bound and satisfy the homogeneous Neumann condition.

In interior coordinates or a common reflected chart, a uniform normalized gradient bound on a backward box gives subsequential convergence of the states and spatial gradients on each smaller box with a fixed spatial and lower-time buffer. The convergence includes the upper face. This local conclusion applies to fixed-scale boxes and to dilated boxes with the uniform normalized coefficient bounds above.

Suppose now that strong solutions are dilated at poles in \(\bar\Omega\) or in the whole space. If the spatial scales tend to zero, the dilated time intervals exhaust the past and the normalized gradients are uniformly bounded on compact backward boxes up to their upper faces, the local convergence gives a bounded entire ancient limit along a subsequence with converging parameters. At a wall this limit is obtained by reflection if the dilated wall stays at finite distance. A unit gradient at the record pole persists, and the limiting prescribed density in the frozen case is constant.

On the original bounded domain, a bounded gradient on a finite terminal strip with state range as above gives an \(X_\gamma\) terminal trace and strong restart. For the autonomous constant-density whole-space equation the interior compactness and entire-limit conclusions also apply to dilations at arbitrary spatial scales when the corresponding box and time hypotheses hold.

\end{lemma}

\begin{proof}
Write the wall locally as $\varsigma(z')$ with inward unit normal $N(z')$. In a sufficiently small chart let
\[
 X(z',s)=\varsigma(z')+sN_s(z'),\qquad
 N_s(z')=\int N(z'-s\xi)\eta(\xi)\,\dd\xi,\qquad s>0,
\]
where $\eta$ is an even smooth mollifier of compact support and integral one. The first derivatives of $N_s$ in $(z',s)$ are bounded and $s|D^2N_s|\le C\|\nabla N\|_\infty$. The chart and its inverse are $C^{1,1}$. At the wall $\partial_s X=N$, so the mixed entries $g^{ad}$ of the inverse metric vanish for $a<d$. With $\mathcal J=|\det DX|$, extend the tangential block of $g$, $g^{dd}$ and $\mathcal J$ evenly across $s=0$, and the mixed entries oddly. The resulting metric is Lipschitz and uniformly positive. Extend the state and prescribed density evenly. In the full case the reflected equation is
\[
 \mathcal J U_t=\partial_a\big(\mathcal J g^{ab}A_{\rm full}(U)\partial_b U\big).
\]
For the frozen case the time weight is $\mathcal J\rho^0$ and the flux is the pullback of \eqref{eq:frozen}. The variational reflection follows by using the sum of the two half-chart test functions and the no-flux condition. Tangential first derivatives match across the wall. The normal derivative is odd and has zero trace by the Neumann condition. Integration by parts on the two halves therefore produces no interface measure in the weak second derivatives of the even extension.

On a backward box with bounded gradients the reflected equation has the form $w(X)U_t=\Div\mathcal F$, with $w$ positive, time independent and Lipschitz, and with bounded $\mathcal F$. In the interior take the analogous time weight without $\mathcal J$. With $\eta_\ell$ a unit-mass spatial mollifier supported in a ball of radius $\ell$, test with $\eta_\ell/w$. The change in $\int U\eta_\ell$ is bounded by $C|t-s|/\ell$, and the spatial averaging error by $C\ell$. For small $|t-s|$, take $\ell=|t-s|^{1/2}$ inside a buffered spatial box. This gives a uniform time modulus, including at the upper time face.

The expanded principal tensor is $g^{ab}A(U)$, with $A_{\rm full}$ in the full case. It has a common modulus of continuity and a uniform normal-ellipticity bound by the compact parameter hypothesis. All remaining terms are bounded by the gradient control, the Lipschitz metric bounds and the $C^2$ bounds of the prescribed profile. Freeze the principal tensor on small backward boxes. After a linear spatial change, put the frozen species matrix in complex Schur triangular form. Its diagonal entries $a$ satisfy $\Re a\ge c>0$ and $|a|\le C$. For each scalar diagonal equation the heat kernel with Fourier multiplier $e^{-ta|\xi|^2}$ gives the parabolic singular-integral bounds for second derivatives, uniformly on this compact set of $a$. The equation then estimates the time derivative. Solve from the last row to the first, treating the bounded off-diagonal coefficients times previously estimated second derivatives as forcing. Only first derivatives of the reflected metric occur in the lower-order terms, so its Lipschitz bound suffices. In a cutoff problem with zero data at the lower time face the coefficient oscillation can be absorbed in both $W^{2,1}_p$ and $W^{2,1}_q$. Extend the principal perturbation $\mathcal E$ by zero outside the localization region. The localized forcing is bounded with compact spatial support on a finite time interval. Write $P_0$ for the zero-initial inverse of the constant-coefficient equation. The series $\sum_{k\ge0}(P_0\mathcal E{:}D^2)^k P_0 f$ converges in both strong spaces when $\|\mathcal E\|_\infty$ is small enough for both exponents. The constant-coefficient inverse agrees on common forcings by its heat kernel representation, so the partial sums, and hence the distributional limits, coincide. Uniqueness in the original space identifies the localized solution. This gives the estimate for any chosen finite $q>p$. The cutoffs require a spatial and lower-time buffer and may be nonzero through the top. The estimates on shorter backward boxes are uniform as their upper faces approach that top. Parabolic Sobolev embedding gives common moduli of continuity for the states and spatial gradients up to the upper face, hence compactness there. The uniform local \(W_q^{2,1}\) estimates on the rescaled boxes also pass to each limit by weak compactness.

If the pole-to-wall distance divided by the scale tends to infinity, these estimates are interior. Otherwise center the chart at a nearest wall point and choose an orthonormal tangent frame there. In the rescaled coordinates centered at the pole the wall tends to a plane at offset equal to the limit of the normalized distance. The chart derivatives tend to the corresponding orthogonal change of frame and the reflected metric tends to the identity. The first and second scaled derivatives of the prescribed profile tend to zero, and the scaled profile tends locally to a positive constant. Passing to the reflected divergence identity gives the entire limiting equation. Gradient convergence in these coordinates also preserves the physical gradient magnitude at the pole.

At fixed scales the same estimates hold on buffered terminal strips. A finite covering gives the strong estimates through their upper time face. For $q>p$, the trace in $W^{2-2/q,q}$ is uniformly bounded on the shorter terminal strip and this space embeds compactly into $X_\gamma$. The time modulus identifies a unique terminal value and implies convergence in $X_\gamma$. Local theory at that value uses the smooth normally elliptic extension, and negative-part testing preserves nonnegativity. Scalar Harnack on the joined strip yields positivity at the joining time for each species with positive conserved integral. This proves the restart assertion.
\end{proof}

\paragraph{Local regularity of bounded-gradient weak solutions.}
For the autonomous whole-space equation $U_t=\Div(A(U)\nabla U)$, let $U$ be a local-energy solution with bounded spatial gradient and state range as above. The bounded flux gives the preceding time modulus by spatial mollification. On a small cylinder fix $S_0=D^2H_\gamma(U_0)$ at its central state. By continuity and Lemma~\ref{lem:shifted-entropy}, $S_0A(U)$ is coercive there, on the tangent plane for the frozen equation. For a coordinate direction $e_k$ put $v_h=(U(\cdot+he_k)-U)/h$. Its divergence equation has principal matrix $A(U(\cdot+he_k))$ and additional flux $(D_h A(U))\nabla U$ bounded independently of $h$. With a spatial and lower-time cutoff $\zeta$, test by $\zeta^2 S_0v_h$, projecting onto the tangent plane in the frozen case. For fixed $h$ the time derivative belongs locally to $L^2_tH^{-1}_x$, so time averaging justifies this test. Coercivity and Young's inequality give
\[
 \sup_{s\le t}\int\zeta^2 v_h^\top S_0v_h+
 c\iint^{t}\zeta^2|\nabla v_h|^2
 \le C\iint^{t}(\zeta^2+|\nabla\zeta|^2+|\zeta\partial_t\zeta|)
                        (1+|v_h|^2).
\]
The integrals start at the lower face. The bound is uniform for small $h$ and for $t$ approaching the upper face. Since $v_h$ is bounded, $D^2U\in L^2$ on smaller backward cylinders, and expansion of the equation gives $U_t\in L^2$ there. The common-forcing argument above now applies at exponents $2$ and $q$, giving $U\in W_q^{2,1}$ on smaller cylinders for every $1<q<\infty$. For $q>d+2$ the spatial gradient has a continuous representative through the upper face.

\paragraph{From a H\"older bound to a gradient bound.}
We also use the local consequence of the preceding proof on buffered backward cylinders. Suppose a family has a uniform parabolic H\"older bound on a larger cylinder and has gradients unbounded on a strictly smaller cylinder in its interior, allowing both to meet the physical wall and have the same upper time face. Choose an intermediate cylinder with positive spatial and lower-time buffers in the larger one and containing the smaller one with such buffers. Truncate before the common upper time face and maximize \(|\nabla U|D\) on the truncated cylinder, where \(D\) is the minimum of the distance to the intermediate cylinder's artificial lateral boundary and the square root of the time since its lower face. The maximizing products tend to infinity because $D$ has a positive lower bound on the smaller cylinder. The maximizing points can be chosen with a fixed comparison factor if the maximum is not attained. On a backward box of radius \(cD\) around such a point, with fixed sufficiently small \(c>0\), \(D\) is at least a fixed multiple of its value at the point, so the gradients there are bounded by a constant times the selected gradient. Dilating by the reciprocal of this gradient gives a unit gradient and bounded gradients on boxes exhausting the past. The inherited H\"older modulus gives a constant state limit on each compact box. The compactness just proved, applied in the interior or with reflection, contradicts the unit gradient. In radius-one coordinates this argument gives a uniform bound on the smaller cylinder.

\begin{lemma}[A weight for no-flux localization]\label{lem:neumann-weight}
For centers $x_0\in\bar\Omega$ and scales $0<h\le1$, set $\Omega_h=(\Omega-x_0)/h$. There are positive weights $\phi$ with uniform mass and local lower bounds, $|D\phi|+|D^2\phi|\le C\phi$, and $\partial_n\phi\ge0$ on the physical wall. They have uniform exponential tails at small scales. If $v\ge0$, $v(0)\ge a>0$ and $|\nabla v|\le M$, then $v\ge a/2$ on a sufficiently small intersected ball about the center. The weighted integral of $v$ there is bounded below by a positive constant times the intersected volume. The assertion applies at wall centers.
\end{lemma}

\begin{proof}
The signed distance $d_\Omega$ is negative inside and $C^2$ in a fixed collar \cite{foote}. Let $\eta_{\mathrm{wall}}$ be smooth and bounded, equal to the identity near zero and constant on $(-\infty,-2]$. Use
\[
 \phi(X)=\exp\!\left[-\sqrt{1+|X|^2}+\eta_{\mathrm{wall}}(d_\Omega(x_0+hX)/h)\right].
\]
At small $h$ continue by the plateau beyond the collar. Then $\partial_n\phi=\phi(1-X\cdot n/\sqrt{1+|X|^2})\ge0$. The derivative, mass, and lower bounds follow directly. At larger fixed scales the constant weight suffices. Uniform boundary charts give the requisite lower volume estimates and paths of controlled length from a center to the points of a small truncated ball. For a scalar nonnegative component with no normal flux, the second integration by parts of its diffusion term against $\phi$ has boundary contribution $-\int a u\,\partial_n\phi\le0$. For entropy tests only the first integration is needed and the physical boundary flux vanishes.
\end{proof}

For a pole $(x_0,t_0)$, $x_0\in\bar\Omega$, let
\[
 D_R=\Omega\cap B_R(x_0),\qquad Q_R^\Omega=D_R\times(t_0-R^2,t_0).
\]
Write $Q_R$ for the Euclidean cylinder. We take cylinders meeting the wall below a uniform chart radius. Their volumes are comparable to $R^{d+2}$.

Let $\mathcal H_R$ be the closure in $H^1(D_R)$ of functions vanishing near the artificial boundary, with free trace at the physical wall. The norm $\|\cdot\|_{-1;D_R}$ is dual to the gradient seminorm on these tests. For a divergence source this gives
\[
 \|\Div F\|_{-1;D_R}
   =\sup_{\substack{\zeta\in\mathcal H_R\\\|\nabla\zeta\|_2\le1}}
            \left|\int_{D_R}F\cdot\nabla\zeta\right|.
\]
In a wall chart the source includes the volume density in the variational identity. Pairing its reflection against a compact test is equivalent to pairing the original source against the sum of the two half-chart restrictions pulled back to the physical domain. Their gradient norms are bounded by the gradient norm of the reflected test times a constant from the chart. This applies on buffered subcylinders whose half-chart restrictions stay away from the original artificial boundary.

Under $x=x_*+hX$, a divergence source $g$ becomes $g_h=h^2g(x_*+hX)$, in the distributional sense. If $x_*+hD'\subset D_R$, extension by zero across the artificial boundaries gives
\[
 \|g_h\|_{-1;D'}\le h^{1-d/2}\|g\|_{-1;D_R}.
\]
Indeed the source pairing has factor $h^{2-d}$ and the pullback of a test of unit $X$-gradient norm has $x$-gradient norm $h^{d/2-1}$. The squared norm integrated in rescaled time has factor $h^{-d}$. We use ordinary $H^{-1}$ norms on compact interior subsets for local time derivatives.

\begin{lemma}[Small full-state excess]\label{lem:small-excess}
Use the state ranges and compact positive parameter families of Lemmas~\ref{lem:shifted-entropy} and~\ref{lem:wall-record}. For sufficiently small mean squared space-time deviation of a strong solution $U$ from its mean on a cylinder $Q_R^\Omega$,
\[
 R\sup_{Q_{R/2}^\Omega}|\nabla U|\le K_{\mathrm{ex}} .
\]
For cylinders meeting the wall, the upper bound on \(R\) depends on the domain charts and the uniform state and coefficient bounds. In the unforced whole-space equation \(R\) is arbitrary. For the frozen equation with variable prescribed density, \(R\) is taken sufficiently small, with the bound also depending on the common positive bounds and \(C^2\) bound of the profile family.
\end{lemma}

\begin{proof}
Work first in radius-one coordinates for $Q_R^\Omega$ and let $a$ be a constant state in the band or simplex. In the frozen equation write $\rho_R(X)=\rho^0(x_0+RX)$. Its rescaled form is
\[
 \rho_R U_s=\Div_X\big(\rho_R A(U)\nabla_X U+K(U)\nabla_X\rho_R\big).
\]
The extra flux is $O(R)$ and the time weight is bounded above and away from zero. For the full equation the time weight is one and there is no extra flux. Testing with $\zeta^2(DH_\gamma(U)-DH_\gamma(a))$, with $\zeta$ supported away from the artificial lateral boundary and vanishing at the lower face, gives the energy inequality for $H_\gamma(U\mid a)$. The physical wall contributes no term. In the frozen case the time weight stays on the entropy term and the test uses the tangential coercivity of Lemma~\ref{lem:shifted-entropy}. Young's inequality bounds the product of the extra flux with the differentiated test by a fraction of the gradient term, a term in $|U-a|^2$, and $CR^2$. Thus on a concentric subcylinder of any fixed relative radius $\sigma<1$, the time supremum of $\int|U-a|^2$ and the space-time integral of $|\nabla_X U|^2$ are bounded by
\[
 C_\sigma\left(\iint_{Q_1\ \text{in the rescaled domain}} |U-a|^2+R^2\right).
\]
Omit $R^2$ in the unforced case.

For a sequence of radii $R_j\to0$ with a wall at bounded normalized distance, use a nearest-wall chart as in Lemma~\ref{lem:wall-record}. In coordinates based at the pole the limiting wall may be offset from the origin. Reflection in that limiting plane maps the exterior side toward the pole, so both half-chart restrictions of tests supported in a strictly smaller box fit inside the outer cap, for large $j$. The near-identity chart distortion is absorbed by the fixed margin between the boxes. The local energy bounds therefore hold for the even extensions on fixed smaller boxes. The reflected divergence equation bounds $\partial_s U$ locally in $L^2(H^{-1})$ by dividing the test by the positive, time-independent Lipschitz weight and using the $L^2$ flux bound. Local compactness in $L^2$ follows from the space-gradient and time-derivative bounds. The estimates include the interval up to its upper face in the space-time norm, by the uniform time-supremum bound on a smaller box. Interior sequences use the same estimates directly.

Let $e(R)$ be the root mean squared space-time excess at radius $R$. Fix $0<\beta<1$ and choose $\vartheta<1/4$ so that $C\vartheta<\vartheta^\beta/2$ for the constant-coefficient interior caloric excess bound. We claim that for sufficiently small $e(R)$ and, where needed, sufficiently small wall or prescribed-profile scale $R$,
\[
 e(\vartheta R)\le \tfrac12\vartheta^\beta e(R)
 \quad\text{or}\quad
 e(\vartheta R)\le C_\vartheta R .
\]

The second alternative is used only for variable prescribed density. If no cutoffs make the dichotomy valid, choose a failing sequence with \(e_j=e(R_j)\le 1/j\) and \(R_j\le 1/j\) in the cases requiring a small absolute radius. In the forced case take \(C_{\vartheta,j}=j\). Failure of the first inequality implies \(e_j>0\). Restriction and volume comparison give \(e(\vartheta R_j)\le C(\vartheta)e_j\). Thus in the forced case \(jR_j<e(\vartheta R_j)\le C(\vartheta)e_j\), and \(R_j/e_j\to0\).

Let \(a_j\) be the mean on the outer physical space-time cylinder and set \(V_j=(U_j-a_j)/e_j\) in radius-one coordinates. The energy and time-derivative estimates give strong local \(L^2\) convergence of \(V_j\) and weak local \(L^2\) convergence of \(\nabla V_j\), using reflection at bounded normalized wall distance. Along a subsequence \(a_j\to a\), and \(U_j\to a\) almost everywhere on those boxes. Take convergent subsequences of the parameters and, in the frozen case, of the prescribed density values at the poles. Denote the diffusion matrix at the limiting parameters by \(A_\infty\). In the reflected divergence equation the principal flux tensor for \(V_j\) is \(w_j g_j^{ab}A_j(U_j)\), with \(w_j=\mathcal J_j\) for the full system and \(w_j=\mathcal J_j\rho_{R_j}\) for the frozen system. In interior coordinates use the identity metric and omit \(\mathcal J_j\). The wall scales and the variable-density scales tend to zero in the respective cases, so the metric and the weight converge uniformly on fixed boxes to constant factors \(g_\infty^{ab}\) and \(w_\infty>0\). For each smooth compact test \(\varphi\), with summation over spatial indices, boundedness and almost-everywhere convergence give
\[
 \left\|\big(w_jg_j^{ab}A_j(U_j)-w_\infty g_\infty^{ab}A_\infty(a)\big)^\top
                \partial_a\varphi\right\|_2\longrightarrow0 .
\]
Together with weak convergence of \(\partial_b V_j\) this identifies the normalized principal flux in the weak formulation. The time term passes by strong convergence of \(w_j V_j\). The normalized extra flux in the forced case is \(O(R_j/e_j)\), so it vanishes. After cancelling the constant weight and using the orthonormal frame for the limiting metric, the limit solves the constant-coefficient heat system on a fixed interior box.

The limit has a structural local \(L^2\) bound. At the limiting state \(a\) the constant Hessian \(D^2H_\gamma(a)\) at the limiting parameters gives a coercive quadratic form when multiplied by \(A_\infty(a)\) on the state space. Testing successive spatial difference quotients with this fixed Hessian and cutoffs gives interior \(L^2\) bounds for spatial derivatives of every order on smaller boxes. The constant equation supplies time derivatives, and Sobolev embedding gives the linear interior derivative bound. The derivative bounds near the pole give excess at most \(C\vartheta\) on the limiting radius-\(\vartheta\) cap. Uniform convergence of the rescaled boundary graphs on compact sets gives almost-everywhere convergence of the indicators of the smaller physical cylinders, since the limiting plane and sphere have zero volume. In a common buffered box write these cylinders as \(Q_j\) and \(Q_\infty\), and the normalized deviations as \(V_j\to V\) in \(L^2\). Then
\[
 \|\mathbf1_{Q_j}V_j-\mathbf1_{Q_\infty}V\|_2
 \le \|V_j-V\|_2+\|(\mathbf1_{Q_j}-\mathbf1_{Q_\infty})V\|_2
 \longrightarrow 0.
\]
The second term tends to zero by dominated convergence. The normalized chart Jacobians converge uniformly to one, so the cap volumes and the first and second moments converge. This contradicts failure of the first inequality.

For centers in $Q_{3R/4}^\Omega$ start the iteration at $R/8$. Every starting cylinder has a fixed buffer in $Q_R^\Omega$, and restriction bounds its excess by a structural multiple of the outer excess. On shrinking concentric cylinders of radii $r_j=\vartheta^j R/8$ the dichotomy gives
\[
 e(r_{j+1})\le \tfrac12\vartheta^\beta e(r_j)+C_\vartheta r_j,
\]
with the last term absent without forcing. Choose the outer excess threshold and, in the forced case, the initial radius small enough to keep all iterates in the region of validity. Iteration and comparison of intermediate radii yield the uniform Campanato bound $e(r)\le C(r/R)^\beta(e(R)+R)$ in the forced case and the bound with $e(R)$ alone without forcing. Volume comparability and the comparison of averages on intersecting cylinders give a parabolic H\"older modulus in radius-one coordinates on $Q_{3R/4}^\Omega$, through its upper face. Apply the H\"older-to-gradient argument of Lemma~\ref{lem:wall-record} on this larger cylinder, with the smaller one equal to $Q_{R/2}^\Omega$. This gives the claimed bound. Bounds on a binary face and in fixed affine coordinates follow from the corresponding restricted entropy estimate.
\end{proof}

\begin{lemma}[From a vanishing cylinder defect to ancient rigidity]\label{lem:one-cylinder}
Consider a nonnegative local defect on cylinders for the equations of Lemma~\ref{lem:small-excess}. Assume that restriction and rescaling to a smaller cylinder of scale ratio $h$ increases the normalized defect by at most $C h^{-\beta_{\mathrm{sc}}}$ for some fixed $\beta_{\mathrm{sc}}>0$. Suppose that vanishing defect in the rescaled sequence imposes a constraint on its local-energy ancient limit, and that the only bounded ancient solutions in the state range with that constraint are constant. Then sufficiently small defect on $Q_R^\Omega$ gives $R\sup_{Q_{R/2}^\Omega}|\nabla U|\le K$. In particular a vanishing rescaled $L^2_t H^{-1}_x$ residual for a bounded observed scalar in a scalar comparison equation with uniformly positive bounded measurable conductivity forces the limiting scalar to be constant.
\end{lemma}

\begin{proof}
If no constants work, use outer defects $\delta_j\to0$ and poles in the inner half with gradients of normalized magnitude $M_j\to\infty$. Measure distances in units of the outer radius and then dilate by
\[
 h_j=\max\{M_j^{-1/2},\,\delta_j^{1/(2\beta_{\mathrm{sc}})}\}.
\]
Both $h_jM_j\to\infty$ and $h_j^{-\beta_{\mathrm{sc}}}\delta_j\to0$. The rescaled cylinders exhaust the past, and on every fixed cylinder the restricted defects tend to zero. Lemma~\ref{lem:small-excess} gives a positive lower bound for the full-state variance in a fixed unit-size cylinder centered at each pole, since its normalized gradient is unbounded on the inner half. The local entropy estimates and compactness from its proof apply to the bounded states. When a wall stays at bounded normalized distance, reflect in the chart on strictly smaller boxes and pass to the limiting plane. The variance lower bound survives the limit, by strong local $L^2$ convergence and convergence of the caps in measure. There is consequently a nonconstant entire weak ancient limit. Applying the energy estimate on a cylinder of radius $2R$ and scaling it back gives $\iint_{Q_R}|\nabla U|^2\le C R^d$ for this limit.

For the scalar assertion, write on a fixed large reflected cylinder the observed scalar as $V_j$, the equation as
\[
 w_j\partial_s V_j-\Div(B_j\nabla V_j)=g_j ,
\]
and let $v_j$ solve the homogeneous equation with the same parabolic boundary data. Here $w_j$ is time independent, bounded above and away from zero, and $B_j$ is the scalar conductivity times the uniformly positive geometric tensor and volume factor, or the conductivity times the identity in the interior. The difference $e_j=V_j-v_j$ has zero parabolic data. Its energy estimate is
\[
 \sup_s\int w_j |e_j|^2+\iint|\nabla e_j|^2
 \le C\int\|g_j(s)\|_{H^{-1}}^2\,\dd s .
\]
This follows by testing with $e_j$, and also gives the variational construction of $v_j$ by solving the zero-data problem for the difference. By the scalar maximum principle $v_j$ is uniformly bounded. The scalar oscillation estimate is uniform in the ellipticity bounds and upper and lower bounds of $w_j$. In the truncation estimates the time-independent weight multiplies the squared truncations on the time slices, where the weighted and unweighted integrals are equivalent. On each fixed smaller backward cylinder $v_j$ converges in $L^2$ to the same observed limit as $V_j$. The oscillation bound therefore passes to the limit, with the radius ratio of the smaller to the larger cylinder raised to a fixed positive power. Letting the larger radius increase proves constancy. This comparison requires no identification of a limit conductivity.

Under the stated rigidity hypothesis the constraint now forces the entire state limit to be constant, contradicting its positive variance.
\end{proof}

\section{Global solutions at finite physical density}
\label{sec:finite-density}

Here $r^0\in X_\gamma=W_{\mathrm N}^{2-2/p,p}(\Omega;\mathbb R^3)$, $p>d+2$, is strictly positive on $\overline\Omega$, and $\lambda,\nu>0$. All inverse masses and frictions are positive. We write $\langle r\rangle$ for the normalized spatial mean. All constants may depend on the specified coefficient set, $\Omega$ and the prescribed density or lens ranges. The thresholds below have no dependence on the initial gradients or on the individual species minima.

\begin{theorem}[Finite density]\label{thm:finite-density}
The following classes give global solutions with the classical positive-time regularity specified above.
\begin{enumerate}
\item Fix $0<L_0<U_0<\infty$ and a species $i$. Fix positive inverse masses and frictions. There exists $\delta_0>0$ such that
\[
             L_0\le \rho^0\le U_0,\qquad
             \left\|\sum_{j\ne i}r_j^0\right\|_\infty\le\delta_0
\]
suffice.
\item Fix $m_0>0$ and a compact set $\mathcal K\subset(0,\infty)^4$ of allowed values of $(m_3,f_{12},f_{13},f_{23})$. For the same initial density band there exists $\delta_0>0$, uniform over $\mathcal K$, such that
\[
 \|r_3^0\|_\infty+|m_1-m_0|+|m_2-m_0|\le\delta_0
\]
suffices. The two abundant species can occur in arbitrary proportions. The coefficient class includes triples with all three inverse masses different.
\item Fix positive inverse masses with $m_1=a<m_2=b$ and positive frictions. Set
\[
        k=\frac{\lambda f_{12}}{\nu},\qquad
        (X,Y)=k\big(r_1+r_2,\,ar_1+br_2\big).
\]
Take two apices $P_i=(X_*,Y_i)$ and $P_o=(X_*,Y_o)$ with $aX_*<Y_i<Y_o<bX_*$. Write $\mathcal L(P)$ for the closed convex region defined in Section~\ref{sec:binary-lens} by the two characteristic arcs through $P$. There exists $\delta_0>0$ such that $(X^0,Y^0)\in\mathcal L(P_i)$ and $\|r_3^0\|_\infty\le\delta_0$ suffice.
\end{enumerate}
In each case the total density stays in a compact subinterval of $(0,\infty)$, and $\|\nabla r(t)\|_\infty\le C(1+t^{-1/2})$ for $t>0$. These bounds are uniform over the permitted initial profiles and coefficients.
\end{theorem}

\subsection{The binary characteristic region}
\label{sec:binary-lens}

On the binary face $r_3=0$, the resistance on the zero-sum current of species 1 and 2 is $f_{12}$. In time $s=(\nu/f_{12})t$ the equation for $U=(X,Y)$ becomes
\[
 U_s=\operatorname{div}(A_{\mathrm b}(U)\nabla U),\qquad
 A_{\mathrm b}=\begin{pmatrix}0&X\\-ab&D\end{pmatrix},
 \qquad D=Y+a+b-Y/X .
\]
Indeed $X$ is proportional to the total density, $Y$ to the pressure, and $a^2r_1+b^2r_2=(a+b)(ar_1+br_2)-ab(r_1+r_2)$. The physical wedge is $X>0$, $aX\le Y\le bX$. In its interior set $\mu=Y/X$ and let $\eta_-<\eta_+$ be the roots of
\begin{equation}\label{eq:binary-characteristic}
             X\eta^2-D\eta+ab=0 .
\end{equation}
Evaluating at $\mu$ gives $(\mu-a)(\mu-b)<0$, so $0<\eta_-<\mu<\eta_+$. The graph $Y=F(X)$ solving $F'=\eta_\pm(X,F)$ has tangent eigenvector $(1,\eta_\pm)$ with eigenvalue $X\eta_\pm$, and differentiating \eqref{eq:binary-characteristic} along the graph gives
\[
       F''=\frac{\eta_\pm(\mu-\eta_\pm)}{X(2X\eta_\pm-D)}<0 .
\]
Follow the plus arc to the left of the apex $P$ and the minus arc to its right. Along either arc, $\mu'=(\eta_\pm-\mu)/X$. Thus $\mu$ decreases in the direction from $P$ toward the relevant intercept and stays below $b$. For $X\le1$ we have $D\ge a$ and $\eta_+\ge a/(2X)$. For $X\ge1$ we have $D\ge aX$ and $\eta_-\le 2b/X$. Integration of these two bounds, using $aX\le Y\le bX$, shows that each arc meets the lower ray, at $X=A_0(P)$ for the plus arc and $X=B_0(P)$ for the minus arc. The same estimates keep the intercepts on any compact family of apices away from zero and infinity. Indeed, with $x_b=\min(1,X(P))$ and $x_a=\max(1,X(P))$, they give
\[
 \tfrac a2\log(x_b/A_0)\le b x_b \quad\text{if }A_0<x_b,\qquad
 aB_0-bx_a\le 2b\log(B_0/x_a) \quad\text{if }B_0>x_a.
\]
These bounds also hold if the apices tend to a point of the lower ray with $X>0$.

The two slopes on the lower ray are $a$ and $b/X$. Away from $X_c=b/a$ the roots are simple. Uniqueness of the corresponding scalar ODE prevents an interior arc from meeting the lower ray on a branch whose slope there is $a$. To rule out a contact at $X_c$, write $h=X-X_c$, $w=Y-aX$ and $\sigma=\eta-a$. Equation \eqref{eq:binary-characteristic} becomes
\[
          \sigma(ah+X\sigma)=(1-1/X)w(a+\sigma).
\]
A plus arc reaching $X_c$ from the right has $h,\sigma>0$. A minus arc reaching it from the left has $h,\sigma<0$ near the contact, since the characteristic polynomial at $a$ equals $-(1-1/X)aw<0$ there. In both cases $|ah+X\sigma|\ge a|h|$. Since $\sigma\to0$ and $1-1/X_c<1$, for some $\vartheta<1$ we have
\[
 \frac{|w'|}{w}=\frac{(1-1/X)(a+\sigma)}{|ah+X\sigma|}
       \le \frac{\vartheta}{|h|}
\]
near the contact. Integration along the arc gives $w\ge c|h|^\vartheta$, whereas $w'\to0$ and $w\to0$ give $w=o(|h|)$. Thus the contacts are transverse and $A_0(P)<X_c<B_0(P)$.

Define $\mathcal L(P)$ as the region above the lower ray and below the two arcs, including their boundaries. At the apex the slope drops from $\eta_+$ to $\eta_-$, so the upper envelope is concave and the region is convex. The lens is compact, has a positive lower bound for $X$, and stays strictly below the upper ray. The characteristic coordinates in the interior of the wedge are
\[
 u_+(U)=X_c-A_0(U),\qquad u_-(U)=B_0(U)-X_c .
\]
Their extensions to the lower ray are $(X_c-X)_+$ and $(X-X_c)_+$, respectively. To see this for the plus coordinate, take $U_j=(X_j,Y_j)\to(X_0,aX_0)$ and a subsequence on which $A_0(U_j)\to A$. The intercept bounds just proved give $0<A\le X_0$, and $A\le X_c$. Since $\mu$ increases from the intercept to $U_j$, the intervening segment tends uniformly to the lower ray. If $A<X_c$ and $A<X_0$, continuous dependence at the transverse contact gives positive separation from that ray at some fixed $X\in(A,X_0)$, a contradiction. Thus any limit with $A<X_c$ equals $X_0$. If $X_0<X_c$, the inequality $A\le X_0$ forces this case. If $X_0\ge X_c$, every limiting intercept equals $X_c$. For the minus coordinate, a subsequential limit $B$ of $B_0(U_j)$ satisfies $B\ge X_0$ and $B\ge X_c$. The same argument now uses continuous dependence to the left of a transverse contact to rule out $B>X_c$ and $B>X_0$. This gives $B=\max(X_0,X_c)$.

On compact subsets where $u_\pm$ is positive, the corresponding coordinate is smooth, including at the transverse lower-ray contacts, and $\partial_Yu_\pm>0$. For the plus graph parametrized by its intercept $A$, differentiation of the ODE in $A$ gives a homogeneous linear variational equation with initial value $a-\eta_+(A,aA)<0$. Its sign persists. The analogous initial value for the minus graph, parametrized by its right intercept $B$, is $a-\eta_-(B,aB)>0$. Along the plus arc through $P$, $u_+$ is fixed, while $u_-$ increases with $X$ by
\[
 Du_-\,(1,\eta_+)=(\partial_Yu_-)(\eta_+-\eta_-)>0
\]
in the interior. Along the minus arc, $u_-$ is fixed and $u_+$ decreases with $X$. Consequently, both sublevel conditions $u_\pm\le u_\pm(P)$ hold on the upper envelope and below it. For $X\notin[A_0(P),B_0(P)]$, at least one condition fails already on the lower ray. For $X\in[A_0(P),B_0(P)]$, points above the envelope violate the condition for the family defining the envelope there, by monotonicity in $Y$. This identifies $\mathcal L(P)$ with the two sublevel conditions below the upper ray. At the common $X_*$, increasing $Y$ from $Y_i$ to $Y_o$ increases both coordinates strictly, giving the two strict upper margins in Theorem~\ref{thm:finite-density}.

\begin{figure}[htbp]
\centering
\begin{tikzpicture}[x=.83cm,y=3.0cm,>=stealth]
\fill[teal!13] plot coordinates {(0.5427,0.0000) (0.5783,0.1223) (0.6138,0.2175) (0.6494,0.2930) (0.6849,0.3540) (0.7204,0.4039) (0.7560,0.4452) (0.7915,0.4799) (0.8271,0.5093) (0.8626,0.5345) (0.8982,0.5563) (0.9337,0.5752) (0.9692,0.5919) (1.0048,0.6066) (1.0403,0.6196) (1.0759,0.6313) (1.1114,0.6418) (1.1470,0.6513) (1.1825,0.6598) (1.2180,0.6677) (1.2536,0.6748) (1.2891,0.6814) (1.3247,0.6874) (1.3602,0.6930) (1.3958,0.6982) (1.4313,0.7030) (1.4669,0.7075) (1.5024,0.7116) (1.5379,0.7155) (1.5735,0.7192) (1.6090,0.7226) (1.6446,0.7259) (1.6801,0.7289) (1.7157,0.7318) (1.7512,0.7345) (1.7867,0.7371) (1.8223,0.7395) (1.8578,0.7418) (1.8934,0.7440) (1.9289,0.7461) (1.9645,0.7481) (2.0000,0.7500) (2.1952,0.6700) (2.3905,0.6025) (2.5857,0.5448) (2.7810,0.4949) (2.9762,0.4513) (3.1715,0.4129) (3.3667,0.3789) (3.5620,0.3484) (3.7572,0.3210) (3.9525,0.2962) (4.1477,0.2737) (4.3430,0.2532) (4.5382,0.2343) (4.7335,0.2170) (4.9287,0.2010) (5.1239,0.1862) (5.3192,0.1724) (5.5144,0.1596) (5.7097,0.1476) (5.9049,0.1364) (6.1002,0.1259) (6.2954,0.1160) (6.4907,0.1067) (6.6859,0.0980) (6.8812,0.0897) (7.0764,0.0818) (7.2717,0.0744) (7.4669,0.0673) (7.6621,0.0606) (7.8574,0.0542) (8.0526,0.0482) (8.2479,0.0423) (8.4431,0.0368) (8.6384,0.0315) (8.8336,0.0264) (9.0289,0.0216) (9.2241,0.0169) (9.4194,0.0124) (9.6146,0.0081) (9.8099,0.0040) (10.0051,0.0000)} -- cycle;
\draw[teal!65!black,line width=1pt] plot coordinates {(0.5427,0.0000) (0.5783,0.1223) (0.6138,0.2175) (0.6494,0.2930) (0.6849,0.3540) (0.7204,0.4039) (0.7560,0.4452) (0.7915,0.4799) (0.8271,0.5093) (0.8626,0.5345) (0.8982,0.5563) (0.9337,0.5752) (0.9692,0.5919) (1.0048,0.6066) (1.0403,0.6196) (1.0759,0.6313) (1.1114,0.6418) (1.1470,0.6513) (1.1825,0.6598) (1.2180,0.6677) (1.2536,0.6748) (1.2891,0.6814) (1.3247,0.6874) (1.3602,0.6930) (1.3958,0.6982) (1.4313,0.7030) (1.4669,0.7075) (1.5024,0.7116) (1.5379,0.7155) (1.5735,0.7192) (1.6090,0.7226) (1.6446,0.7259) (1.6801,0.7289) (1.7157,0.7318) (1.7512,0.7345) (1.7867,0.7371) (1.8223,0.7395) (1.8578,0.7418) (1.8934,0.7440) (1.9289,0.7461) (1.9645,0.7481) (2.0000,0.7500) (2.1952,0.6700) (2.3905,0.6025) (2.5857,0.5448) (2.7810,0.4949) (2.9762,0.4513) (3.1715,0.4129) (3.3667,0.3789) (3.5620,0.3484) (3.7572,0.3210) (3.9525,0.2962) (4.1477,0.2737) (4.3430,0.2532) (4.5382,0.2343) (4.7335,0.2170) (4.9287,0.2010) (5.1239,0.1862) (5.3192,0.1724) (5.5144,0.1596) (5.7097,0.1476) (5.9049,0.1364) (6.1002,0.1259) (6.2954,0.1160) (6.4907,0.1067) (6.6859,0.0980) (6.8812,0.0897) (7.0764,0.0818) (7.2717,0.0744) (7.4669,0.0673) (7.6621,0.0606) (7.8574,0.0542) (8.0526,0.0482) (8.2479,0.0423) (8.4431,0.0368) (8.6384,0.0315) (8.8336,0.0264) (9.0289,0.0216) (9.2241,0.0169) (9.4194,0.0124) (9.6146,0.0081) (9.8099,0.0040) (10.0051,0.0000)};
\draw[navy,dashed,line width=.7pt] plot coordinates {(0.4600,0.0000) (0.4976,0.1766) (0.5351,0.3035) (0.5727,0.3975) (0.6103,0.4689) (0.6478,0.5243) (0.6854,0.5682) (0.7229,0.6036) (0.7605,0.6325) (0.7981,0.6565) (0.8356,0.6766) (0.8732,0.6937) (0.9107,0.7084) (0.9483,0.7211) (0.9859,0.7322) (1.0234,0.7420) (1.0610,0.7507) (1.0985,0.7584) (1.1361,0.7653) (1.1737,0.7715) (1.2112,0.7772) (1.2488,0.7823) (1.2863,0.7870) (1.3239,0.7913) (1.3615,0.7952) (1.3990,0.7989) (1.4366,0.8022) (1.4741,0.8054) (1.5117,0.8083) (1.5493,0.8110) (1.5868,0.8135) (1.6244,0.8159) (1.6620,0.8181) (1.6995,0.8202) (1.7371,0.8222) (1.7746,0.8241) (1.8122,0.8258) (1.8498,0.8275) (1.8873,0.8291) (1.9249,0.8306) (1.9624,0.8320) (2.0000,0.8333) (2.2086,0.7398) (2.4171,0.6617) (2.6257,0.5957) (2.8342,0.5391) (3.0428,0.4900) (3.2513,0.4470) (3.4599,0.4090) (3.6684,0.3752) (3.8770,0.3450) (4.0855,0.3178) (4.2941,0.2931) (4.5026,0.2707) (4.7112,0.2501) (4.9197,0.2313) (5.1283,0.2140) (5.3368,0.1980) (5.5454,0.1832) (5.7539,0.1694) (5.9625,0.1565) (6.1710,0.1445) (6.3796,0.1333) (6.5881,0.1227) (6.7967,0.1128) (7.0052,0.1035) (7.2138,0.0947) (7.4223,0.0863) (7.6309,0.0784) (7.8394,0.0709) (8.0480,0.0638) (8.2565,0.0571) (8.4651,0.0507) (8.6736,0.0445) (8.8822,0.0387) (9.0907,0.0331) (9.2993,0.0277) (9.5078,0.0226) (9.7164,0.0177) (9.9249,0.0130) (10.1335,0.0085) (10.3420,0.0042) (10.5506,0.0000)};
\draw[gray,dotted] (2,0)--(2,.8333);
\fill[teal!65!black] (2,.75) circle(1.8pt);
\draw[->] (0,0)--(11.1,0) node[right] {$X$};
\draw[->] (0,0)--(0,1.06) node[above] {$\theta$};
\foreach \x in {2,4,6,8,10}{\draw (\x,0)--(\x,-.022) node[below] {\small\x};}
\foreach \y/\txt in {.25/0.25,.5/0.50,.75/0.75,1/1.00}{\draw (0,\y)--(-.08,\y) node[left]{\small \txt};}
\node[below left] at (0,0) {\small 0};
\draw[teal!65!black,line width=1pt] (5.1,.86)--(5.9,.86) node[right,black] {\small inner boundary};
\draw[navy,dashed,line width=.7pt] (5.1,.68)--(5.9,.68) node[right,black] {\small outer boundary};
\end{tikzpicture}
\caption{Admissible initial states for $m=(1,4,9)$, $(f_{12},f_{13},f_{23})=(2,5,17)$, $\lambda=1$ and $\nu=2$. Here $k=\lambda f_{12}/\nu=1$, so $X=r_1+r_2$ and $\theta=r_2/(r_1+r_2)$. The characteristic arcs have apices $(2,6.5)$ and $(2,7)$ in $(X,Y)$ coordinates. Their intersections with $\theta=0$ occur at $X\approx0.543$ and $10.005$ for the inner region and at $X\approx0.460$ and $10.551$ for the outer region. With these apices the scaled arcs depend only on $m_1,m_2$. The factor $k$ sets the conversion to physical pair density, while $m_3,f_{13},f_{23}$ enter the ternary smallness threshold. Strictly positive data in the shaded inner region with sufficiently small $\|r_3^0\|_\infty$ satisfy Theorem~\ref{thm:finite-density}(3).}
\label{fig:physical-lens}
\end{figure}
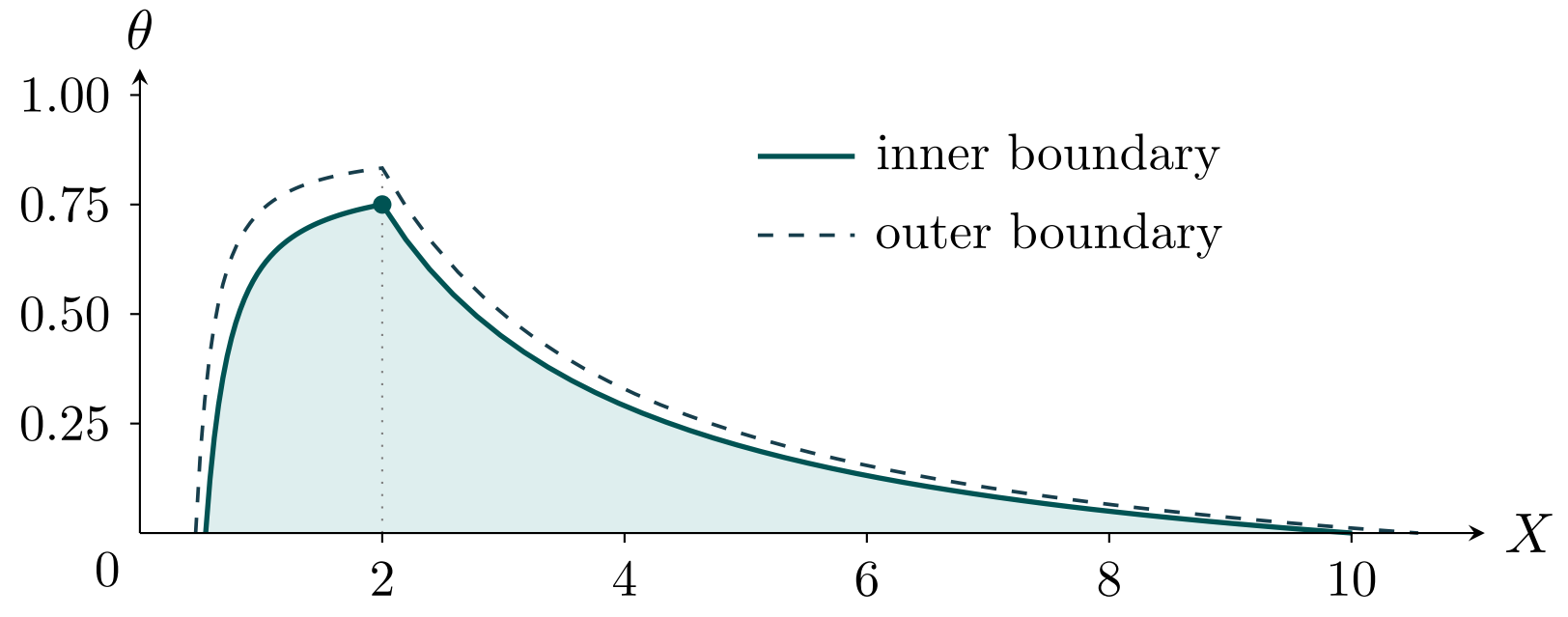

Figure~\ref{fig:physical-lens} shows the two nested regions in total-density and mass-fraction coordinates.

Fix a lens and either coordinate $u=u_\pm$. For $\ell>0$ put
\[
       H(U)=\psi(u(U)),\qquad
       \psi(v)=e^{K(v-\ell)_+}-1-K(v-\ell)_+ .
\]
For every $K>0$ the function is $C^{1,1}$ on the lens. For each fixed $\ell>0$ choose $K=K(\ell)$ sufficiently large. The entropy balance then has nonnegative production and reads
\begin{equation}\label{eq:binary-scalar-entropy}
 H_s-\operatorname{div}(\zeta(U)\nabla H)=-\mathcal P_H,\qquad
 \mathcal P_H=\sum_jD^2H[\partial_jU,A_{\mathrm b}\partial_jU]\ge0.
\end{equation}
Here $\zeta$ is the other eigenvalue from the one tangent to the level of $u$, so $D H A_{\mathrm b}=\zeta D H$ and $\zeta$ stays in a positive finite interval. Assign any value in that interval on the inactive set $u\le\ell$. For the production, in the tangent direction $v=(1,\eta_\pm)$,
\[
        D^2u[v,A_{\mathrm b}v]=-X\eta_\pm u_Y F''>0 .
\]
On $u\ge\ell$ this is uniformly positive. In a transverse frame with $D u\,v=0$, $D u\,w=1$, all other coefficients of the Hessian term with $\psi'$ are bounded. The rank-one Hessian term with $\psi''$ is precisely $\psi''\zeta|\nabla u|^2$. Since $\psi''\ge K\psi'$ on $u>\ell$, Young's inequality proves \eqref{eq:binary-scalar-entropy}, with production bounding $c_\ell|\nabla U|^2$ on $u\ge2\ell$ and $c_\ell\psi''|\nabla u|^2$ on $u>\ell$. The same argument gives $|\psi'|^2/\psi''\le C\psi$ on the active set. Cutoffs can be justified by approximation on that smooth portion of the lens, since $\psi=\psi'=0$ at the activation threshold and locally the chain identity holds for $C^{1,1}$ functions. In particular the upper-level truncations at $u_\pm(P)$ preserve the lens for the binary equation, by integration and a first-exit argument on a slightly larger lens. The required no-flux condition for this entropy is inherited from $\partial_n U=0$.

\begin{lemma}[Ancient binary rigidity in a characteristic region]
\label{lem:binary-ancient}
A bounded ancient binary solution on $\mathbb R^d$ whose range is contained in $\mathcal L(P)$ and that has bounded spatial gradient is constant.
\end{lemma}

\begin{proof}
Test \eqref{eq:binary-scalar-entropy} against the whole-space backward fundamental solution of the scalar uniformly elliptic divergence operator with pole at $(0,0)$. Its mass is one and the bounded entropy has nonnegative production, so $\mathcal P_H$ has finite integral against this kernel over the entire past. This follows first on finite intervals by scalar comparison with the forced solution and then by monotone convergence. The Gaussian lower bound \cite{Aronson} implies, for every fixed $0<a_1<b_1<\infty$ and $M<\infty$,
\[
 R^{-d}\int_{-b_1R^2}^{-a_1R^2}\int_{B_{MR}}\mathcal P_H
       \longrightarrow 0\qquad(R\longrightarrow\infty).
\]
Indeed on each displayed cylinder the kernel is bounded below by $cR^{-d}$, and its weighted production integral is a tail at times $\le -a_1R^2$. The constants are invariant under the parabolic dilation because the scalar coefficient remains in its fixed ellipticity interval.

Set $U_R(x,s)=U(Rx,R^2s)$. The relative shifted entropy on the binary face gives uniform local $L^2(H^1)$ bounds, and the conservation law gives bounds for $\partial_s U_R$ in $L^2(H^{-1})$ on smaller cylinders. Along a subsequence $U_R\to U_\infty$ strongly in local $L^2$ and almost everywhere, and the gradients converge weakly in local $L^2$. For a bounded continuous state cutoff $\eta_{\mathrm{loc}}$ supported where $u_\pm>2\ell$, the preceding tail bound and the quadratic production bound imply $\eta_{\mathrm{loc}}(U_R)\nabla U_R\to0$ in local $L^2$ on every slab lying strictly in the past. Testing the weak gradient limit against $\eta_{\mathrm{loc}}(U_R)$ times a fixed bounded test identifies the limit as $\eta_{\mathrm{loc}}(U_\infty)\nabla U_\infty$. The argument applies for each fixed $\ell>0$ along the same blow-down subsequence. Taking a sequence $\ell\downarrow0$ exhausts the sets $u_\pm>0$ and gives zero gradient except possibly at $U_\infty=(X_c,aX_c)$. The gradient is zero there as well by Sobolev locality applied to the two components. Boundedness and pointwise convergence of $A_{\mathrm b}(U_R)$ allow passage to the limit against the weak gradients in the conservative flux. Thus $\partial_s U_\infty=0$, and $U_\infty$ is constant. The uniformly bounded state range controls the contribution of an arbitrarily short terminal interval, so the strong $L^2$ convergence holds on a fixed cylinder up to the upper face. Its excess now tends to zero. The bounded-gradient solution is locally strong by the upgrade in Section~\ref{sec:boundary-estimates}. Lemma~\ref{lem:small-excess} bounds $R|\nabla U(0,0)|$ along the subsequence, and $R\to\infty$ gives zero gradient at the pole. Translations prove the result.
\end{proof}

\subsection{Confinement of the full ternary state}

In parts 1 and 2 choose a fixed enlarged positive total-density band containing $[L_0,U_0]$ in its interior. Require every minority component to stay below a small outer threshold $\eta$. In part 3 require $U\in\mathcal L(P_o)$ and $r_3\le\eta$, which in particular fixes a positive total-density band. These are first-exit conditions. On the resulting strips the wall record estimate of Lemma~\ref{lem:wall-record} implies
\begin{equation}\label{eq:finite-gradient-bound}
              \|\nabla r(t)\|_\infty\le C(1+t^{-1/2})
\end{equation}
uniformly once the outer minority threshold and, for part 2, the upper bound for $|m_1-m_0|+|m_2-m_0|$ are sufficiently small. If such a uniform bound failed, take these bounds tending to zero while keeping $\mathcal K$ fixed in part 2. The strong trace gives \(\min(1,\sqrt t)\|\nabla r(t)\|_\infty\to0\) as \(t\downarrow0\) for each profile. For each \(j\) take a first time \(t_j>0\) with \(\min(1,\sqrt{t_j})\|\nabla r(t_j)\|_\infty=j\), and set \(M_j=\|\nabla r(t_j)\|_\infty\). Then \(t_j M_j^2\ge j^2\). Rescale at a point attaining \(M_j\) by the spatial scale \(M_j^{-1}\), so the earlier gradients on every fixed rescaled backward cylinder are bounded by \(2\) for large \(j\). After taking a subsequence the coefficients converge to positive limits. In part 2, this uses compactness of $\mathcal K$ for $m_3$ and the frictions. Normal ellipticity is uniform on the fixed band. Lemma~\ref{lem:wall-record} gives a whole-space binary or single-species ancient limit with unit gradient at the selected point. When the rescaled distance to the wall stays bounded, the reflection uses the limiting plane at that distance. In part 1 the total density solves a scalar uniformly parabolic Darcy equation and hence is constant by the scalar interior oscillation estimate on arbitrarily large cylinders. In part 2 the abundant inverse masses both tend to $m_0$, giving $C=m_0\rho$ in the limit. The same scalar argument fixes the total density. The active partial densities then solve heat equations with diffusivity $\nu m_0/f_{12}^{\rm lim}$, where $f_{12}^{\rm lim}>0$ is the limiting binary friction. Thus they are constant. In part 3 the limit lies in the binary lens and Lemma~\ref{lem:binary-ancient} applies. These contradictions prove \eqref{eq:finite-gradient-bound}. Reducing the outer minority threshold or the upper bound for $|m_1-m_0|+|m_2-m_0|$ preserves the same estimate and constant. At a finite endpoint prior to exit the same bound gives a strong terminal trace and restart by Lemma~\ref{lem:wall-record}. We may therefore test a putative exit using the trace on the closed outer band. Species positivity at the trace follows from the scalar positivity argument in Section~\ref{sec:local-theory}.

For every minority component $u=r_j$, smoothness of the density diffusion matrix across its face gives
\[
 u_t=\operatorname{div}(a_j(r)\nabla u+uV_j),\qquad
 a_j\ge c>0,\qquad |V_j|+|\nabla a_j|\le C|\nabla r|.
\]
Here the inactive row has a strictly positive diagonal value and its off-diagonal entries vanish at $r_j=0$, so they divide smoothly by $r_j$. Let $x_*$ be a point of interest at time $t_*\le1$. Use the spatial scale $h=\sqrt{t_*}$ and time scale $h^2$, with the positive weight $\phi$ of Lemma~\ref{lem:neumann-weight} centered at $x_*$. Its derivatives through order two are bounded by $C\phi$ in the scaled domain and its outward normal derivative is zero or nonnegative. Zero normal total flux and integration by parts in the scalar diffusion give
\[
 \frac{\dd}{\dd s}\int_{\Omega_h}\phi u
 \le \int_{\Omega_h}u\big(a_j\Delta\phi+\nabla a_j\cdot\nabla\phi
                       -V_j\cdot\nabla\phi\big)
 \le C(1+s^{-1/2})\int_{\Omega_h}\phi u ,\qquad 0<s\le1.
\]
The only extra wall term is $-\int_{\partial \Omega_h}a_j u\,\partial_n\phi\le0$. Hence the terminal weighted mass is at most $C\delta_0$. Bound \eqref{eq:finite-gradient-bound} makes the rescaled gradients uniformly bounded at the top. The volume and path bounds in the weight lemma then imply $\|u(t_*)\|_\infty\le C\delta_0^{1/(d+1)}$. At $t_*\ge1$ the same conclusion follows from the conserved minority mass and the unscaled Lipschitz bound. This excludes a minority exit if $\delta_0$ is sufficiently small relative to $\eta$.

For clarity, the weighted Fisher bound used to exclude the other exits follows directly from the physical entropy $h(r)=\sum_i m_i r_i(\log r_i-1)$. Put $d_r=M_d\nabla r-y\nabla C$ and use the weighted norm $|v|_{y^{-1}}^2=\sum_i|v_i|^2/y_i$. The resistance quadratic form on $\sum v_i=0$ is
\[
   \langle v,B_y^f v\rangle_{y^{-1}}
    =\sum_{i<j} f_{ij}y_i y_j|v_i/y_i-v_j/y_j|^2 .
\]
Its bounds in this norm are the extremal frictions. The entropy production has the exact form $\lambda|\nabla C|^2+(\nu/\rho)\langle d_r,(B_y^f|_E)^{-1}d_r\rangle_{y^{-1}}$ and controls $\sum_i|\nabla r_i|^2/r_i$ on the outer band, since $|d_r|_{y^{-1}}^2+|\nabla C|^2=\sum_i m_i^2|\nabla r_i|^2/y_i$. Its entropy flux has size at most a constant times the square root of production, by the boundedness of $\sum_i y_i|\log r_i|^2$ and $\sum_i r_i|\log r_i|$ on that band. Consequently, for a putative exit $t_*\le T_0$, rescale by $h=\min(1,\sqrt{t_*})$ and use the weight above to obtain
\begin{equation}\label{eq:finite-weight-fisher}
              \int_0^{t_*/h^2}\int_{\Omega_h}\phi
                  \sum_i\frac{|\nabla r_i|^2}{r_i}\le C(T_0).
\end{equation}
In this estimate the derivatives are in rescaled coordinates. Write $\bar r=\langle r^0\rangle$. For each species, apply the scalar logarithmic Sobolev inequality on the bounded connected domain to $\sqrt{r_i}$ to obtain
\[
 \int_\Omega r_i\log(r_i/\bar r_i)\le C_\Omega\int_\Omega|\nabla\sqrt{r_i}|^2 .
\]
The constant is independent of $\bar r_i$, by homogeneity. The entropy production therefore gives exponential decay of
\[
 \mathcal E(t)=\int_\Omega h(r(t))-|\Omega|h(\bar r)
\]
on every surviving strip. The band gives a uniform bound for $\mathcal E(0)$ and a uniform quadratic lower bound for $\mathcal E(t)$ in terms of $\|r(t)-\bar r\|_2$. By \eqref{eq:finite-gradient-bound}, the late-time gradients are uniformly bounded. The volume and path bounds in Lemma~\ref{lem:neumann-weight} turn this $L^2$ decay into uniform convergence to $\bar r$ along surviving strips. The mean total lies in the initial band in parts 1 and 2. In part 3 the mean projected pair lies in the convex inner lens, so its two upper characteristic levels have a strict margin from the outer levels. Positivity excludes exit through the shared lower ray. Thus a density or upper lens exit must occur by a fixed $T_0$.

First fix provisional outer thresholds giving the uniform gradient bound. The entropy decay and the fixed inner-to-outer margins then determine $T_0$ uniformly for all smaller thresholds. With $T_0$ fixed, decrease $\eta$ to control the minority contribution to a density or lens exit and finally choose $\delta_0$ small relative to $\eta$ for the minority bound and the deviations of the two abundant inverse masses.
For a density exit the total flux is $\lambda\rho(m_0\nabla\rho+F)$, where $F=\sum_i(m_i-m_0)\nabla r_i$ and $m_0=m_i$ for the dominant species in part 1. On the rescaled exit cylinder, \eqref{eq:finite-weight-fisher} gives $\iint\phi|F|^2\le C(T_0)(\eta+\delta_0^2)$. For a minority row, use $|\nabla r_i|^2\le\eta|\nabla r_i|^2/r_i$ and the uniform upper bound for $|m_i-m_0|$. For the abundant rows use $|m_i-m_0|\le\delta_0$ and the upper density bound. Test the upper and lower excesses about $U_0$ and $L_0$ with the excess times $\phi$. The squared excess starts at zero. Young's inequality and Gronwall give terminal weighted squared excess at most $C(T_0)(\eta+\delta_0^2)$. At an exit the fixed height of the outer threshold and the scaled gradient bound give a uniform lower bound for the same weighted excess by the volume and path estimates. Choosing $\eta$ and then $\delta_0$ small excludes this exit.

For a lens exit, in binary time the full projected flux has the form $A_{\mathrm b}(U)\nabla U+E_{\mathrm b}$, with
\[
 |E_{\mathrm b}|\le C(r_3|\nabla r|+|\nabla r_3|),
 \qquad \iint\phi|E_{\mathrm b}|^2\le C(T_0)\eta .
\]
This is smooth division at $r_3=0$ in the coefficients of $\nabla r_1,\nabla r_2$, followed by \eqref{eq:finite-weight-fisher} for $\nabla r_3$. The fluxes have vanishing normal component because $\partial_n r=0$. Apply the two scalar entropy tests \eqref{eq:binary-scalar-entropy} at levels $u_\pm(P_i)$. The differentiated weight in the principal term can be absorbed by the production and $C\int\phi H$, using $|\psi'|^2/\psi''\le C\psi$. The remaining errors are bounded by $C\int\phi(1+|\nabla U|)|E_{\mathrm b}|$, since the first and second derivatives of $H$ are bounded. Their time integral tends to zero with $\eta$ by \eqref{eq:finite-weight-fisher}. Thus the weighted terminal entropy tends to zero from zero initial entropy. At an exit at least one of the two entropies has a fixed positive value. Again gradient and volume bounds give a uniform lower bound on its weighted integral, a contradiction. All exits are excluded, completing the proof.

\paragraph{Admissible ranges in the binary region.}
For any fixed pair total $\rho_*>0$ and any fixed fraction bound $0<\theta_0<1$, take the inner apex at $X_*=k\rho_*$ with ratio above $a+(b-a)(1-\theta_0)$. All profiles with pair density in a sufficiently small interval about $\rho_*$ and $r_2/(r_1+r_2)\le1-\theta_0$ then lie in that lens. The interval need not be narrow when the density avoids a particular level. On the upper ray the two slopes are $b$ and $a/X$. Thus any fixed compact positive $X$ interval strictly above $a/b$ is covered by the left plus arc from an apex near the upper ray at its right endpoint. The right minus arc from the left endpoint provides the corresponding result for an interval strictly below $a/b$. This follows by continuous dependence of the characteristic ordinary differential equation, since the discriminant is uniformly positive along the specified upper-ray interval. Any fixed upper bound for $Y/X$ strictly below $b$ is then admitted.

\subsection{Initial mean oscillation and perturbations of a global reference}
\label{sec:finite-bmo}

Use the seminorm of Section~\ref{sec:introduction}. For large radii its averages include the full domain.

\begin{theorem}[Small initial mean oscillation]\label{thm:bmo-global}
Fix positive inverse masses and frictions. If $0<L_0\le \rho^0\le U_0<\infty$, then sufficiently small $[r^0]_{\mathrm{BMO}_2(\Omega)}$ gives a global solution. The threshold is independent of the initial gradient and of every individual species minimum.

More generally, let $R$ be a fixed positive global reference solution whose trace at some finite time belongs to an open subset of $X_\gamma$ consisting of data with global classical solutions. Then the same initial density band and sufficient smallness of
\[
           [r^0-R^0]_{\mathrm{BMO}_2(\Omega)}
                      +|\langle r^0-R^0\rangle|
\]
give a global solution. The threshold can also depend on the reference.
\end{theorem}

\begin{proof}
Let $q$ be the average of $r^0$ on $\Omega\cap B_h(x_*)$. For the weight of Lemma~\ref{lem:neumann-weight} we have
\[
 \int_{\Omega_h}\phi\,|r^0(x_*+hX)-q|^2\,\dd X
                    \le C[r^0]_{\mathrm{BMO}_2(\Omega)}^2.
\]
Indeed the truncated-ball measures are uniformly doubling and the difference of each two successive dyadic averages is bounded by the seminorm times the square root of the doubling constant. Summing these bounds and the corresponding variances against the exponential tails proves the estimate.

Choose an outer total-density band strictly containing the initial band. By Lemma~\ref{lem:shifted-entropy}, relative shifted entropy about any constant state in this band is comparable to the square deviation, and controls the rescaled weighted test:
\[
 \frac{\dd}{\dd s}\int_{\Omega_h}\phi H_\gamma(r\mid q)
                     \le C\int_{\Omega_h}\phi H_\gamma(r\mid q).
\]
Bootstrap the outer band and $\min(1,\sqrt t)\|\nabla r(t)\|_\infty\le1$. The strong trace starts the bootstrap, and the bounded-gradient trace conclusion of Lemma~\ref{lem:wall-record} allows us to include any finite first exit time. Consider a sequence of such exits with initial seminorm tending to zero. At an exit $t_*\le1$, choose an exit point $x_*$, rescale by $h=\sqrt{t_*}$, and use the weighted estimate above. The weighted variance about the corresponding local initial mean tends to zero through the upper face, and the rescaled gradients on the final half-interval are uniformly bounded. The finite-cylinder conclusion of Lemma~\ref{lem:wall-record} gives convergence of the states and gradients through that face, in the interior or in a reflected chart, both at vanishing scales and at scales bounded away from zero. The local limit is constant by the variance bound and the uniform local lower bound of the weight. The local initial means have their totals in $[L_0,U_0]$, so neither a density-band exit nor a unit normalized gradient can persist. At an exit $t_*\ge1$, the gradients are uniformly bounded on the last half time unit, and the unweighted relative entropy about $\langle r^0\rangle$ is nonincreasing and initially tends to zero. The fixed-scale compactness yields the same contradiction. Thus the bootstrap persists for all time.

For the second assertion, the same early rescaling gives a common initial strip of gradient smoothing and total-density confinement. Indeed, for contrary first exits $t_j\downarrow0$ take the rescaled weighted variance about the actual initial localized mean. The contribution of $r^0-R^0$ is small by the preceding dyadic bound. The contribution of the fixed reference tends to zero by uniform continuity and the fixed tails of the weights. Thus the argument leading to a contradiction applies to the early rescalings.

On that common strip choose the shifted entropy on a compact positive total-density band containing the ranges of both solutions. Write $e=r-R$, $F_r=A(r)\nabla r$, $F_R=A(R)\nabla R$ and $M_R=D^2H_\gamma(R)$, where $A$ is the full density diffusion matrix. Differentiating the relative entropy and integrating with the no-flux condition gives the entropy flux
\[
 (D H_\gamma(r)-D H_\gamma(R))F_r-(M_R e)F_R
\]
and the production integrand
\[
 -\big(D^2H_\gamma(r)\nabla r-M_R\nabla R\big)\cdot F_r
 +\big(M_R\nabla e+D^3H_\gamma(R)[\nabla R,e]\big)\cdot F_R .
\]
Write $\nabla r=\nabla e+\nabla R$ and expand
\[
 D^2H_\gamma(r)-M_R=D^3H_\gamma(R)[e]+O(|e|^2).
\]
The terms linear in $e$ and quadratic in derivatives of $R$ cancel by the symmetry of $D^3H_\gamma$ and $A(r)-A(R)=O(|e|)$. The mixed terms are bounded by $C|e||\nabla R||\nabla e|$ and the remaining reference-gradient terms by $C|e|^2|\nabla R|^2$. Coercivity of $D^2H_\gamma(r)A(r)$ in the $\nabla e$ term and Young's inequality give
\[
 \frac{\dd}{\dd t}\int_\Omega H_\gamma(r\mid R)
         +c\int_\Omega|\nabla e|^2
 \le C\|\nabla R\|_\infty^2\int_\Omega H_\gamma(r\mid R).
\]
The assumption makes the initial $L^2$ difference small, and $\|\nabla R\|_\infty$ is bounded on the common strip by the strong trace. Fix a positive comparison time in that strip, earlier than the basin entry time if entry occurs at a positive time. Gradient smoothing and the fixed-scale parabolic estimates on the preceding buffered interval give compactness in $X_\gamma$ using a trace exponent larger than $p$. The relative entropy estimate and this compactness yield closeness in $X_\gamma$ to the trace of $R$ at the comparison time. If entry occurs at time zero, choose the comparison time small enough that the reference remains in the open basin by trace continuity. Otherwise finite-time dependence gives closeness at the entry time. Openness of the basin proves the claim.
\end{proof}

Small mean oscillation does permit a fixed pointwise excursion. For example, for an interior ball $B_{R_*}(x_*)$ take smoothed versions of
\[
 \ell_j(x)=\min\left\{1,\frac1j\log_+\!\frac{R_*}{|x-x_*|}\right\},
\]
extended by zero. The seminorm of $\log|x-x_*|$ on Euclidean balls is bounded by scaling for balls near $x_*$ and by its gradient estimate on balls far from $x_*$. Truncation preserves the bound up to a fixed factor. Truncated-ball volume bounds transfer it to $\Omega$, so $[\ell_j]_{\mathrm{BMO}_2(\Omega)}\le C/j$. Smooth at a radius much smaller than $R_*e^{-j}$ to retain the unit plateau and interior support. For any two fixed positive vectors $v_0,v_1$ in the initial total-density band, $r^0=(1-\ell_j)v_0+\ell_j v_1$ is compatible since it is constant near the wall. Its seminorm tends to zero while its pointwise excursion from $v_0$ is $|v_1-v_0|$.

\section{Affine pressure coordinates and continuation}\label{sec:full-pressure}
For molar composition \(z\), set
\[
 S(z)=\sum_i z_i/m_i,\qquad
 y(z)=\frac{M_d^{-1}z}{S(z)},\qquad
 \omega(z)=m^\top L_{y(z)}.
\]
Along a solution \(S(z)=\rho/C\). The one-form \(\omega\) is smooth
on the closed simplex.
Throughout this section the physical frictions are positive. We allow coincidences among the effective frictions
\(g_{ij}=f_{ij}/(m_i m_j)\). The inverse masses are positive and not
all equal. Summing the full species balances with weights \(1\) and
\(m_i\), respectively, gives
\begin{equation}\label{eq:pressure-pair}
 \begin{aligned}
  \rho_t&=\operatorname{div}(\lambda\rho\nabla C),\qquad
  C_t=\operatorname{div}F_C,\\
  F_C&=(\lambda C+\nu k)\nabla C+F,\qquad
  F=\nu C\,\omega(z)[\nabla z],\qquad k(z)=\omega(z)[z-y(z)]\ge0 .
 \end{aligned}
\end{equation}
Indeed \(M_d\nabla r-y\nabla C=C\nabla z+(z-y)\nabla C\). For positive
states, \(B_y^f|_E\) is positive self-adjoint in the metric
\(\langle v,w\rangle_{y^{-1}}=\sum_i v_iw_i/y_i\). Since
\(\omega(z)[z-y]=(C/\rho)\langle z-y,L_y(z-y)\rangle_{y^{-1}}\), the asserted
sign extends to the closed simplex by continuity. In particular
\(\lambda C+\nu k\) is uniformly positive on each compact positive
total-density band.

\subsection{The pressure integrating factor}
Let \(T(z)=S(z)(B_z^g|_E)^{-1}\). Substitution gives
\(B_y^f M_d^{-1}=S^{-1}B_z^g\). Since \(B_y^f y=0\), this implies
\begin{equation}\label{eq:pressure-projection}
 (M_d-z\otimes m)L_y=T(z),\qquad DS(z)\,T(z)=-S(z)\omega(z).
\end{equation}
For instance, if \(\xi\in E\), the vector \((M_d-z\otimes m)L_y\xi\)
lies in \(E\), and applying \(B_z^g/S\) gives \(\xi\).
Multiplying the first identity by \((1/m_i)_i^\top\) proves the second.

Suppose first that the opposite effective frictions \(\ell_i=g_{jk}\),
\(\{i,j,k\}=\{1,2,3\}\), are distinct. Define
\[
 p_z(X)=\sum_i z_i\!\prod_{j\ne i}(X-\ell_j)
       =X^2-\Lambda X+\Pi,\qquad
 S=s_0+s_1\Lambda+s_2\Pi.
\]
The change of coordinates \(z\mapsto(\Lambda,\Pi)\) is affine and
invertible on the composition plane. Indeed the three polynomials
multiplying the \(z_i\) are distinguished by evaluating at the
\(\ell_i\). Both \(\Lambda\) and \(\Pi\) have positive lower bounds
on the closed simplex. To obtain the matrix in these coordinates, put
\(H_g=\tfrac12\sum_i\ell_i\), \(\widehat h_i=H_g-\ell_i\), and
\(U_j=\sum_i z_i \widehat h_i^j\). Then \(g_{ij}=\widehat h_i+\widehat h_j\),
\(U_1=\Lambda-H_g\), \(U_2=\Pi+H_g^2-\sum_{i<j}\ell_i\ell_j\), and on \(E\)
\[
 B_z^g=U_1 I+(I-z\mathbf1^\top)\operatorname{diag}(\widehat h_i).
\]
Applying \(D\Lambda=\widehat h^\top\) and \(D\Pi=(\widehat h_i^2)_i^\top\) on \(E\) to
this formula, and using the cubic satisfied by the \(\widehat h_i\), gives the
matrix \(\bigl(\begin{smallmatrix}0&1\\-\Pi&\Lambda\end{smallmatrix}\bigr)\).
Consequently, \eqref{eq:pressure-projection} gives
\begin{equation}\label{eq:pressure-moments}
 T=\frac S\Pi\begin{pmatrix}\Lambda&-1\\ \Pi&0\end{pmatrix},
 \qquad
 \omega=-\left(s_2+\frac{s_1\Lambda}{\Pi}\right)d\Lambda
              +\frac{s_1}{\Pi}\,d\Pi .
\end{equation}
These formulas hold on the closed simplex, including where the roots
of \(p_z\) coalesce. Thus the affine coordinates avoid any root-gap
condition. If \(s_1\ne0\), set \(\chi=s_2/s_1\) and take
\begin{equation}
 \Phi(\Lambda,\Pi)=e^{-\chi\Lambda}\Pi-\int_0^\Lambda \xi e^{-\chi\xi}\,\dd\xi,\qquad
 \omega=\alpha\,d\Phi,\qquad
 \alpha=s_1e^{\chi\Lambda}/\Pi .
\end{equation}
If \(s_1=0\), take \(\Phi=-s_2\Lambda\) and \(\alpha=1\).
In both cases \(\alpha\) has constant sign and its absolute value
is bounded above and away from zero.

For a coincident effective pair, relabel so that
\(g_{12}=g_{13}=a\), \(g_{23}=b\ne a\). Set \(v=z_2+z_3\), \(w=z_2\),
\(q=a+(b-a)v\), and write \(S=S_0+ev+d_0w\) in these coordinates.
A direct calculation of the two rows on \(E\) gives
\[
 T=S\begin{pmatrix}1/a&0\\(b-a)w/(a q)&1/q\end{pmatrix},
 \qquad
 \omega=-\left(\frac e a+\frac{d_0(b-a)w}{a q}\right)dv-\frac{d_0}{q}\,dw .
\]
If \(d_0\ne0\), put \(\ell=(b-a)/a\). The potential
\[
 \Phi(v,w)=e^{\ell v}w+(e/d_0)\int_0^v e^{\ell\xi}q(\xi)/a\,\dd\xi
\]
satisfies \(d\Phi=-(q/d_0)e^{\ell v}\omega\). If \(d_0=0\) then
\(\omega=-(e/a)dv\) is exact. Finally, if all effective frictions
equal \(a\), then \(T=(S/a)I\) and \(\omega=-dS/a\). In each case
\(\omega\) is a bounded signed multiplier times an exact form,
with a multiplier bounded away from zero in absolute value.
These coincidences concern the matrices in molar coordinates.

\subsection{Weak ancient rigidity in the positive total-density band}
\begin{lemma}[Pressure rigidity]\label{lem:full-pressure-rigidity}
Let \(r\) be a bounded, nonnegative ancient local-energy solution of the full
equations on \(\mathbb R^d\times(-\infty,0)\) with total density in a fixed
compact subinterval of \((0,\infty)\). Assume the uniform local energy bounds
\(\iint_{Q_R}|\nabla r|^2\le C_{\mathrm{est}}R^d\) on backward cylinders, as in
the proof of Lemma~\ref{lem:one-cylinder}. If \(C\), \(\rho\), \(S\), or \(C-G(\rho)\)
is spatially homogeneous on almost every time slice, then \(r\) is constant.
Here \(G\) is \(C^1\) on the density band with \(G'\) uniformly positive.
The conclusion also holds if \(\operatorname{div}(\omega(z)[\nabla z])=0\)
or if \(\operatorname{div} F=0\) in distributions.
\end{lemma}

\begin{proof}
First assume that \(C\) is constant in space and time. From
\eqref{eq:pressure-pair}, \(\rho_t=0\) and
\(\operatorname{div}(\omega(z)[\nabla z])=0\). Since the flux belongs to local
\(L^2\), testing with products of spatial and temporal tests gives its
elliptic equation on almost every slice. A countable dense family of
spatial tests on an exhaustion by balls gives one set of full time measure
on which all spatial tests are admissible. Thus the integrating factor
gives a scalar homogeneous elliptic equation for the bounded function
\(\Phi(z)\) with uniformly elliptic coefficient \(\lvert\alpha(z)\rvert\).
Elliptic oscillation decay on expanding balls makes \(\Phi(z)\) spatially
constant. Thus \(\omega(z)[\nabla z]=0\) and
the projected compositional balance is
\[
 z_t=\nu\operatorname{div}(T\nabla z),\qquad S_t=0.
\]
This follows directly by multiplying the species conservation law by
\(M_d/C\). Here \(C\) is constant and \(M_d L_y\nabla z=T\nabla z\)
because \(m^\top L_y\nabla z=0\).
Only spatial Sobolev chain rules are used.

In the distinct-friction case with \(s_1\ne0\) the zero force gives
\(\nabla\Pi=(\Lambda+\chi\Pi)\nabla\Lambda\), so the two moment
balances are
\[
 \Lambda_t=-\nu \chi\,\operatorname{div}(S\nabla\Lambda),\qquad
 \Pi_t=\nu\operatorname{div}(S\nabla\Lambda).
\]
For \(\chi<0\) bounded ancient scalar parabolic oscillation decay gives
constant \(\Lambda\). The force relation and the second equation
then give constant \(\Pi\). For \(\chi\ge0\),
\[
 \nabla S=s_1\bigl(1+\chi(S-s_0)/s_1\bigr)\nabla\Lambda
          =s_1(1+\chi\Lambda+\chi^2\Pi)\nabla\Lambda .
\]
The coefficient in the middle is a fixed bounded function of \(S\)
and is separated from zero, by the last expression. Hence
\(\nabla\Lambda\) is independent of time. The second compositional balance
has a time-independent right hand side. Testing against a spatial
test function and integrating over arbitrarily long past intervals
shows that this right hand side vanishes, since \(\Pi\) is bounded.
Elliptic oscillation decay for \(\operatorname{div}(S\nabla\Lambda)=0\)
and the two balances now make both moments constant.
If \(s_1=0\), then \(s_2\ne0\), and the zero force gives
\(\nabla\Lambda=0\). The quantity \(\Pi\) is frozen because
\(S=s_0+s_2\Pi\). The first moment equation reads
\(\Lambda_t=-\nu\operatorname{div}((S/\Pi)\nabla\Pi)\).
The same long-time test and elliptic oscillation argument give
constant moments. If the effective pair coincides, the first
coordinate \(v\) solves a scalar uniformly parabolic equation with
coefficient \(\nu S/a\). After \(v\) is constant, \(w\) solves one
with coefficient \(\nu S/q\). The all-coincident case is scalar in
each coordinate. Thus the rigidity holds for constant pressure in all cases.

A spatially homogeneous conserved affine moment \(U\) is independent of
time. To see this, test its conservation law between two times against
\(\eta_R=R^{-d}\eta(\cdot/R)\), where \(\eta\in C_c^\infty(B_1)\),
\(\eta\ge0\), and \(\int\eta=1\). For a fixed interval \([a,b]\), the
increment in the average is at most
\[
 C_{\mathrm{est}}(b-a)^{1/2}\|\nabla\eta_R\|_2
     \|\nabla r\|_{L^2(B_R\times(a,b))}
 \le C_{\mathrm{est}}(b-a)^{1/2} R^{-1}
\]
for \(R\) large enough, by using a backward cylinder containing the interval.
The average is the value of \(U\) on almost every slice, and the bound tends
to zero. In particular homogeneous \(C\) is constant. If \(\rho\)
is homogeneous, \(\rho\) is constant and its balance makes \(C\)
bounded harmonic on each slice. If \(S\) or \(C-G(\rho)\) is
spatially homogeneous, the mass equation is scalar uniformly
parabolic for \(\rho\), with coefficient \(\lambda\rho/S(t)\) or
\(\lambda\rho G'(\rho)\), respectively. Bounded ancient scalar
regularity makes \(\rho\) constant, reducing to the previous case.

If \(\operatorname{div}(\omega(z)[\nabla z])=0\), the first elliptic
argument directly makes \(\omega(z)[\nabla z]=0\). If
\(\operatorname{div}F=0\), the multiplier in that argument is
\(\nu C\lvert\alpha\rvert\), also uniformly elliptic. In either case
\(C\) solves a homogeneous uniformly scalar parabolic equation with
coefficient \(\lambda C+\nu k\). It is constant by ancient oscillation
decay, so the result follows.
\end{proof}

\subsection{Full-system tests on a cylinder and at a terminal time}
Set \(D_R=\Omega\cap B_R(x_0)\) and
\(Q_R^\Omega=D_R\times(t_0-R^2,t_0)\), allowing \(x_0\in\partial\Omega\).
We use the mixed \(H^{-1}\) norm \(\|\cdot\|_{-1;D_R}\) dual to the
spatial gradient seminorm on tests vanishing at the artificial
boundary, as in Section~\ref{sec:boundary-estimates}.
Fix a total-density band \(I=[c_-,c_+]\subset(0,\infty)\) and
\(0<b_-<b_+<\infty\). Put
\[
 \mathcal P_b(Q_R^\Omega)=R^{-d}\int_{t_0-R^2}^{t_0}
       \|\operatorname{div}(F_C-b\nabla C)\|_{-1;D_R}^2\,\dd t,
 \qquad
 \mathcal V_U(Q_R^\Omega)
  =\frac1{|Q_R^\Omega|}\inf_{c(t)}
    \iint_{Q_R^\Omega}|U-c(t)|^2.
\]
Measurable conductivities \(b\) range over \([b_-,b_+]\). The divergence
in the norm is the variational action of the flux, with no boundary
term at the physical wall.

\begin{theorem}[Pressure continuation tests]\label{thm:full-pressure-continuation}
Fix the coefficients, \(I\), and the conductivity bounds \(b_-,b_+\).
There exist \(R_0,\epsilon_0,K>0\)
such that the following assertions hold through the top of \(Q_R^\Omega\),
\(0<R\le R_0\), on any strong branch with \(\rho\in I\). No individual
species lower bound is needed.
\begin{enumerate}
\item If \(\inf_b\mathcal P_b(Q_R^\Omega)<\epsilon_0\), then
\(R\sup_{Q_{R/2}^\Omega}|\nabla r|\le K\).
The same bound holds if \(\mathcal V_U(Q_R^\Omega)<\epsilon_0\)
for \(U=C,\rho,S\), or \(C-G(\rho)\), with \(G\) as in
Lemma~\ref{lem:full-pressure-rigidity}. Constants for a variance
test may depend on the chosen \(U\).
\item In \(d=2\) the same bound holds if, for \(g=\partial_t\rho\)
or \(g=\operatorname{div}F\),
\[
 R^{-2}\int_{t_0-R^2}^{t_0}
       \min\left\{\int_{D_R}g_+,\int_{D_R}g_-\right\}\,\dd t
     <\epsilon_0 .
\]
\item A branch on a finite interval \([0,T)\) with strictly positive
compatible initial data as in Definition~\ref{def:classical} continues strongly through
\(T\) if
\[
 \Theta=\nu\operatorname{div}(\omega(z)[\nabla z]),\qquad
             \int_0^T\|\Theta(t)\|_{L^\infty(\Omega)}\,\dd t<\infty .
\]
This last condition supplies a positive total-density band, so the band
need not be assumed separately.
\end{enumerate}
\end{theorem}

\begin{proof}
Apply Lemmas~\ref{lem:small-excess}--\ref{lem:one-cylinder} to the
full diffusion matrix on \(I\). After normalizing the outer radius to one,
let the defect of a contrary sequence be \(\delta_j\to0\) and the pole
gradient be \(M_j\to\infty\). For the flux defect take conductivities with
\(\mathcal P_{b_j}\) within \(o(1)\) of the infimum and use their actual defect
values in the choice of scale. The intermediate scales are
\(h_j=\max(M_j^{-1/2},\delta_j^{1/(2\beta_{\mathrm{sc}})})\), with \(\beta_{\mathrm{sc}}=d\)
for the flux defect and \(\beta_{\mathrm{sc}}=d+2\) for the state variance.
On each fixed rescaled cylinder the defect is \(O(h_j^{-\beta_{\mathrm{sc}}}\delta_j)\).
At a wall record the mixed dual norm controls the reflected source by
testing against the sum of the two half-chart restrictions. For the
variance defect, restriction to the smaller cylinders and strong local
\(L^2\) convergence show that the observed scalar in the limit is
spatially homogeneous. Reflection preserves this property.

For the flux defect, use the homogeneous scalar comparison of
Lemma~\ref{lem:one-cylinder} on each fixed rescaled backward cylinder.
It has conductivity \(b_j\) and the same parabolic boundary values as
the rescaled pressure. In a reflected chart, use the associated
elliptic tensor and the positive, time-independent volume density.
The zero-trace energy inequality controls the squared difference in
\(L^\infty_tL^2_x\) and the integrated squared gradient difference by
the squared residual norm. The comparison function is uniformly bounded
by the maximum principle and obeys scalar oscillation decay with
constants uniform over \(b_j\). Letting the comparison cylinder radius
tend to infinity after passage to the limit gives constant \(C\) in the
ancient limit. This argument does not require convergence of \(b_j\).
Since \(h_jM_j\to\infty\), small-excess regularity ensures nonzero
full-state variance on a fixed cylinder, which persists by local
entropy compactness. In both cases Lemma~\ref{lem:full-pressure-rigidity}
makes the limiting solution constant, a contradiction.

For the planar signed-source bound, write
\(g=\operatorname{div}(a(x,t)\nabla V)\) with \(V\) bounded,
\(a\) uniformly positive, and zero normal flux on the wall.
For \(\partial_t\rho\) take \(V=C\) and \(a=\lambda\rho\).
For \(\operatorname{div}F\) take \(V=\operatorname{sgn}(\alpha)\Phi\)
and \(a=\nu C|\alpha|\), since the sign of \(\alpha\) is fixed.
On almost every time slice, test in \(B_{R_{\log}}\cap\Omega\)
with \(\eta_{\log}^2(V-\inf V)\) and \(\eta_{\log}^2(\sup V-V)\), where
\(\eta_{\log}\) is the logarithmic cutoff from radius \(R_{\log}\) to \(r<R_{\log}\).
The two tests respectively bound the source term by the integrals of
\(g_-\) and \(g_+\), up to the uniform oscillation bound for \(V\).
Absorbing the cutoff cross-terms by Young's inequality gives
\[
 \int_{\Omega\cap B_r}|\nabla V|^2
 \le \frac{C_{\mathrm{est}}}{\log(R_{\log}/r)}
      +C_{\mathrm{est}}\min\left\{\int_{\Omega\cap B_{R_{\log}}}g_+,
                    \int_{\Omega\cap B_{R_{\log}}}g_-\right\}.
\]
The cutoff cost is bounded by its whole-plane value
\(2\pi/\log(R_{\log}/r)\). Both source integrals are monotone under spatial
restriction, so the smaller integral is bounded by the corresponding
minimum on the original cap. In a contradictory
sequence normalize \(R=1\). Write \(\delta_j\to0\) for the signed-source
defects and \(M_j\to\infty\) for the gradients at the selected points.
Take \(R_{\log}=1/4\) about each point and
\[
 h_j=\max\{M_j^{-1/2},\,\delta_j^{1/4}\}.
\]
For any fixed \(L>0\), use inner radius \(Lh_j\) and integrate over
the corresponding backward time interval. The squared gradient
integral of the rescaled \(V\) is bounded by
\[
 \frac{C_{\mathrm{est},L}}{\log(R_{\log}/(Lh_j))}+C_{\mathrm{est}}h_j^{-2}\delta_j\longrightarrow0 .
\]
Meanwhile \(h_jM_j\to\infty\), so the small-excess estimate
and local entropy compactness give a nonconstant ancient limit.
The state maps \(C\) and \(\Phi(z)\) converge strongly in local \(L^2\).
Their scaled spatial gradients tend to zero in the respective cases,
so their limits have zero spatial weak gradient. For \(V=C\) this gives
spatially constant pressure. For \(V=\Phi\) it gives zero microscopic
flux \(F=\nu C\alpha(z)\nabla\Phi(z)\) in the limit. Reflection gives the same conclusions at the
wall. Both cases contradict Lemma~\ref{lem:full-pressure-rigidity}.

For the time-integrated criterion use \eqref{eq:pressure-pair}
in the form
\[
 C_t=\operatorname{div}((\lambda C+\nu k)\nabla C)
     +\nu\omega(z)[\nabla z]\cdot\nabla C+C\Theta,\qquad
 \partial_n C=0 .
\]
Apply scalar Neumann comparison first from a positive time and then
let that time decrease to zero using continuity of the strong trace.
The resulting bounds are
\[
 (\min_{\overline\Omega} C(0))\exp\!\left[-\int_0^t\|\Theta_-(\tau)\|_\infty\,\dd\tau\right]
 \le C(t)\le
 (\max_{\overline\Omega} C(0))\exp\!\left[\int_0^t\|\Theta_+(\tau)\|_\infty\,\dd\tau\right].
\]
The same bounds for \(\rho\) up to fixed factors follow from positivity
of the inverse masses. If gradients were unbounded near \(T\), take
gradient records as in Lemma~\ref{lem:wall-record}. For a dilation
around a record with spatial scale \(h\), the rescaled source is
\(\Theta_h(x,s)=h^2\Theta(x_*+hx,t_*+h^2s)\). Its \(L^1_sL^\infty_x\)
norm on any fixed rescaled time interval is bounded by the original
integral of \(\|\Theta(t)\|_\infty\) over a shrinking interval, and hence
tends to zero. In a reflected chart the weak source
is reflected with the volume density, which is uniformly bounded
and converges to a constant. The microscopic flux has zero normal
component and thus produces no interface term. The limiting compositional
one-form has vanishing divergence, by strong state convergence and
weak gradient convergence. Lemma~\ref{lem:full-pressure-rigidity} then
excludes a nonzero gradient at the record. Bounded gradients and the
fixed band give the restart.
\end{proof}

Assume the branch has a common positive total-density band on a
terminal strip \([t_1,T)\). Call \(x_0\in\overline\Omega\) singular at
\(T\) if the gradient is unbounded in every backward neighborhood of
\((x_0,T)\). The integrated shifted entropy identity on \([t_1,T)\), with zero normal flux,
bounds \(\iint|\nabla r|^2\) there. The flux \(F\) is a linear function of
\(\nabla r\) with smooth bounded coefficients on the band, so \(F\in L^2\)
on that strip.

Choose \(b_-,b_+\) to bracket \(\lambda C+\nu k\) on the band, including
its species faces. For \(T-R^2>t_1\), apply the first part of
Theorem~\ref{thm:full-pressure-continuation} to the cylinders with poles
\((x_0,T-\tau)\), \(\tau>0\), and take \(b=\lambda C+\nu k\). The estimate
\(\|\operatorname{div}F\|_{-1;D_R}\le\|F\|_{L^2(D_R)}\) bounds the defect.
If the integral of \(|F|^2\) on the cylinder with pole \((x_0,T)\)
were less than \(\epsilon_0R^d\), continuity of the moving time integral
would give the same strict bound for all sufficiently small \(\tau>0\).
The resulting uniform gradient bound on the shifted half-cylinders
would bound the gradient on \(D_{R/2}\times(T-R^2/4,T)\), contrary to
singularity at \((x_0,T)\). Thus
\(\iint_{Q_R^\Omega}|F|^2\ge\epsilon_0R^d\) on every sufficiently small
terminal cylinder centered at a singular point.

A maximal family of disjoint balls of radius \(R\) centered at terminal
singular points has cardinality \(N_R\) with
\[
 \epsilon_0N_RR^d\le\int_{T-R^2}^T\int_\Omega |F|^2\longrightarrow0 .
\]
Their triple enlargements cover the terminal singular set, which therefore
has zero \(d\)-dimensional Hausdorff measure, or equivalently zero Lebesgue
measure.

\subsection{An exact finite-density global family}
Suppose now that the three inverse masses and the \(\ell_i\) are
each pairwise distinct, and relabel so that
\(\ell_1<\ell_2<\ell_3\). The following family gives finite-amplitude
global profiles without restricting the overall density scale.

\begin{theorem}[Invariant heat lines]\label{thm:pressure-heat-line}
A smooth regular invariant curve in the strictly positive density
cone exists if and only if \(1/m_2\) is a strict extremum among
\(1/m_1,1/m_2,1/m_3\). Invariance here means local preservation for
every smooth scalar profile on the curve. If this condition holds,
the maximal curve at each pressure \(C_*>0\) is the positive segment
\[
 C=C_*,\qquad
 p_z(\kappa)=0,\qquad
 \kappa=-s_1/s_2,\qquad
 S_*=s_0-s_2\kappa^2 .
\]
Every compatible profile with range compactly inside that segment
gives an exact global solution, with \(\Lambda\) solving the scalar
Neumann heat equation of diffusivity \(\nu S_*/\kappa\).
Each such reference has an open neighborhood of global positive data
in the compatible strong trace space of the full system.
\end{theorem}

\begin{proof}
Write the full operator as \(r_t=\operatorname{div}(A(r)\nabla r)\).
For a regular curve \(r=R(s)\), the operator applied to \(R(s(x))\)
is
\[
 A(R)R'\Delta s+(A(R)R')'|\nabla s|^2.
\]
The values of \(\Delta s\) and \(|\nabla s|^2\) can be chosen independently at an interior point. Invariance thus gives \(A(R)R'=aR'\) and tangency of \((A(R)R')'=a'R'+aR''\). Normal ellipticity gives \(a>0\), so \(R''\) is tangential. The curve is contained in an affine line of fixed direction \(v\). If \(\sum_i v_i\ne0\), normalize this sum to one and write
\(r=\rho v+w\), \(\sum_i w_i=0\). With \(p=m\cdot v\), summation
of the flux makes \(a=\lambda\rho p>0\) and the microscopic flux
per unit \(\nabla\rho\) equals \(J_0=-\lambda p w/\nu\).
Then
\(B^f_{v+w/\rho}J_0=M_dv-pv-pw/\rho\), while the coefficient
of \(1/\rho\) on the left is zero since \(B_w^fw=0\).
Thus \(w=0\), and \(M_dv=pv\). Positivity would force equal masses,
which are excluded here. Therefore \(\sum_i v_i=0\) and the summed
flux also gives \(m\cdot v=0\).

The curve lies at fixed \(C,\rho\), and its \(z\)-direction is tangent
to \(S=\mathrm{const}\). By \eqref{eq:pressure-projection}, flux tangency
is exactly the eigenvector condition for \(T\) in that direction.
For a line \(A_0\Lambda+B_0\Pi=c\), testing the normal \((A_0,B_0)\)
against \(T(-B_0,A_0)^\top\) in \eqref{eq:pressure-moments} gives
\(-S(A_0^2+B_0c)/\Pi\). If \(B_0=0\), tangency also forces \(A_0=0\), which gives no line.
Otherwise the tangent lines have
\(\Pi-\kappa\Lambda+\kappa^2=0\), \(\kappa=-A_0/B_0\).
On a line meeting the simplex, \(\kappa\) is a positive root of \(p_z\),
and the tangent diffusivity is \(S/\kappa\). Constancy of \(S\) requires
\(s_1+s_2\kappa=0\). If \(s_2=0\) then \(s_1\ne0\), since distinct masses
make \(S\) nonconstant, and there is no such line. Evaluating \(S\) at the
vertices gives
\[
 \frac{1/m_i-1/m_j}{\ell_i-\ell_j}
       =-s_1-s_2\ell_k,\qquad\{i,j,k\}=\{1,2,3\}.
\]
The outer values of the affine function \(s_1+s_2\ell\) have
opposite signs precisely when \(1/m_2\) is a strict extremum.
Distinct masses also make all three vertex differences nonzero, so its
root cannot equal any \(\ell_i\). A root in either interlacing interval
gives both signs for the vertex values of \(p_z(\kappa)\). Their affine
combination therefore vanishes somewhere in the open simplex. If
\(\kappa\) is outside \([\ell_1,\ell_3]\), all three vertex values of
\(p_z(\kappa)\) are positive. This proves
the order condition and identifies the segment.

Along the segment, \(\omega(z)[\nabla z]=0\) and
\(J=(C_*S_*/\kappa)M_d^{-1}\nabla z\); the Darcy flux vanishes.
All species therefore solve the claimed heat equation. It preserves
the compact segment and converges to the positive spatial mean in
the strong trace. At the homogeneous state \(D^2H(\bar r)\) symmetrizes the diffusion
matrix, with all its eigenvalues positive. Each eigenvalue multiplies
a scalar Neumann Laplacian, so the Neumann Poincar\'e inequality
gives a spectral gap on the specieswise mean-zero subspace. The local Neumann
maximal-regularity estimate in a fixed exponentially weighted space
then gives stability. Indeed the principal difference from the
constant matrix and the quadratic gradient terms cost
\(O(\|r-\bar r\|_{C_t W^{2-2/p,p}})\) times the weighted
maximal-regularity norm. Absorb these terms for sufficiently small
\(X_\gamma\) distance to equilibrium, using conservation of the
species means. The estimates vary continuously with nearby positive
means. Continuous dependence up to a fixed time then puts sufficiently
small compatible perturbations of the reference into this neighborhood.
\end{proof}

\paragraph{Initial pressure flux.}
For distinct effective frictions with $s_1s_2\ne0$, normalize $\Phi=e^{-\chi\Lambda}p_z(-1/\chi)$. Choose an interior state $z_*$ with $\phi_0=\Phi(z_*)\ne0$ and write $S_*=S(z_*)$, $\Lambda_*=\Lambda(z_*)$. Take $C^0$ constant and $\Lambda^0=\Lambda_*+\delta u$, with $u\in C_c^\infty(\Omega)$ and $\nabla\Delta u(x_*)\ne0$. Solve $\Phi=\phi_0$ for $\Pi^0$ near the state, since $\Phi_\Pi=e^{-\chi\Lambda}$. At time zero the moment balances give
\[
 C_t(0)=0,\qquad \Phi_t(0)=\nu \chi^2\phi_0\Div(S^0\nabla\Lambda^0).
\]
The gradient of this divergence at $x_*$ is $\delta S_*\nabla\Delta u(x_*)+O(\delta^2)$. For small nonzero $\delta$, interior smoothness through time zero gives a nonzero linear term in $F(x_*,t)$ and $\nabla C(x_*,t)=o(t)$.

\paragraph{The conventional constant-diffusivity normalization.}
The pressure criteria in this section also apply if the microscopic
flux in the species equations is \((\nu/S)J\), with the same
constrained resistance law. Put \(\mu=1/S=C/\rho\), as in Corollary~\ref{cor:constant-d}. Then
\[
 F_C=(\lambda C+\nu\mu k)\nabla C+F,\quad
 F=\nu\mu C\,\omega(z)[\nabla z],\quad
 \Theta=\nu\operatorname{div}(\mu\,\omega(z)[\nabla z]).
\]
The potential \(\Phi\) is unchanged and \(\mu\omega=(\mu\alpha)\,d\Phi\).
For constant ancient \(C\), \(S_t=0\) and the same compositional moment
proof uses \(\mu T\) in place of \(T\). The coefficient in the
time-independent divergence for \(\chi\ge0\) is \(\mu S\), which is
constant for this normalization. For \(\chi<0\) the pressure-constrained
compositional equation remains scalar uniformly parabolic. The no-band criterion and planar signed-source test use the indicated
\(F,\Theta\). The concentration conclusion uses conductivity bounds
bracketing \(\lambda C+\nu\mu k\).
The exact heat lines at fixed \(\rho,C\) are the same, and their
scalar diffusivity becomes \(\nu/\kappa\). To rule out varying-\(\rho\)
invariant lines, repeat the tangent-flux argument with
\(r=\rho v+w\), \(q=m\cdot v>0\), \(D=m\cdot w\), \(\xi=1/\rho\).
The left hand side \(B^f_{v+w \xi}(-\lambda q w/\nu)\) is constant in
\(\xi\), while the right hand side is
\((q+D \xi)(M_dv-qv-qw \xi)\).
If \(w\ne0\) the quadratic coefficient forces \(D=0\), and the
linear term then forces \(w=0\). Normal stability is unchanged.
The binary characteristic-region theorem uses its stated constant normalized resistance.

\section{Composition tests on one frozen cylinder}
\label{sec:frozen-affine-cylinder}

Let \(K=\{y\geq0:\mathbf1^\top y=1\}\). For constant prescribed density, \(y=y(x,\tau)\) solves the following equation in fast time.
\begin{equation}
 y_\tau=\operatorname{div}(A(y)\nabla y),\qquad
 A(y)=(B_y^f|_E)^{-1}P_yM_d|_E .
 \label{eq:frozen-cylinder-operator}
\end{equation}
Write \(I_R=(\tau_0-R^2,\tau_0)\) and \(Q_R=D_R\times I_R\). Use the mixed dual norm of Section~\ref{sec:boundary-estimates} on wall cylinders and \(\dot H^{-1}(B_R)\) in the whole space.

Write \(y_{D_R}(\tau)=|D_R|^{-1}\int_{D_R}y(x,\tau)\,\dd x\) and
\[
 \mathcal C_R=\frac{1}{|Q_R|}
 \int_{Q_R}(y-y_{D_R}(\tau))\otimes(y-y_{D_R}(\tau))\,\dd x\,\dd\tau .
\]
For a unit covector \(\ell\in E\) and a measurable scalar \(b\) with fixed positive bounds \(b_-\leq b\leq b_+\), set
\[
 \mathcal D_{\ell,b}(Q_R)
 = R^{-d}\int_{I_R}
 \left\|\operatorname{div}\big((\ell^\top A(y)-b\ell^\top)\nabla y\big)
                    \right\|_{-1;D_R}^{\,2}\,\dd\tau .
\]
The divergence uses the variational no-flux action at the physical wall.

\begin{theorem}[Frozen one-cylinder tests]\label{thm:frozen-affine-cylinder}
Fix positive masses and frictions. In the bounded-domain case take \(\Omega\) as in Section~\ref{sec:introduction} and take \(R\leq R_0(\Omega)\), with \(R_0(\Omega)>0\) sufficiently small. For a \(K\)-valued strong frozen solution on \(Q_R\), smallness of either
\[
 \lambda_{\min}(\mathcal C_R|_E)
 \qquad\text{or}\qquad
 \inf_{\ell,b}\mathcal D_{\ell,b}(Q_R)
\]
implies \(R\sup_{Q_{R/2}}|\nabla y|\leq C\). Each test has its own positive threshold, independent of individual species minima. For the projected-flux test the threshold and gradient bound may depend on \(b_-,b_+\). If \(m\) is nonconstant, write \(p=m\cdot y\). Each of the following is also a sufficient smallness test:
\begin{equation}
 \begin{gathered}
 R^{-d}\int_{I_R}\|p_\tau\|_{-1;D_R}^{\,2}\,\dd\tau,\qquad
 R^{-d}\int_{Q_R}|p_\tau|\,\dd x\,\dd\tau,\\
 \frac{1}{|Q_R|}\inf_{c=c(x)}\int_{Q_R}|p(x,\tau)-c(x)|^2\,\dd x\,\dd\tau .
 \end{gathered}
 \label{eq:frozen-pt-tests}
\end{equation}
In two dimensions the same conclusion follows from smallness of
\begin{equation}
 R^{-2}\int_{I_R}\min\left\{\int_{D_R}(p_\tau)_+\,\dd x,\
                                  \int_{D_R}(p_\tau)_-\,\dd x\right\}\,\dd\tau .
 \label{eq:frozen-pt-signed}
\end{equation}
The covariance test extends to a fixed smooth positive nonconstant prescribed density having zero normal derivative, at sufficiently small absolute radii. The pressure-time and projected flux tests above concern constant prescribed density.
\end{theorem}

\paragraph{Reduction to ancient equations.}
For the frozen equation the shifted entropy in Lemma~\ref{lem:shifted-entropy} gives, on the tangent plane, the production and relative-flux estimates used by Lemma~\ref{lem:small-excess}. In particular on a bounded ancient limit it gives
\(\int_{Q_R}|\nabla y|^2\leq C R^d\), using a larger cylinder and boundedness of \(K\). The local gradient bound gives weak \(L^2\) compactness, and the conservation law gives a local \(L^2_\tau H^{-1}_x\) bound for the time derivatives. After localization the states converge strongly in \(L^2\). The bounded matrix \(A(y_j)\) then converges strongly in \(L^2\) on fixed cylinders, so its product with \(\nabla y_j\) converges against smooth flux tests. The same argument in a frozen wall chart uses the bounded uniformly positive metric and the time-independent volume density.

Suppose a test fails along a sequence with normalized defect \(\delta_j\to0\) and normalized gradient \(M_j\to\infty\) in the inner half-cylinder. Set the outer radius to one and dilate by \(h_j\) in space and \(h_j^2\) in time about the selected points, where
\[
 h_j=\max\{M_j^{-1/2},\,\delta_j^{1/(2\beta_{\mathrm{sc}})}\},\quad
 \beta_{\mathrm{sc}}=\begin{cases}
 d,&\text{derivative or flux defect},\\
 d+2,&\text{state variance}.
 \end{cases}
\]
The scaled gradients \(h_jM_j\) diverge. Lemma~\ref{lem:small-excess} implies a strictly positive lower bound on the full-state space--time variance on a fixed cylinder about the selected point. If the normalized distance to the wall diverges, use interior compactness. Otherwise flatten the wall, reflect as in Section~\ref{sec:boundary-estimates}, and pass to a half-space whose boundary may lie at a finite distance from the pole. On each fixed reflected cylinder the metric and volume density tend to constants. A fixed linear change gives the Euclidean equation, and the reflected limit solves \eqref{eq:frozen-cylinder-operator} on the whole space. The strong \(L^2\) convergence holds up to the upper time face in space--time, by discarding a last time interval of small length and using boundedness of the states there. For the moving physical caps, the indicators converge almost everywhere and the normalized volumes converge to positive limits. Together with strong \(L^2\) convergence this shows that the variance lower bound persists in the limit. The trials \(c(x)\) and \(c(\tau)\) restrict to every smaller cylinder and reflect in the same charts.

On each fixed normalized cylinder the scaled defect is at most \(C h_j^{-\beta_{\mathrm{sc}}}\delta_j\). For the negative norm this follows by extending tests by zero beyond the artificial boundary and, on a reflected chart, pairing the source with the sum of the two half-chart restrictions. For the projected flux test, compare \(\ell_j\cdot y_j\) to the scalar solution with the same parabolic boundary values and conductivity \(b_j\), in the geometric divergence form of the rescaled equation. The difference has squared \(L^\infty_tL^2_x\) norm controlled by the squared \(L^2_tH^{-1}_x\) norm of the source. The scalar oscillation estimate on expanding cylinders makes the limiting projection constant, as in Lemma~\ref{lem:one-cylinder}. This step does not require convergence of \(b_j\). For the covariance test choose minimizing unit covectors \(\ell_j\) and a subsequence converging to \(\ell\). Spatially constant functions on each fixed cylinder form a closed subspace of \(L^2\), so vanishing distance to this subspace and strong local \(L^2\) convergence give spatial constancy of \(\ell\cdot y\). For prescribed variable density its derivative flux is \(O(R)\) at normalized radius \(R\), as in Lemma~\ref{lem:small-excess}. It vanishes under the intermediate dilation, and the chart metric separately converges to a constant. In \eqref{eq:frozen-pt-tests}, the scaled derivative tends to zero in the indicated dual or \(L^1\) norm for the first two tests. The third gives a vanishing \(L^2\) distance to the closed subspace of time-independent functions on each fixed cylinder. Each of the three therefore gives \(p_\tau=0\) in distributions in the ancient limit. It remains to establish both required rigidity statements.

\paragraph{Spatially constant affine observation.}
For a conserved affine scalar \(\ell\cdot y\), let \(\eta_R(x)=R^{-d}\eta(x/R)\), with \(\eta\ge0\) smooth, supported in the unit ball and of integral one. On any fixed interval \([\tau_1,\tau_2]\) the flux bound gives, as \(R\to\infty\),
\begin{equation}\label{eq:frozen-affine-increment}
 \left|\int_{\R^d}\eta_R\,\ell\cdot
       \big(y(\tau_2)-y(\tau_1)\big)\,\dd x\right|
 \le C \|\nabla\eta_R\|_2(\tau_2-\tau_1)^{1/2} R^{d/2}
 \le \frac{C (\tau_2-\tau_1)^{1/2}}{R}.
\end{equation}
Indeed, the spatial cutoff has gradient norm \(O(R^{-d/2-1})\) and the energy on the interval and its spatial support is bounded by \(CR^d\) using the enclosing backward cylinder. The scalar pairing has an absolutely continuous representative for each cutoff, so the estimate holds at almost every pair of times for the original representative, and simultaneously along \(R=2^j\) after removal of a single null set.

If \(\ell\cdot y\) is spatially homogeneous on almost every slice, \eqref{eq:frozen-affine-increment} makes it constant in time. Thus \(y=y_*+v\xi\) for a fixed \(y_*\in K\) and a unit \(v\in E\) perpendicular to \(\ell\), with \(\xi\) in a compact interval. The \(\ell\)-row becomes
\[
 \Delta G(\xi)=0,\qquad G'(\xi)=\ell^\top A(y_*+v\xi)v .
\]
This follows from the spatial Sobolev chain rule. Testing with products of spatial and temporal test functions, and then taking a countable dense family of spatial tests, gives the elliptic equation for almost every time. Bounded harmonic functions are constant on each whole-space slice. Thus \(\ell^\top A(y)v\,\nabla\xi=0\) almost everywhere. Wherever \(\nabla\xi\neq0\), \(A(y)v\) is a scalar multiple of \(v\) with eigenvalue uniformly positive and bounded on \(K\). Extend this scalar on the zero-gradient set inside the same bounds. The \(v\)-row is then a scalar uniformly parabolic divergence equation for \(\xi\). Scalar oscillation decay on expanding backward cylinders gives constancy, including when the state line lies on a species face.

\subsection{Stationary pressure in the frozen system}

Here \(p_\tau=0\) by itself is sufficient for bounded ancient rigidity, with a priori no spatial constancy of \(p\). Let \((X,p)\) be affine coordinates on the simplex with \(v=\partial_X y\in E\) and \(m\cdot v=0\). The exact additive representation of the three frictions is \(f_{ij}=h_i+h_j\), with real \(h_i\), and its restricted resistance is
\[
 B_y^f|_E=(h\cdot y)I+P_y\operatorname{diag}(h)|_E .
\]
Thus \(D B^f_y[v]v=0\) and \(P_yM_dv=M_dv\). Cramer's rule in \((X,p)\) gives the pressure row of \(A(y)\) as the one-form
\begin{equation}
 m^\top A\,dy=a(y)\,dp+\frac{G(p)}{\det(B_y^f|_E)}\,\dd X .
 \label{eq:frozen-pressure-affine-row}
\end{equation}
Here \(G(p)=\det(B_y^f v,M_dv)\) in the stated coordinates. The numerator is affine in \(y\) and constant in direction \(v\), so depends affinely on \(p\) alone. At a zero of \(G\), \(a\) is a diagonal entry of a triangular normally elliptic matrix, so is positive. If \(G\equiv0\), \(p\) itself solves a uniformly parabolic scalar conservation equation, and Liouville for that equation followed by the affine observation argument gives the result.

Otherwise choose \(L>0\) and \(X'=X-L\Gamma(p)\), with \(\Gamma'=G\). The new \(dp\) coefficient \(a'=a+LG^2/\det B\) is uniformly positive on \(K\) for large enough \(L\), where \(B=B_y^f|_E\). Using a cutoff in a neighborhood of the transformed compact set, extend \(D=(a'\det B)^{-1}\) to a smooth function bounded above and away from zero with bounded first derivatives. Extend \(G\) smoothly with bounded support, keeping it unchanged on a neighborhood of the pressure interval. The vector field of \(dp/\dd X'=-G(p)D(X',p)\) is then bounded with bounded first derivative in \(p\), so its flow is defined for every \(X'\). Label each solution by its \(p\) value at \(X'=0\). Differentiating the flow with respect to the initial \(p\) gives a positive exponential, so the label \(\Psi(X',p)\) has \(\Psi_p>0\), uniformly bounded away from zero on the compact range. Differentiation along the flow gives \(\Psi_{X'}=GD\Psi_p\). On the state space,
\begin{equation}
 p_\tau=\operatorname{div}(\vartheta(y)\nabla\Psi(y)),\qquad
 \vartheta=a'/\Psi_p ,
 \label{eq:frozen-pressure-potential}
\end{equation}
with \(\vartheta\) bounded above and away from zero. If \(p_*\) is a zero of \(G\) in the pressure interval, then \(\Psi(X',p_*)=p_*\); also \(\Psi-p_*\) has the sign of \(p-p_*\) and magnitude at most \(C|p-p_*|\).

When \(p_\tau=0\), \(p\) has a time-independent representative. The same distributional slicing as for the affine observation gives the scalar elliptic equation in \eqref{eq:frozen-pressure-potential} on almost every time slice. Bounded elliptic oscillation decay on expanding balls makes \(\Psi\) spatially constant on each such slice, with value \(q(\tau)\). If \(G\) has a zero \(p_*\) in the pressure interval and \(q(\tau)=p_*\) on an admissible slice, the strict increase of \(\Psi\) in \(p\) gives \(p=p_*\), which is covered by the affine observation. Otherwise, evaluating \(\Psi-p_*\) at a single admissible slice and using \(|\Psi-p_*|\le C|p-p_*|\) puts the essential range of this stationary \(p\) on one side of \(p_*\), a positive distance away. If \(G\) has no zero in the pressure interval, it is uniformly separated from zero there.

Consequently \(\Psi_{X'}\) has a common nonzero sign on the transformed pressure slices attained by the solution. For two admissible times, the two states at almost every \(x\) have the same pressure and are joined by a segment in that pressure slice of the transformed simplex. Hence \(X'(x,\tau_2)-X'(x,\tau_1)\) has a single sign across space and has magnitude at least \(|q(\tau_2)-q(\tau_1)|/\sup|\Psi_{X'}|\), where the supremum is over the transformed compact simplex. This is also the increment in \(X\). Applying \eqref{eq:frozen-affine-increment} to \(X\) shows \(q(\tau_1)=q(\tau_2)\). At a fixed \(p\), the label is strictly monotone in \(X'\), so the whole state is stationary.

The selected level is locally the graph of a smooth \(y=Y(p)\) by the implicit function theorem, since \(\Psi_{X'}\ne0\). These local graphs agree where they overlap, by strict monotonicity in \(X'\) for the extended flow on the pressure interval separated from the root. They cover the essential range of \(p\) by compactness of the state range. This range is an interval: if it had a gap \((p_a,p_b)\) with values on both sides, the spatially Sobolev function \(\min\{(p-p_a)_+,p_b-p_a\}\) would have zero gradient almost everywhere by locality and would take two different values, contrary to connectedness of \(\R^d\). Thus \(Y\) is smooth on a neighborhood of the compact essential range and the spatial Sobolev chain rule applies. In the remaining row the flux is \(\beta(p)\nabla p\), where \(\beta(p)=D X\,A(Y(p))Y'(p)\) and \(D X\) denotes the constant covector of the affine variable \(X\). Along the level, \(m^\top A Y'=0\), while \(m\cdot Y'=1\) and \(A\) is invertible. The two covectors \(D X\) and \(m^\top|_E\) span \(E^*\), so \(\beta(p)\ne0\) on the attained range. A primitive of \(\beta\) composed with \(p\) is bounded and harmonic on \(\R^d\), hence its weak spatial gradient vanishes. By the chain rule and \(\beta(p)\ne0\) this gives \(\nabla p=0\), and the affine observation completes the argument.

This proves the three tests \eqref{eq:frozen-pt-tests} through the ancient argument. To prove the two-dimensional test \eqref{eq:frozen-pt-signed}, use \eqref{eq:frozen-pressure-potential} at each time, with \(\Psi=p\) in the triangular case. Test once with \(\eta_{\log}^2(\Psi-\inf\Psi)\) and once with the opposite sign after subtracting \(\sup\Psi\). Here \(\eta_{\log}\) is the logarithmic cutoff between radii \(r<R_{\log}\), equal to one on the inner ball, with squared gradient integral \(2\pi/\log(R_{\log}/r)\). Zero normal physical flux makes these tests admissible even when the wall intersects the balls. Uniform ellipticity and boundedness of \(\Psi\) give
\[
 \int_{D_r}|\nabla\Psi|^2\,\dd x
 \leq \frac{C}{\log(R_{\log}/r)}
       +C\min\left\{\int_{D_{R_{\log}}}(p_\tau)_+\,\dd x,\
                    \int_{D_{R_{\log}}}(p_\tau)_-\,\dd x\right\}.
\]
Normalize the outer radius to one in a contradictory sequence and
write \(\delta_j\) for its signed defect and \(M_j\) for the selected
gradients. Use \(R_{\log}=1/4\) and \(h_j=\max\{M_j^{-1/2},\delta_j^{1/4}\}\)
as in the full-system proof. On the ball of radius \(Lh_j\), the pressure flux in
\eqref{eq:frozen-pressure-potential} bounds the squared spatial dual
norm of \(p_\tau\) by \(C\int_{D_{Lh_j}}|\nabla\Psi|^2\) using the
mixed tests. Integrating over a time interval of length \(O_L(h_j^2)\)
and rescaling the mixed norm in \(d=2\) gives a bound
\(C_L/\log(1/(4Lh_j))+C h_j^{-2}\delta_j\) on each fixed rescaled
cylinder. Here each source integral on the radius-\(R_{\log}\) subdomain is
bounded by the corresponding outer integral of that sign. This bound
tends to zero.
Since \(h_jM_j\to\infty\), the small-excess lemma supplies a nonconstant
ancient limit, contradicting stationary-pressure rigidity.

At a terminal singular point \((x_0,T)\) the smaller covariance eigenvalue is bounded below on every sufficiently small centered cylinder. For nonconstant masses each pressure-time defect in the theorem also has its own positive lower bound there, including the signed defect in two dimensions. To see this, fix an admissible \(R\) with \(T-R^2\) strictly inside the strong branch. Keep \(D_R\) fixed and shift the time window to \(I_\sigma=(T-\sigma-R^2,T-\sigma)\), \(\sigma>0\). The covariance entries are integrals over \(I_\sigma\) of bounded slice integrands, and the smallest eigenvalue is continuous in these entries. In the time-variance test the minimizing function is \(x\mapsto R^{-2}\int_{I_\sigma}p(x,\tau)\,\dd\tau\). Its squared \(L^2(D_R)\) norm and the integral of \(p^2\) over the window converge to their terminal values as \(\sigma\downarrow0\), by boundedness of \(p\). If a derivative defect is finite on the terminal window, its time integrand is integrable there and on a neighborhood just below \(T-R^2\) by local strong regularity. For the mixed dual norm one can also use the pressure flux bound \(\|p_\tau\|_{-1;D_R}\le C\|\nabla y\|_{L^2(D_R)}\) on that neighborhood. Hence the derivative defects converge by absolute continuity of the integral. An infinite derivative defect already exceeds the threshold. If a terminal defect were strictly below its threshold, the same strict bound would hold for all sufficiently small \(\sigma>0\). The corresponding shifted half-cylinders cover \(D_{R/2}\times(T-R^2/4,T)\) and the theorem gives a common gradient bound there, contradicting singularity.

\section{Frozen interior rigidity and the species-face issue}

\begin{theorem}[Interior rigidity for the frozen equation]\label{thm:interior-rigidity}
Let $h$ be the additive generator of the physical frictions, $f_{ij}=h_i+h_j$, and suppose $m,h,\mathbf1$ are linearly independent. Every bounded ancient local-energy solution of \eqref{eq:frozen-cylinder-operator} on the whole space with bounded spatial gradient and range in a compact subset $K_0$ of the open simplex is constant.
\end{theorem}

\begin{proof}

Use \(B=B_y^f|_E\), \(D=P_yM_d|_E\) and the metric
\(\langle v,w\rangle_{y^{-1}}=\sum_i v_iw_i/y_i\) in the interior. Both \(B,D\) are positive and self-adjoint in this metric. The commutator of their projected diagonal factors is, up to sign, the antisymmetric tensor generated by \(P_yM_d y\) and \(P_y\operatorname{diag}(h)y\). These are independent because \(m,h,\mathbf1\) are independent. Hence \(B,D\) never commute. The diffusion \(B^{-1}D\) has simple positive eigenvalues \(d_1,d_2\), since it is similar to the self-adjoint \(B^{-1/2}DB^{-1/2}\) at each state. If \(v,w\) are corresponding right eigenvectors then \(\langle w,Bv\rangle_{y^{-1}}=0\) and \(\langle w,v\rangle_{y^{-1}}\ne0\): simultaneous orthogonality would diagonalize both \(B,D\) in the same basis. Differentiating \(Dv=d_1 Bv\) along \(v\), and using \(B'[v]v=0\), \(D'[v]v=-(m\cdot v)v\), gives
\[
 D_v d_1=-\frac{(m\cdot v)\langle v,v\rangle_{y^{-1}}}{\langle v,Bv\rangle_{y^{-1}}},
 \qquad
 (d_2-d_1)\langle w,B D_vv\rangle_{y^{-1}}
       =(m\cdot v)\langle w,v\rangle_{y^{-1}} .
\]
Thus the characteristic tangent to \(v\) has nonzero curvature precisely when \(m\cdot v\neq0\).

In the affine coordinates \((x,z)=(m\cdot(y-y_0),h\cdot(y-y_0))\), the matrices have the forms
\[
 D=\begin{pmatrix}p_0-x&q_0\\r_0-z&s_0\end{pmatrix},
 \qquad B=\begin{pmatrix}a_0+z&b_0-x\\c_0&d_0\end{pmatrix}.
\]
The vertical direction is an eigendirection precisely where \(q_0d_0-s_0(b_0-x)=0\). The exceptional set is empty, a single pressure line, or the entire simplex. Whenever it is nonempty the vertical direction belongs to one and the same family along its connected portion in the interior, by simplicity. It follows that the other family (or either family in the empty case) has nonzero pressure component throughout the open simplex. Normalize it as \(v=(1,P(x,z))\). Its graph curvature \((\partial_x+P\partial_z)P\) has strict constant sign there.

Extend the smooth slope field from a neighborhood of \(K_0\) by a smooth compactly supported cutoff. Transport the coordinate \(z\) from one fixed vertical line along this extended field to obtain a first integral \(g\), with \(g_z>0\) by the derivative of the ODE flow. On \(K_0\), \(Dg\,v=0\) and
\[
 D^2g[v,v]=-g_z(\partial_x+P\partial_z)P .
\]
Reverse the sign of \(g\) if needed. The state function \(\mathcal H_{\mathrm{int}}=e^{\beta_g g}\) for large enough \(\beta_g\) then satisfies
\[
 D\mathcal H_{\mathrm{int}}\,A=d_2\,D\mathcal H_{\mathrm{int}},\qquad D^2\mathcal H_{\mathrm{int}}[\xi,A\xi]\geq c_{K_0}|\xi|^2 .
\]
Indeed, \(Dg\) annihilates \(v\) and is a left eigenform for \(d_2\). The \(v,v\) component of \(D^2g\) is uniformly positive. In the production the rank-one Hessian term contributes \(\beta_g^2 \mathcal H_{\mathrm{int}} d_2|Dg\,\xi|^2\), which controls the other direction after Young's inequality absorbs the mixed terms. The local-energy chain rule, obtained by time regularization with the localized test \(D\mathcal H_{\mathrm{int}}(y)\varphi\) for smooth \(\varphi\), gives the scalar entropy balance. Its flux is \(d_2(y)\nabla \mathcal H_{\mathrm{int}}(y)\), and the production is at least \(c_{K_0}|\nabla y|^2\).

Fix a terminal time \(\tau_*\) in the ancient interval and a spatial pole. Translate both to zero. The bounded nonnegative function \(W=\sup_{K_0}\mathcal H_{\mathrm{int}}-\mathcal H_{\mathrm{int}}(y)\) satisfies
\[
 W_\tau-\operatorname{div}(d_2(y)\nabla W)
 \ge c_{K_0}|\nabla y|^2.
\]
The scalar conductivity is bounded and uniformly positive. For a finite past interval \((-T,0)\), compare \(W\) with the zero-initial solution of this scalar equation with a nonnegative compactly supported truncation of \(c_{K_0}|\nabla y|^2\) as the source. Equivalently one may perform the comparison in increasing balls with zero values on their artificial boundary. The latter boundaries carry nonnegative data for \(W\). Passing first to the entire space and then increasing the source to the full production gives
\[
 c_{K_0}\int_{-T}^{0}\!\int_{\R^d}
    \mathcal G(0,0;x,\tau)|\nabla y(x,\tau)|^2\,\dd x\,\dd\tau
 \le W(0,0)\le \operatorname{osc}_{K_0}\mathcal H_{\mathrm{int}} ,
\]
where \(\mathcal G\) is the scalar fundamental solution. The comparison is justified by the scalar maximum principle and the kernel representation on bounded time intervals. The Gaussian lower bound for the drift-free scalar divergence equation \cite{Aronson}, with constants invariant under parabolic scaling, gives \(\mathcal G(0,0;x,\tau)\ge c_M R_j^{-d}\) for \(x\in B_{MR_j}\) and \(-2R_j^2<\tau<-R_j^2\). The slabs are disjoint for \(R_j=2^j\), so, by letting \(T\) grow, one obtains
\[
 \sum_j R_j^{-d}
 \int_{-2R_j^2}^{-R_j^2}\int_{B_{MR_j}}|\nabla y|^2\,\dd x\,\dd\tau<\infty .
\]
The blow-downs \(y_j(x,s)=y(R_jx,R_j^2s)\) have a strong local \(L^2\) limit by the shifted relative entropy and conservation equation. The estimate makes the limiting spatial gradient zero on \((-2,-1)\) and the limiting conservation law makes the state constant there, say \(y_*\). For the shifted relative entropy \(H_\gamma(y_j\mid y_*)\) use the weight \(\phi(x)=\exp(-\sqrt{1+|x|^2})\). The convergence on fixed balls and the uniform bound on the states make its weighted integral tend to zero in time average on \((-2,-1)\), by first cutting off the integrable tail. There are thus times \(s_j\in(-2,-1)\) with \(E_j(s_j)\to0\), where \(E_j(s)=\int\phi\,H_\gamma(y_j(s)\mid y_*)\). The relative flux and coercivity bounds in Lemma~\ref{lem:shifted-entropy} give, after Young's inequality and using \(|\nabla\phi|\le\phi\),
\[
 E_j'(s)+c\int_{\R^d}\phi|\nabla y_j(s)|^2\le C E_j(s).
\]
One can first use compactly supported cutoffs and pass to the weight by the local energy bounds on expanding balls. Gronwall's inequality carries \(E_j\to0\) uniformly to the final unit cylinder. The bounded-gradient solution and its dilations are locally strong by the upgrade in Section~\ref{sec:boundary-estimates}. Lemma~\ref{lem:small-excess} then bounds the rescaled gradient through its upper face. The original gradient at the pole is at most \(C/R_j\). Translations prove constancy.

\end{proof}

When \(m,h,\mathbf1\) are independent, a positive floor for every species on a constant-density frozen cylinder bounds the normalized gradient on its inner half. At a bounded-domain wall take radii below the chart scale. Indeed, an unbounded sequence there gives normalized gradient records with bounded gradients on each fixed backward cylinder and unit gradient at the pole. These remain in a compact interior state range. Interior compactness, together with Neumann reflection at a possible wall, gives a whole-space ancient solution in that same range, contrary to the preceding rigidity argument. If \(0\le y_i\le\delta\) on a cylinder, the unit projection
\((e_i-\mathbf1/3)/|e_i-\mathbf1/3|\) gives a covariance eigenvalue at most \(3\delta^2/8\), since the variance of a variable in $[0,\delta]$ is bounded by its mean squared distance to $\delta/2$. Theorem~\ref{thm:frozen-affine-cylinder} therefore gives the gradient bound on the inner half-cylinder for sufficiently small \(\delta\).

\paragraph{Species-face collision.}
At \(y=(0,1/2,1/2)\), for \(m=(4,1,3)\) and
\((f_{12},f_{13},f_{23})=(19/20,61/20,1)\), the missing-species and
binary tangential diffusivities both equal \(2\). In the tangent basis
\((e_1-e_3,e_2-e_3)\), the diffusion matrix is
\[
 A=\begin{pmatrix}2&0\\-11/20&2\end{pmatrix}.
\]
The repeated eigenvalue has a one-dimensional eigenspace. The scalar
entropy above is constructed on compact interior sets, where the
eigenvalues are simple.

\section{Uniform low-density comparison on an impermeable domain}
\label{sec:singular-comparison}
This section transfers a positive frozen trajectory to a global solution of the physical system. The estimates use second-order Neumann spaces, including on the slow time scale.

\subsection{Gibbs coordinates and the scalar reference}
Fix positive normalized species means \(\bar r=\langle r\rangle\), with \(\langle\cdot\rangle\) denoting the spatial average. Write
\[
 q=P_E(m_i\log r_i)_{i=1}^3,\qquad
 r=R(\rho,q),\qquad
 R_i(\rho,q)=\exp\!\big((q_i+\xi)/m_i\big),\qquad
 \sum_iR_i=\rho .
\]
The last equality determines \(\xi\) uniquely. Put \(G=R_q|_E\) and
\[
 \mathcal M(r)=(B_y^f|_E)^{-1}
       \big(\operatorname{diag}(r)-rr^\top/\rho\big)\big|_E .
\]
At every positive \((\rho,q)\), both \(G(\rho,q)\) and \(\mathcal M(R(\rho,q))\) are symmetric positive definite on \(E\). For \(G\) this follows by differentiation of the total constraint. With \(v_i=R_i/m_i\) one obtains
\[
 G=\big(\operatorname{diag}(v)-vv^\top/\textstyle\sum_i v_i\big)\big|_E,
 \qquad R_\rho=\frac{v}{\sum_i v_i},\qquad
 C_\rho=\frac{\rho}{\sum_i v_i}.
\]
For \(\mathcal M\), the second factor is \(\rho\) times the inverse of the \(y^{-1}\) metric on \(E\), and \(B_y^f|_E\) is positive and self-adjoint in that metric. All matrix bounds are uniform on compact positive state sets. Writing the chemical potential as \(q+\xi\mathbf1\), we have
\[
 d_r=M_d\nabla r-y\nabla C
     =\big(\operatorname{diag}(r)-rr^\top/\rho\big)\nabla q .
\]
Thus the normalized microscopic flux is \(\mathcal M\nabla q\).

For a positive density field \(u\) of the prescribed total mean let \(a_{\mathrm G}[u]\in E\) solve
\[
 \frac1{|\Omega|}\int_\Omega R(u(x),a_{\mathrm G}[u])\,\dd x=\bar r .
\]
The parameter-to-mean map is injective by strict monotonicity and open by positivity of \(\langle G\rangle\). It is relatively proper into the interior of the mean simplex when \(u\) belongs to a positive compact density band. In fact \(r_i\le\sup u\), and \(\langle r_i\rangle\ge c>0\) yields a point where \(r_i\ge c\). Evaluating \(m_i\log r_i-m_j\log r_j=q_i-q_j\) at such points bounds both signs of every difference for constant \(q\). This also gives a compact parameter set for the fixed positive mean, and the image is relatively closed. Thus the parameter exists uniquely. For two fields \(u,v\) on the positive band with the same prescribed mean, write \(q=a_{\rm G}[u]\) and \(q'=a_{\rm G}[v]\). The uniform positive lower bound on \(\langle G\rangle\) along the segment from \(q\) to \(q'\) gives
\[
 c|q-q'|\le
 |\langle R(u,q')-R(v,q')\rangle|
 \le C\|u-v\|_{L^1(\Omega)} .
\]
This proves continuity in \(L^p\) and \(L^2\) on the band. In \(X_\gamma\), the positive profiles form a relatively open subset of the affine space \(\{u:\langle u\rangle=\mathbf1^\top\bar r\}\). On that space the residual \(\langle R(u,q)\rangle-\bar r\) takes values in \(E\), and its \(q\)-derivative is the invertible operator \(\langle G\rangle\) on \(E\). The implicit function theorem gives smooth dependence of \(a_{\mathrm G}[u]\) on that relative chart.

Write \(C(\rho,q)=m\cdot R(\rho,q)\). In slow time \(s=\varepsilon\tau=\delta t_{\rm phys}\), the equations for the normalized densities are
\begin{align}
 \rho_s&=\Div\big(a(\rho,q)\nabla\rho+b(\rho,q)\nabla q\big),\nonumber\\
 \varepsilon q_s
   &=G^{-1}\Div(\mathcal M\nabla q)
       +\varepsilon K(\rho,q)\Delta C
       +\varepsilon\lambda\nabla q\cdot\nabla C,
 \label{eq:comparison-coordinates}\\
 a&=\lambda\rho C_\rho,\qquad b=\lambda\rho C_q,\qquad
 K=\lambda G^{-1}(r-\rho R_\rho).\nonumber
\end{align}
The scalar equation follows by summation. For the compositional equation subtract \(R_\rho\rho_s\) from the vector equation. The first-order part of the Darcy flux then contains \((\nabla r-R_\rho\nabla\rho)\cdot\nabla C=(G\nabla q)\cdot\nabla C\), while the coefficient of \(\Delta C\) is \(\lambda(r-\rho R_\rho)\). In fast time the resulting Darcy terms carry the factor \(\varepsilon\). In the microscopic operator \(G^{-1}\Div(\mathcal M\nabla q)\), second spatial derivatives fall only on \(q\). The boundary conditions become \(\partial_n\rho=\partial_nq=0\).

For the scalar reference, set \(q=a_{\mathrm G}[u(\cdot,s)]\), constant in space at each \(s\). The problem is
\begin{equation}\label{eq:comparison-scalar-reference}
\begin{aligned}
 u_s&=\Div\big(a(u,a_{\mathrm G}[u])\nabla u\big),&
 a(u,a_{\mathrm G}[u])&=\frac{\lambda u^2}{\sum_i R_i(u,a_{\mathrm G}[u])/m_i},\\
 \partial_nu&=0,& u(\cdot,0)&=\rho^0 .
\end{aligned}
\end{equation}
It is globally defined for a positive compatible prescribed profile \(\rho^0\) as in the frozen problem. For local scalar Neumann theory on the full \(X_\gamma\) space, extend the parameter map near a positive compatible profile of the prescribed mean by
\[
 \widetilde a_{\mathrm G}[u]=a_{\mathrm G}\!\left[u+\langle\rho^0-u\rangle\right].
\]
The corrected argument stays positive on a sufficiently small trace neighborhood and has the prescribed mean. Use \(a(u,\widetilde a_{\mathrm G}[u])\) in the divergence equation, which has a positive scalar principal coefficient and the common Neumann domain. Integration conserves \(\langle u\rangle=\langle\rho^0\rangle\), so the local solution satisfies \eqref{eq:comparison-scalar-reference}. The extension applies at any positive restart trace of this mean. Comparison preserves the initial density band. The coefficient in \eqref{eq:comparison-scalar-reference} is bounded above and below there. Scalar local oscillation estimates for the divergence equation extend to the wall by the boundary charts and reflection. The reflected volume density \(\mathcal J\) is positive, time independent and Lipschitz. In these coordinates, dividing by \(\mathcal J\) gives
\[
 u_s=\Div(ag\nabla u)+(ag\nabla\log\mathcal J)\cdot\nabla u,
\]
where \(g\) is the reflected inverse metric. The drift coefficient is bounded. The scalar oscillation estimate gives a uniform parabolic H\"older modulus on backward cylinders in every strip separated from the initial time. A finite spatial covering and the $L^1$ Lipschitz bound for the mean map give a uniform time modulus for $q_*(s)=a_{\mathrm G}[u(s)]$ there.

Apply the buffered gradient-record argument of Section~\ref{sec:boundary-estimates} to this scalar equation. For times $s\ge s_0>0$, fix an outer radius below the chart scale and $\sqrt{s_0}/4$. Unbounded gradients would give buffered records of size $M_j\to\infty$ on such backward cylinders, possibly meeting the wall. At scale $h_j=M_j^{-1}$ the rescaled states $v_j$ have bounded gradients on every fixed backward box and unit gradient at the pole. The principal tensor is $a(v_j,q_*(s_j+h_j^2\sigma))g_j$, with $g_j$ the scaled inverse metric. It has a common modulus of continuity. Since $q_*$ is spatially constant, the lower terms in the expanded equation are bounded by $C(|\nabla v_j|^2+h_j|\nabla v_j|)$, including the Lipschitz metric and volume factors. Coefficient freezing as in Lemma~\ref{lem:wall-record} gives compactness of the gradients through the upper face. On each such box the H\"older estimate makes the oscillation of $v_j$ tend to zero. A constant limit contradicts the unit gradient. This gives a gradient bound near any finite terminal time and uniformly on late strips. On the original scale the same backward estimates now give uniform $W_q^{2,1}$ bounds on shorter strips. The terminal trace and restart argument of Lemma~\ref{lem:wall-record}, with the smooth mean map in the scalar coefficient, gives global continuation and compact late traces in $X_\gamma$. Its conserved mean and the energy identity
\[
 \frac{\dd}{\dd s}\|u-\langle\rho^0\rangle\|_2^2
       =-2\int_\Omega a(u,a_{\mathrm G}[u])|\nabla u|^2\,\dd x
\]
give exponential \(L^2\) convergence to \(\langle\rho^0\rangle\) by the Neumann Poincar\'e inequality. Trace compactness upgrades the convergence to \(X_\gamma\).

\subsection{Neumann estimates and their rescaling}
Fix \(p>d+2\), and use
\[
 X_0=L^p(\Omega),\qquad X_1=W_{\rm N}^{2,p}(\Omega),\qquad
 X_\gamma=W_{\rm N}^{2-2/p,p}(\Omega).
\]
The spaces may be \(E\)-valued. The subscript denotes the single homogeneous normal condition. The time trace space of \(W^{1,p}(I;X_0)\cap L^p(I;X_1)\) is \(X_\gamma\), and \(X_\gamma\) embeds continuously into \(C^{1,\beta}(\overline\Omega)\) for \(0<\beta<1-(d+2)/p\). For fixed positive reference coefficients the scalar and compositional operators are
\[
 L_a=\Div(a_{\rm ref}\nabla\,\cdot\,),\qquad
 L_q=G_{\rm ref}^{-1}\Div(\mathcal M_{\rm ref}\nabla\,\cdot\,).
\]
For \(L_q(\rho,q)=G(\rho,q)^{-1}\Div(\mathcal M(R(\rho,q))\nabla\,\cdot\,)\) and \(L_a(\rho,q)=\Div(a(\rho,q)\nabla\,\cdot\,)\), expansion of the coefficients on a compact positive state range with bounded \(C^1\) norms gives
\[
 \|L_a(\rho,q)-L_a(\rho_*,q_*)\|_{\mathcal L(X_1,X_0)}
 +\|L_q(\rho,q)-L_q(\rho_*,q_*)\|_{\mathcal L(X_1,X_0)}
 \le C\|(\rho-\rho_*,q-q_*)\|_{C^1}.
\]
Indeed each operator difference contains a difference of second-order coefficients times \(D^2v\) and a difference of first-order coefficients times \(Dv\). The reference operators have maximal \(L^p\) regularity on the common domain \(X_1\). Their principal symbols are normally elliptic. The decaying solution for the frozen boundary symbol with tangential frequency \(\xi\) and \(\Re\zeta\ge0\) is
\[
 v(x_d)=\exp\!\left[-x_d\big(|\xi|^2I+\zeta D^{-1}\big)^{1/2}\right]v(0),
\]
where \(D\) is the principal species matrix. For \(|\xi|+|\zeta|>0\) the stable square root is invertible, so a homogeneous Neumann condition yields \(v=0\). Thus the complementing condition of \cite{dhp} holds, with the same Neumann trace as in \cite[Theorem 3.2]{hmpw}. Maximal-regularity constants on bounded time intervals can be chosen uniformly for a compact reference family in \(\mathcal L(X_1,X_0)\) whose members have this property. This follows from the stability of the inverse under small operator perturbations and a finite covering. The constants also apply to shorter intervals by extending the forcing and restricting the solution.

Let \(\Pi_0 v=v-\langle v\rangle\). On the ordinary mean-zero space \(\Pi_0L_q\) is similar to \(L_q\) on the \(G_{\rm ref}\)-weighted mean-zero space. The similarity from ordinary to weighted mean zero is \(v\mapsto v-\langle G_{\rm ref}\rangle^{-1}\langle G_{\rm ref}v\rangle\), with inverse \(\Pi_0\). This map preserves \(X_1\) and \(X_\gamma\), and the operators agree through the similarity because \(L_q\) kills constants and \(\langle G_{\rm ref} L_qv\rangle=0\). Weighted integration by parts makes \(L_q\) self-adjoint and negative in weighted \(L^2\). The weighted Poincar\'e inequality gives a uniform spectral gap on the weighted mean-zero subspace. The scalar \(L_a\) has the corresponding gap on ordinary mean-zero functions. Elliptic regularity and compactness of the resolvent identify the nonzero eigenvalues in the \(L^2\) and \(L^p\) realizations. For the projected compositional operator the bounded projection preserves maximal regularity up to the bounded similarity. Thus on the indicated mean-zero spaces both estimates hold with a small positive exponential time weight and constants independent of interval length \cite[Section 2]{psz}. For a fixed autonomous operator, first take zero initial trace and shift by a weight strictly below the gap. Exponential semigroup decay and Young's convolution inequality bound the half-line $L^p(X_0)$ norm by the forcing norm. Apply bounded-interval maximal regularity with time cutoffs on overlapping unit intervals and sum the $p$th powers. The cutoff errors are bounded by this $L^p(X_0)$ estimate. For a compatible nonzero trace, solve the homogeneous shifted equation on the first unit interval and multiply that solution by a smooth time cutoff equal to one near zero and supported in the interval. Subtract this extension before applying the zero-trace estimate. Bounded-interval maximal regularity controls the strong norm of the extension and the forcing from the cutoff by the initial $X_\gamma$ norm. Uniformity over the compact reference family follows by the same finite covering argument, now with the stable zero-trace inverse on the half-line. This argument also gives uniformity for shifts in a fixed compact interval below the common spectral gap.

\begin{lemma}[Finite-time transfer and global restart]\label{lem:comparison}
Consider a fixed positive frozen reference \(r^f\) on a finite fast-time interval \([0,T_f]\), with prescribed smooth positive compatible density \(\rho^0\). Assume its trace is continuous in \(X_\gamma\) and it belongs to the strong space on that interval. Starting from the same normalized profile, the physical solution exists up to \(T_f\) in fast time for small \(\varepsilon>0\), and its trace differs by \(O(\varepsilon^{1/p})\) from the frozen trace on this interval. In slow time, fix a finite comparison strip and take initial data with species means \(\bar r\) sufficiently close in \(X_\gamma\) to \((\rho^0,a_{\mathrm G}[\rho^0])\). For small \(\varepsilon>0\) the solution exists through that strip and its trace distance to \((u,a_{\mathrm G}[u])\) from \eqref{eq:comparison-scalar-reference} is bounded by a constant times the initial trace error plus \(\varepsilon^{1-1/p}\). On the same fixed-mean set, the global restart radius about the homogeneous state \(\bar r\) can be chosen independently of small \(\varepsilon>0\).
\end{lemma}

\begin{proof}
On a fast-time slab \(I=[\tau_0,\tau_0+h]\) write \(e=\rho-\rho^0\), \(z=q-q^f\), and let \(q^f=P_E(m_i\log r_i^f)_i\). Define
\begin{align*}
 W_f&=\|e\|_{C(I;X_\gamma)}
          +\varepsilon^{1/p}\|e\|_{L^p(I;X_1)}
          +\varepsilon^{1/p-1}\|e_\tau\|_{L^p(I;X_0)},\\
 Z_f&=\|z\|_{C(I;X_\gamma)}
          +\|z\|_{L^p(I;X_1)}+\|z_\tau\|_{L^p(I;X_0)}.
\end{align*}
Freeze each row at \((\rho^0,q^f(\cdot,\tau_0))\). Use \(s=\varepsilon\tau\) in the scalar estimate and \(\tau\) in the compositional estimate. Let \(\eta\) bound the trace difference in the bootstrap, \(\varpi(h)\) be the fixed reference modulus in \(C(X_\gamma)\), and
\(c_I=h^{1/p}+\|q^f\|_{L^p(I;X_1)}\).
Expansion of the two fluxes in \eqref{eq:comparison-coordinates} and the \(C^1\) trace embedding give
\begin{align}
 W_f&\le C\|e(\tau_0)\|_{X_\gamma}
       +C(\varpi(h)+\eta)W_f+C\varepsilon^{1/p}(1+Z_f),\nonumber\\
 Z_f&\le C\|z(\tau_0)\|_{X_\gamma}
       +C(\varpi(h)+\eta)Z_f+C c_I(W_f+Z_f)\nonumber\\[-1mm]
    &\hspace{29mm}+C\varepsilon(1+Z_f)+C\varepsilon^{1-1/p}W_f.
\end{align}
In the first line the term \(b\nabla z\) contributes at most \(C\varepsilon^{1/p}Z_f\) after the scalar time change. Fixed spatial second derivatives of \(\rho^0\) and the \(L^p\)-in-time \(W^{2,p}\) norm of \(q^f\) give the constant term there. In the second line the difference of the microscopic operator applied to \(q^f\) costs \(C c_I(W_f+Z_f)\), since it differentiates the density difference at most once. The chain rule for \(C=m\cdot R(\rho,q)\) gives on the slab
\[
 \|\Delta C\|_{L^p(I;X_0)}
 \le C_{\mathrm{est}}\big(h^{1/p}+\|\rho\|_{L^p(I;W^{2,p})}
                        +\|q\|_{L^p(I;W^{2,p})}\big)
\]
under the \(C^1\) bound. Only the Darcy part contains \(\Delta e\).
Its factor \(\varepsilon\) gives \(C\varepsilon^{1-1/p}W_f\) in the compositional estimate. First choose \(h\) to make \(c_I\) and \(\varpi(h)\) small uniformly over the slab partition, then choose \(\eta\), and finally \(\varepsilon\). The number of slabs is fixed before $\varepsilon$ is chosen. Iterating the absorbed estimate therefore multiplies constants only a fixed number of times and gives the fast-time claim.

On a fixed slow slab \(I=[s_0,s_0+h]\) use \(e=\rho-u\), \(z=\Pi_0q\), and \(a_0=\langle q\rangle\). By conservation,
\(\langle R(\rho,a_0+z)\rangle=\langle R(u,a_{\mathrm G}[u])\rangle=\bar r\).
The derivative of the first mean map with respect to \(a_0\) is \(\langle G(\rho,a_0+z)\rangle\). Its uniform positivity along the segment between \(a_0\) and \(a_{\mathrm G}[u]\), and smooth dependence on \(\rho,z\) on the comparison neighborhood, give
\[
 |a_0-a_{\mathrm G}[u]|
       \le C\big(\|e\|_p+\|z\|_p\big).
\]
In particular the trace norm of the full composition difference is bounded by \(C(\|e\|_{X_\gamma}+\|z\|_{X_\gamma})\). Define
\begin{align*}
 W&=\|e\|_{C(I;X_\gamma)}+\|e\|_{L^p(I;X_1)}
                              +\|e_s\|_{L^p(I;X_0)},\\
 Z&=\|z\|_{C(I;X_\gamma)}
      +\varepsilon^{-1/p}\big(\|z\|_{L^p(I;X_1)}
                                +\varepsilon\|z_s\|_{L^p(I;X_0)}\big).
\end{align*}
Freeze at \((u(\cdot,s_0),a_{\mathrm G}[u(\cdot,s_0)])\) and use the projected compositional operator on the corresponding ordinary mean-zero space. Use $\tau=(s-s_0)/\varepsilon$ on this slab. For a slow-time weight $e^{\sigma_{\mathrm w}(s-s_0)}$ the fast-time weight is $e^{\varepsilon\sigma_{\mathrm w}\tau}$. Choose $0<\sigma_{\mathrm w}$ below the common spectral gap of the frozen reference operators, so the shifts $\varepsilon\sigma_{\mathrm w}$ for $0<\varepsilon\le1$ stay in a fixed compact interval below that gap. The case $\sigma_{\mathrm w}=0$ gives the unweighted estimate used here. Extend the forcing after the time change by zero on the half-line, use the stable maximal-regularity estimate with the prescribed initial trace, and restrict to \(0\le\tau\le h/\varepsilon\). Its constant is independent of the length of the rescaled interval. For the equation \(\varepsilon z_s-\Pi_0 L_q z=f\), the rescaled forcing has norm
\[
 \|f(s_0+\varepsilon\,\cdot)\|_{L^p(0,h/\varepsilon;X_0)}
   =\varepsilon^{-1/p}\|f\|_{L^p(I;X_0)}.
\]
The initial trace is unchanged by the time rescaling. The projection removes the derivative of the spatial mean of \(q\). With \(c_I=h^{1/p}+\|u\|_{L^p(I;X_1)}\) and the slow reference modulus \(\varpi(h)\), the estimates are
\begin{align}
 W&\le C\|e(s_0)\|_{X_\gamma}
        +C(\varpi(h)+\eta)W
        +C c_I(W+Z)+C\varepsilon^{1/p}Z,\nonumber\\
 Z&\le C\|z(s_0)\|_{X_\gamma}
        +C(\varpi(h)+\eta)Z
        +C\varepsilon^{1-1/p}(1+W)
        +C\varepsilon Z .
 \label{eq:comparison-slow-bounds}
\end{align}
The reference flux difference in the scalar row is bounded by \(C c_I(W+Z)\) because it multiplies fixed derivatives of \(u\). The compositional source carries \(\varepsilon\) except for the perturbation of the principal diffusive operator, which is absorbed using the trace closeness to the frozen reference. There is no microscopic reference gradient on this scale because \(a_{\mathrm G}[u]\) is constant in space. Choose \(h,\eta,\varepsilon\) in the same order and iterate the absorbed estimates. Here the reference may start from the prescribed \(\rho^0\) even when the physical solution starts the slow strip with a small scalar error. The fixed species means remain the same.

For either finite comparison the bootstrap starts by fixed-\(\varepsilon\) local theory. On fast slabs \(W_f,Z_f\) control the full strong norm for each fixed \(\varepsilon>0\). On slow slabs the projected norm controls \(z_s\) but omits \((a_0)_s\). Differentiate the conserved species means on the existing strong strip to obtain
\[
 (a_0)_s
   =-\langle G\rangle^{-1}\langle R_\rho\rho_s+Gz_s\rangle .
\]
The average on the right lies in \(E\) because \(\mathbf1^\top R_\rho=1\), \(\mathbf1^\top G=0\), and \(\langle\rho_s\rangle=0\). Here the inverse is taken on \(E\). The uniform positive lower bound for \(\langle G\rangle\) in the comparison neighborhood and the \(L^p\) bounds for \(\rho_s,z_s\) give \((a_0)_s\in L^p(I)\) for each fixed \(\varepsilon>0\). The mean itself is bounded and spatially constant, so \(\rho\) and \(q=a_0+z\) have the full mixed strong norm up to a finite slab endpoint. The strong trace theorem gives a compatible \(X_\gamma\) limit in the compact positive range. The absorbed estimates keep this limit in the comparison neighborhood when the initial trace error and \(\varepsilon\) are sufficiently small. Local theory thus continues the branch at the endpoint, and no additional estimate for \((a_0)_s\) uniform in \(\varepsilon\) is needed.

For the restart take \(u=\langle\rho\rangle\) constant and keep \(z=\Pi_0q\). Freeze the mean-zero operators at the constant Gibbs state determined by \(\bar r\). Both operators now have gaps, so add a small positive exponential weight to the slow norms for all \(s\ge0\). Its corresponding fast-time exponent is multiplied by \(\varepsilon\). The constant reference has zero spatial derivatives, so the terms with reference variation and the independent compositional forcing in \eqref{eq:comparison-slow-bounds} disappear. To estimate the remaining terms, write \(W,Z\) for the exponentially weighted norms and \(\sigma_{\mathrm w}>0\) for the time exponent. In the chain rule for \(\Delta C\), pair the weight on one factor of each quadratic gradient term with the unweighted \(C^1\) bound of the other factor. This gives
\[
 \|e^{\sigma_{\mathrm w} s}\Delta C\|_{L^p(I;X_0)}
 \le C_{\mathrm{est}}\big(\|e^{\sigma_{\mathrm w} s}e\|_{L^p(I;X_1)}
         +\|e^{\sigma_{\mathrm w} s}z\|_{L^p(I;X_1)}\big).
\]
The mean \(a_0\) is spatially constant. The coefficient perturbations are \(O(\eta)\) in the unweighted trace bootstrap. The two cross-row costs are \(C\varepsilon^{1/p}Z\) and \(C\varepsilon^{1-1/p}W\). Absorption gives
\[
 W+Z\le C\big(\|e(0)\|_{X_\gamma}+\|z(0)\|_{X_\gamma}\big)
\]
on every finite interval, with a common constant. A sufficiently small restart radius closes the bootstrap on the entire half-line. This also gives exponential convergence in the strong trace.
\end{proof}

Fix a strictly positive \(r^0\in X_\gamma\) with smooth positive compatible total density \(\rho^0\). First choose a finite slow time \(S\) at which the scalar reference from \eqref{eq:comparison-scalar-reference} is inside the uniform equilibrium restart neighborhood with positive margin. With \(S\) fixed, the slow comparison estimate determines a positive tolerance for the initial trace error at \((\rho^0,a_{\mathrm G}[\rho^0])\). If the frozen trajectory enters this tolerance at a finite fast time \(T_f\), choose \(\varepsilon>0\) small enough for the fast approximation at \(T_f\) and for the error on the slow strip of length \(S\). The homogeneous restart then gives a global classical solution. The small upper bound on the physical multiplier may depend on \(\Omega\), the coefficients, and the fixed initial profile and its species means.

\section{Global classical solutions at low physical density}
\label{sec:low-density}

The preceding transfer needs a positive frozen strong trajectory only up to the finite comparison time. In particular a frozen trajectory whose traces converge along a sequence to the Gibbs state yields a global classical solution for each sufficiently small physical multiplier \(\delta\).

\subsection{Whole-space rigidity and two open coefficient regimes}

\begin{theorem}[All compatible initial shapes]\label{thm:all-shapes}
If every bounded-gradient entire ancient solution of \eqref{eq:frozen} at constant prescribed total density, valued in the closed composition simplex, is constant, the global low-density conclusion holds for every smooth positive compatible initial shape. The hypothesis is satisfied in the following two regimes:
\begin{enumerate}
\item A sufficiently small open neighborhood in the six mass and friction parameters of $f_{ij}^*=a m_i^*m_j^*$, with $a>0$ and arbitrary positive $m_i^*$.
\item A sufficiently small open neighborhood of $m_2^*=m_3^*=m_b$, $f_{12}^*=f_{13}^*=a$, $f_{23}^*=b$, with $m_1^*,m_b,a,b>0$, $m_1^*\ne m_b$, and $a\ne b$. This neighborhood includes pairwise distinct inverse masses.
\end{enumerate}
In the first regime any prescribed triple of positive frictions can be fixed. Indeed, the reference inverse masses are
$m_i^*=\sqrt{f_{ij}^*f_{ik}^*/(a f_{jk}^*)}$, where $\{i,j,k\}=\{1,2,3\}$. At fixed frictions, the identities $g_{ij}=g_{ik}$ become $f_{ij}m_k=f_{ik}m_j$, so nearby masses can make the effective frictions pairwise distinct. Any prescribed positive mass triple can also be fixed. In the second regime the mass and physical friction contrasts can be made arbitrarily large independently.
\end{theorem}

\begin{proof}
At constant normalized total density set $p=m\cdot y$ and $s_m=\sum_i y_i/m_i$. With $f_{ij}=a m_i m_j$ the resistance in the compositional equation is
\[
 B_y^fJ=a\big(pM_dJ-(M_dy)(m\cdot J)\big).
\]
Multiplication by $M_d^{-1}$ and summation in its force law yield
\[
 p_\tau=\Div((s_m/a)\nabla p),\qquad
 J=(ap)^{-1}\big(\nabla y+(s_m y-M_d^{-1}y)\nabla p\big)
\]
in the compositional normalization. Bounded scalar ancient oscillation decay gives constant $p$. If the inverse masses all coincide, $p$ is already constant on the simplex. In all cases, once $p$ is constant, each $y_i$ solves a heat equation with diffusivity $1/(ap)$, including at species faces. For the second reference, its first row is $aJ_1$ and its force is $(m_1-(m_1-m_b)y_1)\nabla y_1$. Thus $y_1$ satisfies a uniformly parabolic scalar equation and is constant in any bounded ancient solution. Then $J_1=0$, and the second row gives scalar heat for $y_2$ with diffusivity $m_b/(a y_1+b(1-y_1))$. The third composition is determined by the total.

Ancient rigidity is open under coefficient changes on the closed simplex. Otherwise take a sequence of nonconstant bounded-gradient ancient solutions with coefficients tending to the reference. In each choose a pole where the gradient magnitude \(M_j\) is at least half its global supremum, translate to that pole and dilate with spatial scale \(M_j^{-1}\) and time scale \(M_j^{-2}\). The past gradients are then bounded by two and the pole gradient has magnitude one. Because the constant-density equation on the whole space is invariant under every such dilation, the argument does not require \(M_j^{-1}\) to tend to zero. The positive reference parameters give uniform normal ellipticity on the closed state range. The interior bounded-gradient compactness argument of Lemma~\ref{lem:wall-record} gives convergence through the top in spatial gradients. The limiting equation is that of the reference, violating its rigidity.

For any fixed compatible prescribed positive profile $\rho^0(x)$, apply the gradient-record argument of Lemma~\ref{lem:wall-record} to \eqref{eq:frozen}. If $\min(1,\sqrt{\tau})\|\nabla y(\tau)\|_\infty$ becomes arbitrarily large, a normalized record produces an entire bounded-gradient ancient solution with constant prescribed density and a unit gradient at its pole. This includes a reflected wall record. The hypothesis excludes it and gives continuation and uniform late Sobolev trace bounds. The integrated entropy production bounds $\int_0^\infty \|M_d\nabla r-y\nabla C\|_2^2\,\dd\tau$. Along a sequence where this tends to zero, compact late traces give a limit in $X_\gamma$. It obeys $\nabla r_i=r_i\nabla C/(m_i\rho^0)$. A zero propagates by this linear relation along paths, which would contradict the positive conserved mean. Therefore the limit has all three components positive. Its projected chemical potentials are spatially constant and its conserved means identify the unique Gibbs profile in Section~\ref{sec:singular-comparison}. The transfer to positive physical density finishes the proof.
\end{proof}

\begin{corollary}[Constant conventional binary diffusivities]\label{cor:constant-d}
The all-shape low-density results apply also to the isotropic ideal model with the convention of \cite{brust} and constant conventional binary diffusivities, on the corresponding open diffusivity regimes. More generally, the preceding proof is valid for a microscopic resistance $\kappa_{\mathrm{mic}}(y)B_y^f$ with any fixed smooth strictly positive scalar multiplier on the closed simplex.
\end{corollary}

\begin{proof}
This multiplies every microscopic flux by $\kappa_{\mathrm{mic}}^{-1}$. At each reference the first scalar equation remains uniformly parabolic with its conductivity multiplied by $\kappa_{\mathrm{mic}}^{-1}$. After that scalar is constant, each remaining component solves a scalar divergence equation, still uniformly elliptic even if its conductivity depends on the whole state. Ancient rigidity and openness follow as before. Entropy, gradient compactness and the singular comparison use smooth positive microscopic mobility on the stated compact ranges and therefore hold with the same changes. For constant prescribed $D_{ij}^{\rm eff}=\mathcal D_{ij}$ in \eqref{eq:physicaldiff}, take $\kappa_{\mathrm{mic}}(y)=(m\cdot y)^{-1}$ and $f_{ij}=\nu m_i m_j/\mathcal D_{ij}$. Equation \eqref{eq:physicaldiff} for the multiplied resistance then gives exactly $\mathcal D_{ij}$. At fixed masses, the smooth invertible change between positive frictions and diffusivities preserves openness. The finite-density binary characteristic proof in Section~\ref{sec:finite-density} is formulated for the constant normalized resistance.
\end{proof}

\subsection{Small mean oscillation of Gibbs composition}

Let \(W(r)=w\) be the compact composition coordinate obtained by normalizing the Gibbs representative to total one:
\[
 w_i=r_i e^{\xi_1/m_i},\qquad
 \sum_i r_i e^{\xi_1/m_i}=1.
\]
Its inverse at total $\rho$ is $r_i=w_i e^{\xi_\rho/m_i}$ with $\sum r_i=\rho$. Denote it by $\widehat R(\rho,w)$, to distinguish it from the chart in chemical potentials. Both are smooth through each species face for $0<L\le\rho\le U$, by the positive scalar derivative determining the multiplier. They preserve the support. The profile $\widehat R(\rho^0,a)$ for a constant composition $a$ has zero microscopic flux and homogeneous Neumann condition.

\begin{theorem}[The Gibbs mean-oscillation class]\label{thm:gibbs-bmo}
Fix positive coefficients and a smooth positive prescribed density \(\rho^0\) with homogeneous Neumann data. For strictly positive compatible initial profiles in \(X_\gamma\) with total \(\rho^0\), sufficiently small \( [w^0]_{\mathrm{BMO}_2(\Omega)} \) in the seminorm of Section~\ref{sec:finite-bmo} gives a global frozen strong trajectory. Along a sequence of times $\tau_j\to\infty$, its traces converge in $X_\gamma$ to the Gibbs profile with the prescribed total and the initial species means. Consequently the physical problem has a global classical solution for every sufficiently small common density multiplier.

The mean-oscillation threshold is independent of compositional derivatives and individual species minima. It is uniform on families of prescribed profiles with fixed bounds \(0<L\le\rho^0\le U\) and a fixed \(C^3\) bound. Their oscillations need not be small.
\end{theorem}

\begin{proof}
The frozen system with $\rho=\rho^0$ has stationary comparison profiles $R^a=\widehat R(\rho^0,a)$ for every $a$ in the closed simplex. Write $e=r-R^a\in E$ and use the relative shifted entropy of Lemma~\ref{lem:shifted-entropy}. It is uniformly equivalent to $|e|^2$ and hence to $|w-a|^2$. The difference of the physical frozen fluxes is
\[
 A_f(r)\nabla e+(A_f(r)-A_f(R^a))\nabla R^a,
 \qquad A_f(r)=L_yP_yM_d .
\]
Since \(R^a\) is stationary with zero microscopic flux, the difference is the entire flux of the frozen equation. In the production, \(\nabla(DH_\gamma(r)-DH_\gamma(R^a))\) equals \(D^2H_\gamma(r)\nabla e\) plus a term bounded by \(C|e|\,|\nabla R^a|\). Lemma~\ref{lem:shifted-entropy} controls the leading quadratic on each column of the tangential gradient \(\nabla e\). The other terms are bounded by \(C(|e|\,|\nabla e|+|e|^2)\) when \(\nabla R^a\) is bounded. On heat scales at most one, the scaled \(\nabla R^a\) has a uniform bound. The differentiated weight term satisfies the same estimate because \(DH_\gamma(r)-DH_\gamma(R^a)\) is \(O(|e|)\) and \(|\nabla\phi|\le C\phi\). Absorb the mixed term by Young's inequality. The no-flux condition removes the wall term, and testing with the weights of Lemma~\ref{lem:neumann-weight} gives
\[
 \frac{\dd}{\dd\tau}\int\phi H_\gamma(r\mid R^a)
        \le C\int\phi H_\gamma(r\mid R^a)
\]
on the heat scale, uniformly in $a$ and any prescribed density in the bounded positive $C^3$ class. The relative inequality has no source independent of $e$.

Bootstrap $\min(1,\sqrt{\tau})\|\nabla w(\tau)\|_\infty\le1$. At a first record before time one, use scale $\sqrt{\tau}$ and take $a$ to be the average of the initial $w$ on that truncated ball. By dyadic averaging the initial weighted squared difference is at most $C[w^0]^2_{\rm BMO_2}$, using the same estimate as in the finite-density mean-oscillation proof. On the final half of the dilated interval the bound for the scaled \(\nabla w\), the smooth inverse Gibbs chart and the prescribed density bound in \(C^1\) give bounds for the scaled \(\nabla r\) and \(\nabla y\). The relative entropy forces the rescaled Gibbs composition to converge to a constant when the seminorm tends to zero. The finite-cylinder conclusion of Lemma~\ref{lem:wall-record} contradicts a unit normalized gradient, using the buffered estimates at fixed or vanishing heat scales. In both cases compactness of the gradients of \(r\) and of the prescribed profiles in \(C^2\) gives compactness of the gradients of \(w\) by the chain rule.

For late records compare to the actual Gibbs representative $a_*$ for the conserved means. The mean map $a\mapsto\avg{\widehat R(\rho^0,a)}$ is continuous and injective on the closed simplex. Its support equals that of $a$, and injectivity on each relative interior follows from the positive Gibbs mean map restricted to those species. Its inverse on its image is uniformly continuous. The same is true jointly on a compact family of positive prescribed densities, by subsequences. Small global $L^2$ oscillation of $w^0$ makes its initial means close to those of the profile with $a=\avg{w^0}$ and hence makes $a_*$ close to $\avg{w^0}$. Since \(a_*\) is positive, \(DH(R^{a_*})\) is the sum of a constant vector and a spatial scalar times \(\mathbf1\). Its pairing with \(r-R^{a_*}\) integrates to zero by conservation of the species means and pointwise equality of the totals. The unshifted relative Gibbs entropy is therefore exactly \(\int(H(r)-H(R^{a_*}))\). Its initial value tends to zero by continuity of \(x\log x\) even when some species minima tend to zero. It is nonincreasing. Uniform strict convexity from below then makes \(\|w-a_*\|_2\) small on the entire surviving strip. Compactness on the final unscaled half-time unit excludes a late unit gradient. The gradient bound and strong restart give a global frozen branch, and the entropy argument of Theorem~\ref{thm:all-shapes} gives a sequence of late traces converging to \(R^{a_*}\) in \(X_\gamma\). The comparison of Section~\ref{sec:singular-comparison} gives the physical conclusion.
\end{proof}

\subsection{A frozen pressure strip}

\begin{theorem}[Low-density existence in a pressure strip]\label{thm:frozen-pressure-strip}
Relabel distinct inverse masses so $m_1<m_2<m_3$. Suppose the effective friction $g_{13}$ is strictly above both of $g_{12},g_{23}$ or strictly below both. For any fixed $\rho_c>0$ there exist $p_*\in(m_1,m_3)$ and $\eta_*,\zeta_*>0$ with the following property. Let \(r^0\in X_\gamma\) be strictly positive on \(\bar\Omega\) with homogeneous Neumann trace and let its prescribed total \(\rho^0\) be smooth. Then the following conditions yield global classical solutions at all sufficiently small common physical density multipliers:
\[
  \|\rho^0-\rho_c\|_{C^3(\bar\Omega)}<\zeta_*,
  \qquad \|m\cdot y^0-p_*\|_\infty<\eta_* .
\]
The two widths are independent of compositional derivatives and individual species minima. The coefficient condition covers four of the six strict effective-friction orders.
\end{theorem}

In an affine coordinate $X$ along pressure levels with
$\partial y/\partial X=v=(m_2-m_3,m_3-m_1,m_1-m_2)$, the $p$-row coefficient of $\nabla X$ is $G(p)/\det(B_y^f|_E)$, where $G$ is affine by \eqref{eq:frozen-pressure-affine-row}. Its signs at $e_1,e_3$ are those of $g_{12}^{-1}-g_{13}^{-1}$ and $g_{13}^{-1}-g_{23}^{-1}$ by the inactive diagonal diffusivities. They are opposite, so $p_*$ is the unique zero of $G$ in the pressure range. On \(p=p_*\) the pressure conductivity is a positive eigenvalue of the triangular compositional matrix. It remains uniformly positive in a fixed thin strip about the line.

\begin{proof}
Fix a bounded \(C^3\) neighborhood of \(\rho_c>0\) with common positive lower and finite upper density bounds. Choose the state strip \(|w|\le\eta\), \(w=p-p_*\), with \(a(y)>0\) uniformly there and with \(\eta^2/|P_E m|^2\) below the covariance threshold of Theorem~\ref{thm:frozen-affine-cylinder}, where \(P_E\) is the Euclidean projection onto \(E\). The gradient-record argument of Lemma~\ref{lem:wall-record}, applied to solutions trapped in the strip with prescribed density in this neighborhood, gives a common bound
\[
 \|\nabla y(\tau)\|_\infty\le C(1+\tau^{-1/2}),\qquad \tau>0,
\]
on every surviving strip. Indeed a normalized record of magnitude one would give a whole-space ancient limit in the strip, including at a reflected wall. For the unit covector \(\ell=P_Em/|P_Em|\), the variance on every centered cylinder of this limit is at most \(\eta^2/|P_Em|^2\), by taking the scalar center \((p_*-\sum_i m_i/3)/|P_Em|\). The affine-variance test of Theorem~\ref{thm:frozen-affine-cylinder} therefore bounds its pole gradient by \(C/R\) at every radius \(R\). This contradicts the unit gradient and gives the stated bound independently of the initial gradient.

Write the exact pressure equation as
\[
 \rho^0 w_\tau=\operatorname{div}(\rho^0 a\nabla w+\rho^0 wV)+S_{\mathrm{strip}},\quad
 V=b(y)\nabla X,\quad S_{\mathrm{strip}}=\operatorname{div}\!\bigl((m\cdot K(y))\nabla\rho^0\bigr),
\]
where $X$ is the affine coordinate complementary to $p$ in the strip. The vector $V$ has zero normal component. Kato's inequality gives the analogous inequality for $|w|$ with source $|S_{\mathrm{strip}}|\le C\zeta(1+|\nabla y|)$, where $\zeta=\|\rho^0-\rho_c\|_{C^3}$. For a candidate early exit at $x_*\in\bar\Omega$ and time $\tau_*\le h_0^2$ set $h=\sqrt{\tau_*}$. The fixed $h_0\le1$ is smaller than a quarter of a uniform tubular radius for $\partial\Omega$. Use the weight of Lemma~\ref{lem:neumann-weight} on the resulting domain. The second integration by parts of the diffusion in the weighted Kato inequality contributes $-\int_{\partial\Omega_h}\rho^0 a|w|\partial_n\phi\le0$. Only the diffusion term is integrated twice. In the rescaled spatial variables, \(\nabla(\rho^0a)\) and the scaled \(V\) are bounded by \(C(1+s^{-1/2})\) for \(0<s\le1\), and the scaled source \(h^2S_{\mathrm{strip}}\) is bounded by \(C\zeta(h^2+h s^{-1/2})\). Gronwall yields terminal weighted mass at most $C(\|w(0)\|_\infty+\zeta)$. At an exit the uniform rescaled gradient bound and the interior ball-volume bound $|\Omega_h\cap B_r(0)|\ge c r^d$, for fixed small $r>0$, give a positive lower bound for this mass. To convert a gradient into pointwise control on $\bar\Omega$, use straight lines within a small flattened half ball and the uniformly bi-Lipschitz chart, or an interior ball when the pole is away from the boundary. This excludes early exits for sufficiently small $\|w(0)\|_\infty+\zeta$.

For later times $h_0^2\le \tau_*\le T$, use the unscaled gradient bound and the global Kato estimate $\int_\Omega\rho^0|w(\tau_*)|\le C_T(\|w(0)\|_\infty+\zeta)$. The same local ball-volume argument excludes exit before $T$. At $\zeta=0$ the weighted global $L^1$ mass is nonincreasing and the estimate works for all $T$ at once. A bounded gradient on a finite smooth strip gives strong restart by the same reflected coefficient estimates and standard Neumann regularity, since the total density is prescribed positive and $A$ is smooth normally elliptic across species faces. Strictly positive data stay positive through finite strips by the scalar inactive-row form and homogeneous Neumann conditions.

\paragraph{Uniform time to small mean oscillation.}
On the trapped branch integrate the frozen entropy identity with zero boundary flux. For \(u=W(r)\), the identity \(W((r_i e^{\xi_{\mathrm{sh}}/m_i})_i)=W(r)\) gives \(DW(r)(r_i/m_i)_i=0\). Since
\[
 \nabla r=M_d^{-1}d_r+(r_i/m_i)_i\,\nabla C/\rho^0,
\]
we have \(\nabla u=DW(r)M_d^{-1}d_r\). The uniform bound on \(DW\) and the resistance dissipation estimate give
\[
  \int_0^T\|\nabla u(\tau)\|_{L^2(\Omega)}^2\,\dd\tau\le C,
\]
independent of \(T\), initial derivatives and species minima on any surviving strip. The late gradient bound, Poincar\'e's inequality on the connected domain, and the boundary volume estimate give a trace close in \(L^\infty\) to a constant composition by a common finite time. More precisely, on a strip surviving to \(T>1\), choose \(\tau_1\in[1,T]\) with \(\|\nabla u(\tau_1)\|_2^2\le C/(T-1)\). Subtract its spatial mean and use the Lipschitz bound and the uniform lower volume of small intersected balls to upgrade the resulting small \(L^2\) norm to small \(L^\infty\) norm.

Use the threshold of Theorem~\ref{thm:gibbs-bmo} for the fixed positive bounded \(C^3\) family. First choose \(T\) large enough that the preceding \(L^\infty\) estimate implies the required BMO smallness. Then choose \(\eta_*,\zeta_*\) to ensure survival in the strip to \(T\). Apply the frozen conclusion of Theorem~\ref{thm:gibbs-bmo} with initial trace \(r(\tau_1)\) and the same prescribed density. Finite-time smoothing gives the Neumann strong trace and the regularity needed for the frozen initial-value problem there. Its conclusion and uniqueness extend the original frozen trajectory globally, with the Gibbs trace subsequence for the original species means. The comparison of Section~\ref{sec:singular-comparison} now applies from the original datum.
\end{proof}

\subsection{Low-density existence near aligned masses}

\begin{theorem}[Low-density existence near aligned inverse masses]\label{thm:aligned}
Let $f_{ij}=h_i+h_j>0$ be pairwise distinct reference frictions and let $m_i=a+b h_i>0$ with $a,b\in\mathbb R$ and $b\ne0$. Fix $p>d+2$ and a smooth strictly positive Neumann-compatible reference profile of constant total density. There is a neighborhood of these coefficients in the six positive mass and friction parameters and of the profile in $W^{3,p}(\Omega;\R^3)$, among smooth positive compatible profiles, such that each datum in this neighborhood gives a global classical solution for all sufficiently small common density multipliers. The neighborhood and the density threshold may depend on the reference shape and the threshold may depend on the perturbed datum.
\end{theorem}

The reference shape may be chosen with any strictly positive composition profile of the stated regularity. Nearby profiles may have nonconstant total density.

The proof derives the aligned root equations and their localized entropy estimates. It then adapts the Euclidean level-set and oscillation argument of \cite[Lemmas~6.5--6.7 and Proposition~6.8]{ternary}.

\begin{proof}
We first analyze the frozen equation at the reference coefficients.

\paragraph{Stieltjes roots and diffusivities.}
For constant prescribed density the composition system in $y$ is $y_\tau=\operatorname{div}(A(y)\nabla y)$ on the simplex. Let $P_y=I-y\mathbf1^\top$, $\bar h=h\cdot y$, $H_f=\sum_i h_i$, and $K=P_y\operatorname{diag}(h)|_E$. Then
\[
 A(y)=(\bar h I+K)^{-1}(aI+bK).
\]
Order the generators $h_1<h_2<h_3$. At interior states the Stieltjes roots $u\in(h_1,h_2)$ and $v\in(h_2,h_3)$ are the two roots in $\xi$ of $\sum_i y_i/(\xi-h_i)=0$. For $k\in E$ the covector $l_\xi(k)=\sum_i k_i/(\xi-h_i)$ is a left eigenvector of $K$ of eigenvalue $\xi$, since $l_\xi(y)=0$ and $\sum_i k_i=0$. The differential $D\xi$ is a nonzero multiple of $l_\xi$ by the implicit function theorem. The affine-moment inverse gives $y_i=c_i(h_i-u)(h_i-v)$, $c_i=\prod_{j\ne i}(h_i-h_j)^{-1}$, and $\bar h=H_f-u-v$. Thus the scalar root diffusivities are
\[
 d_u=\frac{a+bu}{H_f-v},\qquad d_v=\frac{a+bv}{H_f-u}.
\]
The numerators are bounded between $\min_i m_i$ and $\max_i m_i$ by interlacing. The denominators lie between the smallest and largest physical friction, since $H_f-\xi$ ranges between sums of two of the $h_i$. Uniform positivity of both diffusivities does not require positive generators.

\paragraph{Global frozen evolution of the reference shape.}
The polynomial inverse satisfies $y_{uu}=y_{vv}=0$ and $y_{uv}=(y_u-y_v)/(v-u)$. Differentiating the divergence-form equation and resolving onto $D u,D v$ gives
\begin{align*}
 u_\tau&=\operatorname{div}(d_u\nabla u)+\frac{d_u+d_v}{v-u}\nabla u\cdot\nabla v,\\
 v_\tau&=\operatorname{div}(d_v\nabla v)-\frac{d_u+d_v}{v-u}\nabla u\cdot\nabla v.
\end{align*}
The roots satisfy homogeneous Neumann conditions. Fix a slightly larger compact rectangle inside \((h_1,h_2)\times(h_2,h_3)\) containing the initial root ranges in its interior. Until a possible exit, test the \(u\) equation by \((u-k)_+\), where \(k=\max_{\bar\Omega}u(\cdot,0)\). If \(c_*>0\) is a lower bound for \(d_u\), the drift contributes at most
\[
 \frac{c_*}{2}\|\nabla(u-k)_+\|_2^2+
 C\|\nabla v\|_\infty^2\|(u-k)_+\|_2^2 .
\]
The strong trace bounds \(\|\nabla v\|_\infty\) on every compact time interval of the branch. The initial truncation is zero and Gronwall's inequality preserves this bound. The other three truncations give the initial lower bound for \(u\) and both bounds for \(v\). Hence the larger rectangle is never exited. The strictly positive reference data give a fixed separated root rectangle for the entire branch.

The interior oscillation proof of Bothe \cite[Sections 5--6]{ternary} applies with $d_u,d_v$ in place of the two root diffusivities there. For the frame derivative identity and the resulting root equations with variable modal diffusivities see \cite[Lemma 3.5, equation (44), and Theorem 3.7]{bothedarken}. We check the needed bounds below. For a smooth state function $F(u,v)$ the production in the balance with entropy flux $d_u F_u\nabla u+d_v F_v\nabla v$ is
\[
 d_uF_{uu}|\nabla u|^2+d_vF_{vv}|\nabla v|^2+
 (d_u+d_v)\left(F_{uv}-\frac{F_u-F_v}{v-u}\right)\nabla u\cdot\nabla v.
\]
No derivatives of $d_u,d_v$ enter this formula. The pure one-sided truncation $F=\mathbf1_{\{\pm(u-k)>0\}}[(v-k)(u-k)^2/2-(u-k)^3/6]$ has the bracket identically zero, $F_{vv}=0$, production $d_u(v-u)|\nabla(u-k)_\pm|^2$ and entropy comparable to $(u-k)_\pm^2$. For the other root use
\[
 \mathbf1_{\{\pm(v-k)>0\}}\big((k-u)(v-k)^2/2+(v-k)^3/6\big),
\]
whose \(vv\) derivative on the active side is \(v-u\) and whose mixed production also vanishes. Both truncation entropies are \(C^{1,1}\) on the separated rectangle, with zero first derivatives on the switching level, and the stabilized entropies below have the same regularity. The balance follows by smooth approximation of the state functions and the Sobolev chain rule. On the switching level the spatial gradient of the root being truncated vanishes almost everywhere by Sobolev locality, so the jump of the second derivative contributes no production. The usual mixing entropy \(\mathcal H=\sum_i y_i\log y_i\) has zero bracket by \(\sum_i c_i=0\), and is smooth on the fixed interior rectangle. Its Hessian on the simplex is the \(y^{-1}\) metric. In the aligned case \(A=(\bar hI+K)^{-1}(aI+bK)\) is a positive function of the self-adjoint \(K\) in this metric, with the positive eigenvalues \(d_u,d_v\). Thus the mixing-entropy production controls \(|\nabla y|^2\), and invertibility of the root chart on the rectangle gives control of both squared root gradients. Set \(r=(u-k)_\pm\). Choose a constant \(C_0\) with \(\mathcal H+C_0\ge1\), and add \(\epsilon r^2(\mathcal H+C_0)/2\) to the pure truncation. The product contributes the mixing production with weight \(\epsilon r^2/2\) and the nonnegative term \(\epsilon d_u(\mathcal H+C_0)|\nabla r|^2\). All other terms are bounded in absolute value by \(C\epsilon r|\nabla r|^2+C\epsilon r|\nabla r||\nabla v|\). Choosing fixed small $\epsilon$ and applying Young's inequality gives production bounded below by $c(|\nabla r|^2+\epsilon r^2|\nabla v|^2)$, with flux bounded by $C(r|\nabla r|+r^2|\nabla v|)$. All constants are uniform over levels in the separated rectangle. For a time-independent cutoff \(\eta\), the flux term paired with \(2\eta\nabla\eta\) is bounded by half of the localized production plus \(C\int r^2|\nabla\eta|^2\) by Young's inequality. Consequently, at times \(s<t\),
\[
 \int\eta^2 r(t)^2+
 c\int_s^t\int\eta^2\bigl(|\nabla r|^2+r^2|\nabla v|^2\bigr)
 \le C\int\eta^2 r(s)^2+C\int_s^t\int r^2|\nabla\eta|^2.
\]
The same estimate holds with the roots interchanged. With a time-dependent cutoff the calculation adds \(C\int r^2|\eta\partial_\tau\eta|\) to the right side. These are the two-time and cutoff inputs of \cite[Lemma 6.3 and Proposition 6.2]{ternary}.

For the logarithmic test, write \(w\) for either root and \(\widetilde w\) for the other. Its equation has the form
\[
 w_\tau=\operatorname{div}(d_w\nabla w)+b_w\nabla w\cdot\nabla\widetilde w ,
 \qquad |b_w|\le C
\]
on the separated rectangle. For either sign \(\sigma\in\{1,-1\}\), the oriented variable \(\sigma(w-k)\) satisfies this same equation with \(d_w,b_w\) evaluated at the original roots. Put \(r=(\sigma(w-k))_+\) and assume \(0\le r\le H\), with \(H>0\). Following \cite[eqs.~(69)--(71)]{ternary}, fix \(0<\sigma_0<1/2\) and \(0<\delta<\sigma_0/2\), and set
\[
 \ell_\delta(r)=\left[\log\frac{(1-\sigma_0+\delta)H}{H-r+\delta H}\right]_+,
 \qquad F_{\log}(r)=H^2\ell_\delta(r)^2,\qquad
 L_\delta=\log\frac{1-\sigma_0+\delta}{\delta}.
\]
For fixed \(H,\delta\) the function \(F_{\log}\) extends as a \(C^{1,1}\) function by taking zero below \(0\) and its tangent line above \(H\). On the active interval \(r>\sigma_0 H\),
\[
 \frac{(F_{\log}')^2}{F_{\log}''}
   =\frac{2H^2\ell_\delta(r)^2}{1+\ell_\delta(r)}
   \le 2\sigma_0^{-2} L_\delta r^2 ,
\]
with the quotient set to zero on the inactive interval. Test the \(w\) equation by \(\sigma\eta^2F_{\log}'(r)\) for a time-independent spatial cutoff \(\eta\). The time term integrates exactly to the change of \(\int\eta^2F_{\log}(r)\). By the lower bound on \(d_w\), Young's inequality absorbs the mixed drift and cutoff terms into part of the diffusion term and an error bounded by
\[
 C L_\delta\int_s^t\int r^2\bigl(|\nabla\eta|^2+\eta^2|\nabla\widetilde w|^2\bigr).
\]
Apply the two-time stabilized estimate at the same level and sign with a spatial cutoff \(\zeta=1\) on a neighborhood of \(\operatorname{supp}\eta\). This gives
\begin{align*}
 \int\eta^2F_{\log}(r(t))
   +c\int_s^t\int\eta^2F_{\log}''(r)|\nabla r|^2
 &\le \int\eta^2F_{\log}(r(s)) + C L_\delta\int\zeta^2 r(s)^2\\
 &\quad + C L_\delta\int_s^t\int r^2\bigl(|\nabla\eta|^2+|\nabla\zeta|^2\bigr).
\end{align*}
The constants are uniform in the level, \(H,\delta,s,t\) for fixed \(\sigma_0\). This is the estimate of \cite[eq.~(72) and Lemma~6.4]{ternary}, with coefficient one on the initial logarithmic entropy.

The complementary level-set propagation, shrinking and critical-mass steps in \cite[Lemmas 6.5--6.7 and Proposition 6.8]{ternary} use these estimates and Euclidean cylinder geometry on interior cylinders. They therefore give the same fixed oscillation contraction for the present root system. The propagation step uses the logarithmic estimate on a time slice, shrinking uses the relative isoperimetric inequality and the two-time energy, and critical mass uses the quadratic cutoff estimate and parabolic Sobolev embedding. All are applied on interior Euclidean cylinders of the reflected ancient limit. Their coefficient inputs are the uniform diffusivity bounds, root separation, and the mixed-term cancellation established above. Parabolic rescaling preserves the equations and the separated rectangle, so the interior estimate is uniform over all Euclidean radii.

If a bounded-domain branch had unbounded \(\min(1,\sqrt{\tau})\|\nabla y(\tau)\|_\infty\), Lemma~\ref{lem:wall-record} would give a reflected or interior whole-space ancient solution inside the same rectangle, with a unit gradient at the pole. The local \(W^{2,1}_q\) estimates for bounded gradients also apply to that limit, so the truncation and logarithmic tests are admissible on its compact strips. Oscillation contraction on arbitrarily large backward cylinders forces both roots to be constant. This contradicts the unit gradient and gives the late bounds and global frozen continuation.

The frozen physical mixing entropy has finite total dissipation. On the separated root rectangle every species is uniformly positive. Late compactness in $X_\gamma$ along times of vanishing dissipation gives convergence along a subsequence to the Gibbs state with the reference species means. For the constant normalized total this Gibbs state is spatially constant.

\paragraph{Perturbation before taking the small-density limit.}
The homogeneous restart radius of Lemma~\ref{lem:comparison} is uniform
for parameters and positive means in sufficiently small compact
neighborhoods of the reference. The proof uses the uniform
constant-coefficient estimates and mean-zero gaps from
Section~\ref{sec:singular-comparison}. Write that radius in the physical-state
trace coordinates using the uniformly smooth Gibbs chart near the
homogeneous states. Choose a finite fast time when the reference frozen
trace is within one fourth of this radius of its homogeneous mean.
Finite-time strong dependence for the frozen equation preserves the
trajectory up to that time. Its principal coefficients vary smoothly in
the parameters and the prescribed total, and its lower-order terms use
at most two derivatives of the prescribed total. Small compositional trace
perturbations and small \(W^{2,p}\) perturbations of the prescribed total
are therefore sufficient. The \(W^{3,p}\) neighborhood in the statement
gives both bounds, as well as convergence of the conserved means. On a
smaller neighborhood the nearby frozen traces are within one half of
the restart radius of their own homogeneous means. For each nearby
datum the fast comparison puts the physical trace within the restart
radius at the chosen time for sufficiently small density multiplier.
Homogeneous restart gives the asserted global solutions.

\end{proof}

For aligned masses, distinctness of effective frictions fails only on the lines \(a=bh_i\), because
\[
 g_{ij}-g_{ik}=\frac{(h_j-h_k)(a-bh_i)}{m_i m_j m_k}.
\]
The theorem includes these effective-friction coincidences. For each fixed triple of pairwise distinct positive frictions, one extreme value \(a+b h_i\) can still approach zero with the other masses bounded away from zero. This gives arbitrarily large mass contrast in the open joint neighborhoods.

The separated Stieltjes-root rectangle, and hence the neighborhood,
depends on the reference shape.

\section*{Funding}
This research received no specific grant from funding agencies in the public, commercial, or not-for-profit sectors.

\section*{Declaration of competing interest}
The author declares no competing interests.

\section*{Data availability}
No external datasets were used. The numerical illustration in
Figure~\ref{fig:physical-lens} is specified by the characteristic equations
and parameters given in the text and caption.

\section*{AI-assisted preparation}
OpenAI ChatGPT was used for exploring proof ideas, checking arguments and references, and editing the text.
The author takes full responsibility for the content and correctness of the manuscript.


\begin{thebibliography}{99}
\raggedright
\bibitem{Aronson} D.~G.~Aronson, \emph{Non-negative solutions of linear parabolic equations}, Ann.\ Scuola Norm.\ Sup.\ Pisa ser.~3, 22 (1968), no.~4, 607--694. \href{https://www.numdam.org/item/ASNSP_1968_3_22_4_607_0/}{Numdam}.
\bibitem{binning} P.~J.~Binning, D.~Postma, T.~F.~Russell, J.~A.~Wesselingh and P.~F.~Boulin, \emph{Advective and diffusive contributions to reactive gas transport during pyrite oxidation in the unsaturated zone}, Water Resour.\ Res.\ 43 (2007), W02414. \href{https://doi.org/10.1029/2005WR004474}{doi:10.1029/2005WR004474}.
\bibitem{botheadditive} D.~Bothe, \emph{Global strong solutions for Maxwell--Stefan diffusion with additive friction coefficients}, preprint (2026). \href{https://arxiv.org/abs/2609.06732}{arXiv:2609.06732}.
\bibitem{bothedarken} D.~Bothe, \emph{Global strong solutions for the affine Darken class of Maxwell--Stefan diffusion systems}, preprint (2026). \href{https://arxiv.org/abs/2610.00166v1}{arXiv:2610.00166v1}.
\bibitem{ternary} D.~Bothe, \emph{Global strong solutions for the ternary Maxwell--Stefan diffusion system}, preprint (2026). \href{https://arxiv.org/abs/2610.00162v1}{arXiv:2610.00162v1}.
\bibitem{bothedreyer} D.~Bothe and W.~Dreyer, \emph{Continuum thermodynamics of chemically reacting fluid mixtures}, Acta Mech.\ 226 (2015), 1757--1805. \href{https://doi.org/10.1007/s00707-014-1275-1}{doi:10.1007/s00707-014-1275-1}.
\bibitem{bothedruet} D.~Bothe and P.-\'E.~Druet, \emph{Mass transport in multicomponent compressible fluids: Local and global well-posedness in classes of strong solutions for general class-one models}, Nonlinear Anal.\ 210 (2021), 112389. \href{https://doi.org/10.1016/j.na.2021.112389}{doi:10.1016/j.na.2021.112389}.
\bibitem{brust} D.~Brust, K.~Hopf, J.~Fuhrmann, A.~Cheilytko, M.~Wullenkord and C.~Sattler, \emph{Transport of heat and mass for reactive gas mixtures in porous media: Modeling and application}, Chem.\ Eng.\ J.\ 516 (2025), 162027. \href{https://doi.org/10.1016/j.cej.2025.162027}{doi:10.1016/j.cej.2025.162027}. The equation numbers refer to WIAS Preprint 3139 (2024), \href{https://doi.org/10.20347/WIAS.PREPRINT.3139}{doi:10.20347/WIAS.PREPRINT.3139}.
\bibitem{dhp} R.~Denk, M.~Hieber and J.~Pr\"uss, \emph{Optimal $L_p$--$L_q$-estimates for parabolic boundary value problems with inhomogeneous data}, Math.\ Z.\ 257 (2007), 193--224. \href{https://doi.org/10.1007/s00209-007-0120-9}{doi:10.1007/s00209-007-0120-9}.
\bibitem{druetideal} P.-\'E.~Druet, \emph{A theory of generalised solutions for ideal gas mixtures with Maxwell--Stefan diffusion}, Discrete Contin.\ Dyn.\ Syst.\ Ser.\ S 14 (2021), 4035--4067. \href{https://doi.org/10.3934/dcdss.2020458}{doi:10.3934/dcdss.2020458}.
\bibitem{druetjungel} P.-\'E.~Druet and A.~J\"ungel, \emph{Analysis of cross-diffusion systems for fluid mixtures driven by a pressure gradient}, SIAM J.\ Math.\ Anal.\ 52 (2020), 2179--2197. \href{https://doi.org/10.1137/19M1301473}{doi:10.1137/19M1301473}.
\bibitem{foote} R.~L.~Foote, \emph{Regularity of the distance function}, Proc.\ Amer.\ Math.\ Soc.\ 92 (1984), 153--155. \href{https://doi.org/10.1090/S0002-9939-1984-0749908-9}{doi:10.1090/S0002-9939-1984-0749908-9}.
\bibitem{hmpw} M.~Herberg, M.~Meyries, J.~Pr\"uss and M.~Wilke, \emph{Reaction--diffusion systems of Maxwell--Stefan type with reversible mass-action kinetics}, Nonlinear Anal.\ 159 (2017), 264--284. \href{https://doi.org/10.1016/j.na.2016.07.010}{doi:10.1016/j.na.2016.07.010}.
\bibitem{jmz} A.~J\"ungel, J.~Miky\v{s}ka and N.~Zamponi, \emph{Existence analysis of a single-phase flow mixture with van der Waals pressure}, SIAM J.\ Math.\ Anal.\ 50 (2018), 1367--1395. \href{https://doi.org/10.1137/16M1107024}{doi:10.1137/16M1107024}.
\bibitem{jungelstelzer} A.~J\"ungel and I.~V.~Stelzer, \emph{Existence analysis of Maxwell--Stefan systems for multicomponent mixtures}, SIAM J.\ Math.\ Anal.\ 45 (2013), 2421--2440. \href{https://doi.org/10.1137/120898164}{doi:10.1137/120898164}.
\bibitem{mason} E.~A.~Mason and A.~P.~Malinauskas, \emph{Gas Transport in Porous Media: The Dusty-Gas Model}, Chemical Engineering Monographs, vol.~17, Elsevier, Amsterdam, 1983.
\bibitem{ostrowski} L.~Ostrowski and C.~Rohde, \emph{Compressible multicomponent flow in porous media with Maxwell--Stefan diffusion}, Math.\ Methods Appl.\ Sci.\ 43 (2020), 4200--4221. \href{https://doi.org/10.1002/mma.6185}{doi:10.1002/mma.6185}.
\bibitem{psz} J.~Pr\"uss, G.~Simonett and R.~Zacher, \emph{On normal stability for nonlinear parabolic equations}, Discrete Contin.\ Dyn.\ Syst.\ Supplement (2009), 612--621. \href{https://doi.org/10.3934/proc.2009.2009.612}{doi:10.3934/proc.2009.2009.612}.
\bibitem{stockie} J.~M.~Stockie, K.~Promislow and B.~R.~Wetton, \emph{A finite volume method for multicomponent gas transport in a porous fuel cell electrode}, Internat.\ J.\ Numer.\ Methods Fluids 41 (2003), 577--599. \href{https://doi.org/10.1002/fld.453}{doi:10.1002/fld.453}.
\bibitem{weber} A.~Z.~Weber and J.~Newman, \emph{Modeling gas-phase flow in porous media}, Int.\ Commun.\ Heat Mass Transfer 32 (2005), 855--860. \href{https://doi.org/10.1016/j.icheatmasstransfer.2004.08.026}{doi:10.1016/j.icheatmasstransfer.2004.08.026}.
\end{thebibliography}
\end{document}